\documentclass[11pt]{article}
\usepackage{latexsym}
\usepackage{amsfonts,amssymb,amsmath}

\newtheorem{thm}{Theorem}[section]
\newtheorem{cor}[thm]{Corollary}
\newtheorem{lem}[thm]{Lemma}
\newtheorem{mlem}[thm]{Main Lemma}
\newtheorem{pro}[thm]{Proposition}
\newtheorem{defn}[thm]{Definition}

\newcommand{\ov }{\overline }
\title{Polynomial-time right-ideal morphisms and congruences}

\author{ J.C.\ Birget  } 

\date{ \footnotesize{24.iv.2018} }

\begin{document}
\maketitle

\begin{abstract}
We continue with the functional approach to the {\sf P}-versus-{\sf NP}
problem, begun in \cite{s1f, BiInfGen}. We previously constructed a monoid 
${\cal RM}^{\sf P}$ that is non-regular iff ${\sf NP} \neq {\sf P}$.
We now construct homomorphic images of ${\cal RM}^{\sf P}$ with 
interesting properties. In particular, the homomorphic image 
${\cal M}^{\sf P}_{\sf poly}$ of ${\cal RM}^{\sf P}$ is finitely generated
and ${\cal J}^0$-simple, and is non-regular iff ${\sf P} \neq {\sf NP}$.
The group of units of ${\cal M}^{\sf P}_{\sf poly}$ is the famous
Richard Thompson group $V$.
\end{abstract}

\section{ Introduction}

In \cite{s1f} we defined the monoids {\sf fP} and ${\cal RM}^{\sf P}$.
The monoid {\sf fP} consists of the partial functions $A^* \to A^*$ that
are computable by deterministic Turing machines in polynomial time, and
that have polynomial I/O-balance (defined below). In \cite{s1f} it was 
proved that {\sf fP} is finitely generated. 
The submonoid ${\cal RM}^{\sf P}$ consists of the elements of {\sf fP} that 
are right-ideal morphisms of $A^*$ (defined below).
This monoid was studied further in \cite{BiInfGen}, where we proved that 
${\cal RM}^{\sf P}$ is not finitely generated.  We saw that the one-way 
functions (in the sense of worst-case time-complexity) are exactly the 
non-regular element of {\sf fP}, and that $f \in {\cal RM}^{\sf P}$ is 
regular in {\sf fP} iff $f$ is regular in ${\cal RM}^{\sf P}$.  So, one-way 
functions exist iff {\sf fP} is non-regular, iff ${\cal RM}^{\sf P}$ is 
non-regular. It is well-known that one-way functions (according to worst-case
time-complexity) exist iff ${\sf NP} \neq {\sf P}$.
For {\sf P}-vs.-{\sf NP}, see e.g.\ \cite{HU,Papadim,DuKo,HemaOgi,Rothe}; 
for worst-case one-way functions, see e.g.\ \cite{DuKo, HemaOgi}.
For definitions related to semigroups and monoids, see e.g.\ \cite{Grillet}.

In this paper we define some congruences on ${\cal RM}^{\sf P}$ that are
algebraic forms of the {\it padding argument}. The padding argument is
often used in computational complexity in order to decrease the complexity
of a problem by lengthening the inputs (since complexity is measured as a
function of the input-length, lengthening the input reduces complexity).
We use the padding argument in the proof of finite generation of {\sf fP} in
\cite{s1f}; for another use, see Lemma \ref{recTimeConstr}.
These congruences lead to infinitely many quotient monoids (i.e.,
homomorphic images of ${\cal RM}^{\sf P}$), some of which have interesting
and unique properties: \\
$\bullet$ \ ${\cal M}^{\sf P}_{\sf poly}$ is regular iff
${\cal RM}^{\sf P}$ is regular (iff ${\sf NP} = {\sf P}$); moreover,
${\cal M}^{\sf P}_{\sf poly}$ is finitely generated, and its group of units
is the well-known finitely presented infinite simple group $V$ of Richard
Thompson \cite{CFP};  \\
$\bullet$ \ ${\cal M}^{\sf P}_{\sf E3}$ is a homomorphic image of
 ${\cal M}^{\sf P}_{\sf poly}$ and is regular;  \\
$\bullet$ \ ${\cal M}^{\sf P}_{\sf bd}$ is a homomorphic image of
 ${\cal M}^{\sf P}_{\sf E3}$, acts faithfully on the Cantor space
 $A^{\omega}$, and has just two non-zero $\cal D$-classes; \\
$\bullet$ \ ${\cal M}^{\sf P}_{\sf end}$ is a homomorphic image of
 ${\cal M}^{\sf P}_{\sf bd}$, is congruence-simple, and has just one
 non-zero $\cal D$-class. \\
More details about these motivations are given at the end of this
Introduction.

\medskip

We now give some definitions.
A function $f$: $A^* \to A^*$ is  {\em polynomially balanced} iff 
there exists a polynomial $p$ such that for all $x \in {\sf Dom}(f)$:
$\, |f(x)| \leq p(|x|)$ and $|x| \leq p(|f(x)|)$.

Here we always use $A = \{0,1\}$ as our alphabet.
In \cite{s1f}, ${\cal RM}^{\sf P}$ was called ${\cal RM}_2^{\sf P}$, where
the subscript 2 indicated the size of $A$; but since here the size of $A$
will always be 2, we drop the subscript 2.

For an alphabet $A$, the set of all words over $A$ is denoted by 
$A^*$; this includes the empty string $\varepsilon$. By a ``word'' or 
``string'' we will always mean a finite word. The set of all 
non-empty words over $A$ is denoted by $A^+$ ($= A^* - \{\varepsilon\}$).
The length of a word $x \in A^*$ is denoted by $|x|$. For $n \ge 0$ we let 
$A^n = \{x \in A^*: |x| = n\}$, and 
$A^{\le n} = \{x \in A^*: |x| \le n\}$. 

For two strings $v, w \in A^*$, when $v$ is a {\em prefix} of $w$ we write 
$v \le_{\sf pref} w$;  i.e., there exists $x \in A^*$ such that $vx = w$. 
The relation $\le_{\sf pref}$ is a partial order on $A^*$, and is called 
the {\em prefix order}.  We write $v <_{\sf pref} w$ when 
$v \le_{\sf pref} w$ and $v \neq u$ (strict prefix order).
We write $v \parallel_{\sf pref} w$ when $v \le_{\sf pref} w$ or
$w \le_{\sf pref} v$, and then we say that $v$ and $w$ are 
{\em prefix-comparable}. One easily proves that $w \parallel_{\sf pref} v$ 
iff there exist $x_1, x_2 \in A^*$ such that $w x_1 = v x_2$.   
A set $P \subset A^*$ is a {\em prefix code} iff no two elements of $P$ are
prefix-comparable. 
A set $R \subseteq A^*$ is a {\em right ideal} iff $R \, A^* = R$. It is
easy to prove (see e.g.\ \cite{s1f}) that for every right ideal $R$ there
exists a unique prefix code $P$ such that $R = P \, A^*$.
For prefix codes and related concepts, see e.g.\ \cite{Codes}.

For a partial function $f$: $A^* \to A^*$, the domain is
${\sf Dom}(f) = \{ x \in A^* : f(x)$ is defined\}, and the image is
${\sf Im}(f) = f(A^*) = f({\sf Dom}(f))$. When we say ``function'', we mean
partial function. When ${\sf Dom}(f) = A^*$, $f$ is a called a total
function. The restriction of $f$ to a set $S \subseteq A^*$ is denoted by
$f|_S$. The identity map on $A^*$ is denoted by {\bf 1} or ${\bf 1}_{A^*}$, 
and its restriction to $S$ is denoted by ${\bf 1}_S$.

A function $h$: $A^* \to A^*$ is a {\em right-ideal morphism} iff
${\sf Dom}(h)$ is a right ideal, and all $x \in {\sf Dom}(h)$ and all
$w \in A^*$: \ $h(xw) = h(x) \, w$.
In that case, ${\sf Im}(h)$ is also a right ideal.
For a right-ideal morphism $h$, let ${\sf domC}(h)$ (called
the {\em domain code}) be the prefix code that generates ${\sf Dom}(h)$ as a
right ideal.  Similarly, let ${\sf imC}(h)$, called the {\em image code}, be
the prefix code that generates ${\sf Im}(h)$ as a right ideal.
In general, ${\sf imC}(h) \subseteq h({\sf domC}(h))$, and it can happen that
${\sf imC}(h) \neq h({\sf domC}(h))$.

It will often be useful to represent any set $S \subseteq \{0,1,\#\}^*$ as 
a prefix code. We choose one way to do that, as follows. Let 
$P = \{00, 01, 11\}$; this is obviously a prefix code. We define 
$\, {\sf code}(0) = 00, \, {\sf code}(0) = 01, \, {\sf code}(\#) = 11$.
For $w = a_1 \ \ldots \ a_1 \ \ldots \ a_n \in \{0,1,\#\}^*$ we define
${\sf code}(w) =$
${\sf code}(a_1) \ \ldots \ {\sf code}(a_i) \ \ldots \ {\sf code}(a_n)$.
Then for every $L \subseteq \{0,1\}^*$, ${\sf code}(L ) \, 11$ (defined to 
be $\{ {\sf code}(x) \, 11 : x \in L\}$) is a prefix code. 
We also encode any function $f$: $A^* \to A^*$ into a right-ideal morphism
$f^C$: $A^* \to A^*$, defined by 
${\sf domC}(f^C) =$ ${\sf code}({\sf Dom}(f)) \ 11$, and 
$f^C({\sf code}(x) \, 11) = {\sf code}(f(x)) \, 11 \, $ (for all 
$x \in {\sf Dom}(f)$). 
The right-ideal morphisms $f^C$, for $f \in {\sf fP}$, have the following 
important {\em normality} property (see Def.\ \ref{normal}): 
  \  $f^C\big({\sf domC}(f^C)\big) = {\sf imC}(f^C)$. 
This follows immediately from the fact that 
$\, f^C\big({\sf domC}(f^C)\big) \, $ $=$ 
$\, \{ {\sf code}(f(x)) \, 11 \ : \ x \in {\sf Dom}(f)\} \, $ 
is a prefix code. 

We define

\medskip

\hspace{.8in} ${\cal RM}^{\sf P} \ = \ \{ f \in {\sf fP} : f$ is a
                right-ideal morphism of $A^*\}$.

\medskip

\noindent By Prop.\ 2.6 in \cite{s1f}, if $f \in {\cal RM}^{\sf P}$ is
regular in {\sf fP} then $f$ is regular in ${\cal RM}^{\sf P}$.  Hence:
{\em The monoid ${\cal RM}^{\sf P}$ is regular iff $\, {\sf P} = {\sf NP}$.}

\medskip

Since ${\sf P} \neq {\sf NP}$ is equivalent to the non-regularity of
${\cal RM}^{\sf P}$, we are interested in approaches towards proving 
non-regularity or regularity of this monoid. In this paper we study some 
congruences on ${\cal RM}^{\sf P}$; these provide us with infinitely many 
homomorphic images of ${\cal RM}^{\sf P}$. Four of these are of particular 
interest; they form a chain 
$\, {\cal RM}^{\sf P}$ $\twoheadrightarrow$ 
${\cal M}^{\sf P}_{\sf poly}$ $\twoheadrightarrow$ 
${\cal M}^{\sf P}_{\sf E3}$ $\twoheadrightarrow$
${\cal M}^{\sf P}_{\sf bd} \twoheadrightarrow  {\cal M}^{\sf P}_{\sf end}$.
The last three are regular monoids.
Moreover, we find a submonoid ${\cal RM}^{n + o(n)}$ of ${\cal RM}^{\sf P}$
which is non-regular, and which maps homomorphically onto 
${\cal M}^{\sf P}_{\sf poly}$; in addition, ${\cal RM}^{n + o(n)}$ is 
$\, \equiv_{\sf poly}$-equivalent to ${\cal RM}^{\sf P}$ (see the Remark at 
the end of the paper for details). Thus we have the following monoid 
homomorphisms (where $\nearrow$ is injective): 

\bigskip

\hspace{1.63in} ${\cal RM}^{\sf P}$ 

\hspace{1.185in} {\footnotesize 1:1} $\nearrow$ \hspace{.15in} $\downarrow$

\vspace{-.15in}

\hspace{1.6in}  \hspace{.15in} $\downarrow$

\smallskip

\hspace{0.7in} ${\cal RM}^{n + o(n)}$ 
 \ $\twoheadrightarrow$ \ ${\cal M}^{\sf P}_{\sf poly}$
 \ $\twoheadrightarrow$ \ ${\cal M}^{\sf P}_{\sf E3}$ 
 \ $\twoheadrightarrow$ \ ${\cal M}^{\sf P}_{\sf bd}$
 \ $\twoheadrightarrow$ \ ${\cal M}^{\sf P}_{\sf end}$

\bigskip

\noindent where ${\cal M}^{\sf P}_{\sf poly}$ is regular iff 
${\cal RM}^{\sf P}$ is regular (Theorem \ref{regofRMequ}), 
and ${\cal M}^{\sf P}_{\sf E3}$ (hence ${\cal M}^{\sf P}_{\sf bd}$ and
${\cal M}^{\sf P}_{\sf end}$) is regular. On the other hand,
${\cal RM}^{n + o(n)}$ is non-regular (Prop.\ \ref{RMo_nonreg}). The 
triangle of maps starting at ${\cal RM}^{n + o(n)}$ is a commutative 
diagram.  These monoids have other interesting properties: 

\smallskip

\noindent $\bullet$ \ ${\cal M}^{\sf P}_{\sf poly}$ (hence its homomorphic 
images) is finitely generated (Theorem \ref{evalfuncRMpolyequ}); on the 
other hand, ${\cal RM}^{\sf P}$ is not finitely generated (see 
\cite{BiInfGen}).

\smallskip

\noindent $\bullet$ \ ${\cal M}^{\sf P}_{\sf end}$ is congruence-simple and 
has only one non-zero ${\cal D}$-class (Theorems \ref{Mcongr_simple} and 
\ref{M_is_reg}); so ${\cal M}^{\sf P}_{\sf end}$ is the end of the chain.
A priori, it was not obvious that ${\cal RM}^{\sf P}$ should have a 
coarsest non-trivial congruence at all.              

\smallskip

\noindent $\bullet$ \ ${\cal M}^{\sf P}_{\sf bd}$ has exactly two non-zero 
${\cal D}$-classes (Theorem \ref{twoDclasses}), and acts faithfully on 
$A^{\omega}$; in fact, ${\cal M}^{\sf P}_{\sf bd}$ is the monoid of the 
action of ${\cal RM}^{\sf P}$ on $A^{\omega}$ (Prop.\ \ref{end_boundActCOR}).

\smallskip

\noindent $\bullet$ \ The group of units of ${\cal M}^{\sf P}_{\sf poly}$, 
${\cal M}^{\sf P}_{\sf E3}$, and ${\cal M}^{\sf P}_{\sf bd}$, is the famous 
{\em Richard Thompson group} $V$, alias $G_{2,1}$ (Theorem
\ref{boundReg}(2)), whereas the group of units of ${\cal RM}^{\sf P}$ is 
trivial (\cite{s1f} Prop.\ 2.12). 

\smallskip

In the above homomorphism chain, the monoid ${\cal M}^{\sf P}_{\sf poly}$ 
(which is regular iff ${\sf P} = {\sf NP}$) is placed between a monoid that 
is proved to be non-regular, and a monoid that is proved to be regular. 
Whether all this brings us closer to an answer to the {\sf P}-vs.-{\sf NP} 
question remains open.

\section{End-equivalence}

We start out with the most basic congruence on ${\cal RM}^{\sf P}$, which 
turns out to be maximal (i.e., it is not contained in any other congruence,
except the trivial congruence).

\begin{defn} \label{defEndEquiv} \ Two sets $L_1, L_2 \subseteq A^*$ are 
{\em end-equivalent} (denoted by $L_1 \equiv_{\sf end} L_2$)  \ iff 
 \ the right ideals $L_1 A^*$ and $L_2 A^*$ intersect the same right ideals 
of $A^*$  (i.e., for every right ideal $R \subseteq A^*$: 
 \ $R \cap L_1 A^* \neq \varnothing$ \ iff 
 \ $R \cap L_2 A^* \neq \varnothing$). 
\end{defn}
Here we say that two sets $S_1$ and $S_2$ {\em intersect} \ iff
 \ $S_1 \cap S_2 \neq \varnothing$.
Note that for $\equiv_{\sf end}$ it is intersection with $L_i A^*$ that 
matters, not just intersection with $L_i$ (unless $L_i$ is already a right 
ideal). The empty set is only end-equivalent to itself.
In the above definition it is sufficient to use intersections with monogenic 
right ideals (a monogenic right ideal is of the form $wA^*$ with $w \in A^*$):

\begin{lem} \label{defEndEquivMonog}
 \ \ $L_1 \equiv_{\sf end} L_2$ \ iff 
 \ $L_1 A^*$ and $L_2 A^*$ intersect the same {\em monogenic right ideals}. 
\end{lem}
{\bf Proof.}  Let $R$ be any right ideal $R$ that intersects $L_1 A^*$, and 
let $x \in R \cap L_1 A^*$. Then $x A^* \subseteq R$; and $x A^*$ intersects 
$L_1 A^*$, hence (by intersection with the same monogenic right ideals), 
$x A^*$ intersects $L_2 A^*$. So, since $x A^* \subseteq R$, $R$ also 
intersects $L_2 A^*$. In the same way one proves that every right ideal that 
intersects $L_2 A^*$ also intersects $L_1 A^*$. 
 \ \ \ $\Box$

\medskip

Note that if $L_1 \equiv_{\sf end} L_2$ and $L_1 \neq \varnothing$, then
$L_1 A^* \cap L_2 A^* \neq \varnothing$.
Indeed, $L_1 A^*$ has a non-empty intersection with itself, so by 
end-equivalence, $L_1 A^*$ intersects $L_2 A^*$ non-emptily too.
However, if $L_1 \equiv_{\sf end} L_2$ it could happen that 
$L_1 \cap L_2 = \varnothing$; e.g., let $L_1 = \{ 1\}$ and 
$L_2 = \{10, 11\}$.

The next Lemma \ref{equivtoEndEquiv}(1) implies that $\equiv_{\sf end}$ is
definable in the first-order logic of $A^*$ with concatenation.

\begin{lem} \label{equivtoEndEquiv} \ 
Let $L_1, L_2 \subseteq A^*$.

\smallskip

\noindent {\bf (1)} \ \ $L_1 \equiv_{\sf end} L_2$ \ iff

\smallskip

 \ \ \  \ \ \ $(\forall x_1 \in L_1, \, w_1 \in A^*)(\exists x_2 \in L_2)[$
    $x_1 w_1 \parallel_{\sf pref} x_2 \, ]$ \ \ {\rm and} 

\smallskip

 \ \ \  \ \ \ $(\forall x_2 \in L_2, \, w_2 \in A^*) (\exists x_1 \in L_1)[$
    $x_1 \parallel_{\sf pref} x_2 w_2 \, ]$.

\medskip

\noindent {\bf (2)} \ If $L_1 \equiv_{\sf end} L_2$ then
 \ $L_1 \equiv_{\sf end} L_2 \equiv_{\sf end} L_1  A^* \cap L_2 A^*$
$\equiv_{\sf end}$ $L_1 \cup L_2$.
\end{lem} 
{\bf Proof.} (1) If $L_1 A^*$ and $L_2 A^*$ intersect the same right 
ideals then for every $x_1 \in L_1$ and $w_1 \in A^*$, the right ideal 
$x_1 w_1 A^*$ intersects $L_2 A^*$; hence there exists $x_2 \in L_2$ and 
$u_1, u_2 \in A^*$ such that $x_1 w_1 u_1 = x_2 u_2$. Hence, $x_1 w_1$ and $x_2$ are prefix comparable.
Similarly, for every $x_2 \in L_2$ and $w_2 \in A^*$ there exists 
$x_1 \in L_1$ such that $x_1$ and $x_2 w_2$ are prefix comparable.

In the other direction, let us assume the prefix comparability condition,
and let $R$ be a right ideal that intersects $L_1 A^*$. We want to show that
$R$ also intersects $L_2 A^*$. Since $R$ intersects $L_1 A^*$, there exists 
$x_1 w_1 \in R$ such that $x_1 \in L_1$ and $w_1 \in A^*$. 
Then let $x_2 \in L_2$ be prefix comparable with $x_1 w_1$. 
If $x_1 w_1 = x_2 z$ for some $z \in A^*$ then $x_2 z = x_1 w_1 \in R$, so 
$R$ intersects $L_2 A^*$.
If $x_2 = x_1 w_1 z$ for some $z \in A^*$ then $x_2 \in R$ (since 
$x_1 w_1 \in R$, and $R$ is a right ideal), so $R$ intersects $L_2$ (and 
$L_2 A^*$).  

\smallskip

\noindent
(2)  Every right ideal that intersects $L_1 A^* \cap L_2 A^*$ obviously
intersects $L_1 A^*$ and $L_2 A^*$. If a right ideal $R$ intersects 
$L_1 A^*$, let $x_1 w_1 \in L_1 A^* \cap R$. Then by end-equivalence, 
the right ideal $x_1 w_1 A^*$ intersects $L_2 A^*$, i.e., 
$x_1 w_1 z \in L_2 A^*$ for some $z \in A^*$. But $x_1 w_1 z$ also belongs
to $R$, so $R$ intersects $L_1 A^* \cap L_2 A^*$. 
So, $L_1 A^*$ and $L_1 A^* \cap L_2 A^*$ are end-equivalent.

If a right ideal $R$ intersects 
$(L_1 \cup L_2) \, A^* = L_1 A^* \cup L_2 A^*$ then $R$ intersects $L_1 A^*$ 
or $L_2 A^*$. If $R$ intersects $L_1 A^*$, then since 
$L_1 \equiv_{\sf end} L_2$, it intersects $L_2 A^*$ too; so in any case, 
$R$ intersects $L_2 A^*$. Hence, $L_1 \cup L_2 \equiv_{\sf end} L_2$.
 \ \ \ $\Box$

\begin{lem} \label{endequPrefCode}

\noindent {\rm (1)} If $f$ is a right-ideal morphism, and if 
$L_1, L_2 \subseteq {\sf Dom}(f)$ are sets such that 
$L_1 \equiv_{\sf end} L_2$, then $f(L_1) \equiv_{\sf end} f(L_2)$.

\smallskip

\noindent {\rm (2)} For any right-ideal morphism $f$,
$\, f({\sf domC}(f)) \equiv_{\sf end} {\sf imC}(f)$.
\end{lem}
{\bf Proof.} 
(1) Let $R$ be any right ideal that intersects $f(L_1) \, A^*$, and let
$y_1 \in R \cap f(L_1) \, A^*$.  So, since $L_1 \subseteq {\sf Dom}(f)$, we
have $y_1 = f(x_1) \, u_1$ for some $x_1 \in L_1$, $u_1 \in A^*$. And
since $y_1 \in R$ we have also $x_1 u_1 \in f^{-1}(y_1) \subseteq f^{-1}(R)$.
Thus, $x_1 u_1 \in f^{-1}(R) \cap L_1 A^*$. Since $f^{-1}(R)$ is a right
ideal and $L_1 \equiv_{\sf end} L_2$, there exists
$x_2 u_2 \in f^{-1}(R) \cap L_2 A^*$ (with $x_2 \in L_2$, $u_2 \in A^*$).
Hence, $f(x_2) \, u_2 \in f(f^{-1}(R) \cap L_2 A^*)$ $\subseteq$
$R \cap f(L_2) \, A^*$; here we use $f(f^{-1}(R)) \subseteq R$ and
$L_2 \subseteq {\sf Dom}(f)$. So, $R$ intersects $f(L_2) \, A^*$.
Similarly, any right ideal that intersects $f(L_2) \, A^*$ intersects
$f(L_1) \, A^*$.  Thus $f(L_1) \equiv_{\sf end} f(L_2)$.

\noindent
(2) Obviously, $L \equiv_{\sf end} L \, A^*$ for any set $L \subseteq A^*$.
We also have $f({\sf domC}(f)) \ A^* = f({\sf domC}(f) \, A^*) =$
$ {\sf imC}(f) \, A^*$. Hence $f({\sf domC}(f)) \equiv_{\sf end} {\sf imC}(f)$.
 \ \ \ $\Box$

\bigskip

Let $A^{\omega}$ be the set of all ${\omega}$-sequences of elements of $A$;
we call the elements of $A^{\omega}$ {\em ends}. See e.g.\ \cite{PerrinPin}
for the study of infinite words. 
For a set $L \subseteq A^*$, the set $L \, A^{\omega}$ is called the 
{\em ends of} $L$, and denoted by ${\sf ends}(L)$; equivalently, 
${\sf ends}(L)$ is the set of ends that have at least one prefix in $L$, 
so ${\sf ends}(L) = {\sf ends}(L A^*)$. 
The {\em Cantor set topology} on $A^{\omega}$ is described by 
$\{ L \, A^{\omega} : L \subseteq A^*\}$ as set of open sets.
Similarly, $A^*$ is a topological space too, with the set of right ideals as 
set of open sets; we call this the {\em right-ideal topology} of $A^*$.  

\smallskip

{\bf Notation:} For $S \subseteq A^{\omega}$, the {\em closure} is denoted 
by ${\sf cl}(S)$, and the {\em interior} by ${\sf in}(S)$. 

\smallskip

In any topological space, two sets intersect the same open sets iff they have
the same closure. Hence we have the following topological characterization of
end-equivalence.

\begin{pro} \label{topolEndEqu}
 \ For all $L_1, L_2 \subseteq A^*$: \ \ $L_1 \equiv_{\sf end} L_2$ \ iff 
 \ ${\sf cl}({\sf ends}(L_1)) = {\sf cl}({\sf ends}(L_2))$ \ in the Cantor
space $A^{\omega}$, \ iff \ ${\sf cl}(L_1) = {\sf cl}(L_2)$ \ in the 
right-ideal topology of $A^*$.  
 \ \ \  \ \ \ $\Box$
\end{pro}
The following example illustrates the importance of closure in the above
Proposition. Let $L_1 = 0^*1$, and $L_2 = \{\varepsilon\}$. Then
$0^*1 \equiv_{\sf end} \{\varepsilon\}$, and 
${\sf cl}(0^*1 \, \{0,1\}^{\omega})$ $=$ $\{0,1\}^{\omega}$ $=$
${\sf cl}(\{\varepsilon\} \, \{0,1\}^{\omega})$. 
But ${\sf ends}(0^*1) \neq {\sf ends}(\{\varepsilon\})$, since 
${\sf ends}(\{\varepsilon\}) = \{0,1\}^{\omega}$, while
${\sf ends}(0^*1) = 0^*1 \, \{0,1\}^{\omega}$ $=$
$\{0,1\}^{\omega} - \{0^{\omega}\}$.

\medskip

The equivalence relation $\equiv_{\sf end}$ can be generalized
to a pre-order: For $L_1, L_2 \subseteq A^*$ we define
 \ $L_1 \subseteq_{\sf end} L_2$ \ iff
 \ every right ideal that intersects $L_1 A^*$ also intersects $L_2 A^*$.
Obviously, $L_1 \equiv_{\sf end}  L_2$ \ iff
$L_1 \subseteq_{\sf end} L_2$ and $L_2 \subseteq_{\sf end} L_2$.
And we have: 
 \ \  $L_1 \subseteq_{\sf end} L_2$ \ iff
 \ ${\sf cl}({\sf ends}(L_1)) \subseteq {\sf cl}({\sf ends}(L_2))$.

\medskip

The following Lemma shows that a one-point change in ${\sf ends}(P)$ does
not change end-equivalence.

\begin{lem} \label{puncturePrefCode}
For any prefix code $P \subset A^*$ and any end $v \in {\sf ends}(P)$
there exists a prefix code, called $P(-v)$, such that 
$\, {\sf ends}(P(-v)) = {\sf ends}(P) - \{v\}$ \ and 
$\, P(-v) \equiv_{\sf end} P$.  
\end{lem}
{\bf Proof.} Let $v = v_1 \ldots v_i \ldots \ \ \in {\sf ends}(P)$, where 
$v_i \in A$ for all $i \ge 1$. Since $v \in {\sf ends}(P)$, there exists 
$i_0 \ge 0$ such that $v_1 \ldots v_{i_0} \in P$. We now define 

\smallskip

 \ \ \ \ \ $P(-v) \ = \ (P - \{v_1 \ldots v_{i_0}\})$ \ $\cup$ 
 \ $\{v_1 \ldots v_{i_0} \ldots v_j \, \ov{v}_{j+1} \, : \, j \ge i_0\}$.

\smallskip

\noindent The set 
$\{v_1 \ldots v_{i_0} \ldots v_j \, \ov{v}_{j+1} \, : \, j \ge i_0\}$ is the
{\em border of the end} $v$ from $v_1 \ldots v_{i_0}$ onwards; indeed,   
$v_1 \ldots v_{i_0} \ldots v_j \, \ov{v}_{j+1}$ is the sibling node of
$v_1 \ldots v_{i_0} \ldots v_j \, v_{j+1}$, and the latter is a prefix of 
$v$.  Clearly, $P(-v)$ is a prefix code and 
$\, {\sf ends}(P(-v)) = {\sf ends}(P) - \{v\}$.

To show that $P(-v) \equiv_{\sf end} P$, it is obvious that every right 
ideal that intersects $P(-v) \, A^*$ also intersects $P A^*$ (since
${\sf ends}(P(-v)) \subseteq {\sf ends}(P)$).
Conversely, let $R \subseteq A^*$ be any right ideal that intersects 
$P A^*$ (at say $u \in R \cap PA^*$); we want to show that $R$ also 
intersects $P(-v) \, A^*$.
If $u \in (P - \{v_1 \ldots v_{i_0}\}) \ A^*$ then $u \in P(-v) \, A^*$.
Alternatively, $v_1 \ldots v_{i_0}$ is a prefix of $u$, i.e., 
$u = v_1 \ldots v_{i_0} \, z$ for some $z \in A^*$. Let 
$v_1 \ldots v_{i_0} \ldots v_j$ be the longest prefix of $u$ that is also a
prefix of the end $v$. If $u = v_1 \ldots v_{i_0} \ldots v_j$, then
$u \, \ov{v}_{j+1} = v_1 \ldots v_{i_0} \ldots v_j \ov{v}_{j+1} \in$
$R \, \cap \, P(-v) \, A^*$; if, on the other hand, $|u| > j$ then $u$ is 
of the form $u = v_1 \ldots v_{i_0} \ldots v_j \ov{v}_{j+1} \ldots \ $ (the 
only alternative would be that 
$u = v_1 \ldots v_{i_0} \ldots v_j v _{j+1} \ldots \ $, 
but that would contradict the maximality of $j$). 
Thus, $u = v_1 \ldots v_{i_0} \ldots v_j \ov{v}_{j+1} \ldots \ $;
this belongs to $P(-v) \, A^*$, since 
$v_1 \ldots v_{i_0} \ldots v_j \ov{v}_{j+1} \in P(-v)$.
 \ \ \ $\Box$

\begin{pro} \label{topolInterior}
 \ For any prefix code $P \subset A^*$, we have:

\smallskip

 \ \ \ $\bigcup_{Q \equiv_{\sf end} P} \, {\sf ends}(Q)$ $ \, = \, $
  ${\sf in}({\sf cl}({\sf ends}(P)))$, \ \ \ and \ \ \ \      
  $\bigcup_{Q \equiv_{\sf end} P} \, Q A^* \ \equiv_{\sf end} \ P$.

\medskip

\noindent Moreover,

\smallskip

 \ \ \ $\bigcap_{Q \equiv_{\sf end} P} \, {\sf ends}(Q) \ = \ \varnothing$ 
 \ $=$ \ $\bigcap_{Q \equiv_{\sf end} P} \, Q A^*$.
\end{pro}
{\bf Proof.} Concerning the union of sets of ends:  

\noindent $[\subseteq]$: 
If $Q \equiv_{\sf end} P$, then (by Prop.\ \ref{topolEndEqu}),
${\sf ends}(Q) \subseteq {\sf cl}({\sf ends}(P))$. And ${\sf ends}(Q)$ is 
open, hence ${\sf ends}(Q) \subseteq {\sf in}({\sf cl}({\sf ends}(P)))$. 

\noindent $[\supseteq]$:  If $x \in A^*$ is such that
${\sf ends}(\{x\}) \subseteq {\sf in}({\sf cl}({\sf ends}(P)))$, 
then ${\sf ends}(\{x\}) \subseteq {\sf cl}({\sf ends}(P))$, hence
${\sf cl}({\sf ends}(\{x\}) \subseteq {\sf cl}({\sf ends}(P))$, hence
${\sf cl}({\sf ends}(\{x\} \cup P)) = {\sf cl}({\sf ends}(P))$. Therefore,
$\{x\} \cup P \equiv_{\sf end} P$ (by Prop.\ \ref{topolEndEqu}).  Thus, 
${\sf ends}(\{x\}) \subseteq \bigcup_{Q \equiv_{\sf end} P} {\sf ends}(Q)$
for every ${\sf ends}(\{x\}) \subseteq {\sf in}({\sf cl}({\sf ends}(P)))$;
thus $\, {\sf in}({\sf cl}({\sf ends}(P))) \subseteq $
$\bigcup_{Q \equiv_{\sf end} P} \, {\sf ends}(Q)$.

Concerning the union of the $QA^*$:  Let $R$ be any right ideal that
intersects $\bigcup_{Q \equiv_{\sf end} P} Q A^*$. Then (by the definition
of union) there exists $Q \equiv_{\sf end} P$ such that $R$ intersects
$Q A^*$.  Since $Q \equiv_{\sf end} P$, every right ideal intersecting
$Q A^*$ intersects $P A^*$, so $R$ intersects $P A^*$.
Conversely, it is obvious that every right ideal that intersects $P$ also
intersects $\bigcup_{Q \equiv_{\sf end} P} Q A^*$. 
Thus, $\bigcup_{Q \equiv_{\sf end} P} \, Q A^* \ \equiv_{\sf end} \ P$.

Concerning the intersection of the $QA^*$, we observe that for any $n > 0$,  
$P \, A^n \equiv_{\sf end} P$. Moreover, 
$\bigcap_{n>0} P \, A^n \, A^*$ $=$ $\varnothing$ (since for any length $n$,
this intersection contains no word of length $< n$).
Hence \ $\bigcap_{Q \equiv_{\sf end} P} \, Q A^*$ \ (which is a subset of
$\bigcap_{n>0} P \, A^n \, A^*$) is $\varnothing$.

Concerning the intersection of the ${\sf ends}(Q)$:
By Lemma \ref{puncturePrefCode} we have 
$\bigcap_{v \in {\sf ends}(P)} {\sf ends}(P(-v)) \ = \ \varnothing$. 
The result follows, since $P(-v) \equiv_{\sf end} P$.
 \ \ \ $\Box$ 

\bigskip

\noindent End-equivalence can also be defined for right-ideal morphisms:

\begin{defn} \label{endequDefMorph} \ 
Two right-ideal morphisms $f_1,f_2$ are {\em end-equivalent} (denoted by 
$ \, f_1 \equiv_{\sf end} f_2$) \ iff 
 \ ${\sf Dom}(f_1) \equiv_{\sf end} {\sf Dom}(f_2)$, and the 
restrictions of $f_1$ and $f_2$ to ${\sf Dom}(f_1) \cap {\sf Dom}(f_2)$ 
are equal.
\end{defn}
{\bf Notation:} 
By $[f]_{\sf end}$ we denote the $\, \equiv_{\sf end}$-equivalence class of 
$f$ in ${\cal RM}^{\sf P}$. So for $f \in {\cal RM}^{\sf P}$ we have 
 \ $[f]_{\sf end}$ \ $=$ 
 \ $\{g \in {\cal RM}^{\sf P}: g \equiv_{\sf end} f\}$ \ $\subset$
 \ ${\cal RM}^{\sf P}$. Note that although $\equiv_{\sf end}$ is defined 
for all right-ideal morphisms, we define $[f]_{\sf end}$ to contain only 
elements of ${\cal RM}^{\sf P}$.

\begin{pro} \label{endequivUnionIntersLatt}
  \ (1) Let $P_1, P_2 \subset A^*$ be prefix codes such that
$P_1 \equiv_{\sf end} P_2$, let $P_{\cap}$ be the prefix code that
generates the right ideal $P_1 A^* \cap P_2 A^*$, and let $P_{\cup}$ be the
prefix code that generates the right ideal $P_1 A^* \cup P_2 A^*$.  Then
$P_1 \equiv_{\sf end} P_2$
$\equiv_{\sf end} P_{\cap} \ \equiv_{\sf end} P_{\cup}$.

\smallskip

\noindent
(2) Let $f_1, f_2$ be right-ideal morphisms such that
$f_1 \equiv_{\sf end} f_2$. Then $f_1 \cap f_2$ and $f_1 \cup f_2$ are
right-ideal morphisms, and $f_1 \equiv_{\sf end} f_2$ \   
$\equiv_{\sf end}$ \ $f_1 \cap f_2 \ \equiv_{\sf end} \ f_1 \cup f_2$.

\smallskip

\noindent
(3) If $f_1, f_2 \in {\cal RM}^{\sf P}$ and $f_1 \equiv_{\sf end} f_2$, then 
$f_1 \cap f_2, \, f_1 \cup f_2 \, \in {\cal RM}^{\sf P}$.
\end{pro}
{\bf Proof.}
(1) This follows from Lemma \ref{equivtoEndEquiv}(2), and the fact that
$L A^* \equiv_{\sf end} L$ for all sets $L \subseteq A^*$. 

\noindent
(2) By (1) we have ${\sf Dom}(f_1) \equiv_{\sf end} {\sf Dom}(f_2)$
$\equiv_{\sf end}$ ${\sf Dom}(f_1) \cap {\sf Dom}(f_2)$ $\equiv_{\sf end}$
${\sf Dom}(f_1) \cup {\sf Dom}(f_2)$.  
Also, ${\sf Dom}(f_1 \cap f_2)$  $=$ ${\sf Dom}(f_1) \cap {\sf Dom}(f_2)$.
And since $f_1 = f_2$ on ${\sf Dom}(f_1 \cap f_2)$ we have
$f_1 = f_1 \cap f_2$ on ${\sf Dom}(f_1 \cap f_2)$. Hence
$f_1 \equiv_{\sf end} f_1 \cap f_2$, and similarly for $f_2$.
Also, ${\sf Dom}(f_1 \cup f_2)$  $=$ ${\sf Dom}(f_1) \cup {\sf Dom}(f_2)$.
And $f_1 = f_1 \cup f_2$ on ${\sf Dom}(f_1)$, and
$f_2 = f_1 \cup f_2$ on ${\sf Dom}(f_2)$, hence
$f_1 \equiv_{\sf end} f_1 \cup f_2 \equiv_{\sf end} f_2$.
 
\noindent
(3) Since the complexity class {\sf P} is closed under $\cap$ and 
$\cup$, ${\sf Dom}(f_1 \cap f_2)$ and ${\sf Dom}(f_1 \cup f_2)$ belong to 
{\sf P}.  Since $f_1 \cap f_2$ is the restriction of $f_1$ to 
${\sf Dom}(f_1 \cap f_2)$, we have $f_1 \cap f_2 \in {\cal RM}_2^{\sf P}$.
To compute $(f_1 \cup f_2)(x)$ in polynomial time, check whether 
$x \in {\sf Dom}(f_1)$, and if so, compute $f_1(x)$; otherwise, check
whether $x \in {\sf Dom}(f_2)$, and compute $f_2(x)$.  
Polynomial balance of $f_1 \cap f_2$ and $f_1 \cup f_2$ follows from the
polynomial balance of $f_1$ and $f_2$.
 \ \ \ $\Box$

\begin{cor} \label{equivclassLattice}
 \ Every $\, \equiv_{\sf end}$-class (within ${\cal RM}^{\sf P}$ or within the
monoid of all right-ideal morphisms) is a lattice under $\subseteq$, 
$\cup$ and $\cap$. In particular, $[f]_{\sf end}$ is a lattice, for every 
$f \in {\cal RM}^{\sf P}$.
\end{cor}
{\bf Proof.} The lattice property follows from Prop.\ 
\ref{endequivUnionIntersLatt}.
  \ \ \ $\Box$

\begin{pro} {\bf (preservation of injectiveness under
$\equiv_{\sf end}$).}  \label{injEndequiv}

If $f, g$ are right-ideal morphisms such that $g \equiv_{\sf end} f$, and
if $f$ is injective, then $g$ is injective.
\end{pro}
{\bf Proof.} From $g \equiv_{\sf end} f$ it follows that $f$ and $g$ agree
on $D = {\sf Dom}(g) \cap {\sf Dom}(f)$. So, $h = g|_D = f|_D$ is injective.
Also, $h \equiv_{\sf end} g$.
To show that $g$ is injective, let $x_1, x_2 \in {\sf Dom}(g)$ be such that
$g(x_1) = g(x_2)$. From $D = {\sf Dom}(h) \equiv_{\sf end} {\sf Dom}(g)$ it 
follows that $x_1 A^*$ intersects ${\sf Dom}(h)$ at $x_1 u$ (for some 
$u \in A^*$).  Hence, since $g$ is a right-ideal morphism, 
$g(x_1 u) = g(x_2 u)$, where $x_1 u \in {\sf Dom}(h)$ and 
$x_2 u \in {\sf Dom}(g)$. 
Again, since ${\sf Dom}(h) \equiv_{\sf end} {\sf Dom}(g)$ it follows that 
$x_2 u A^*$ intersects ${\sf Dom}(h)$ at $x_2 u v$ (for some $v \in A^*$).
Then $g(x_1 uv) = g(x_2 uv)$, where 
$x_1 uv \in {\sf Dom}(h)$ and $x_2 uv {\sf Dom}(h)$.
Since $h(x_1 uv) = g(x_1 uv) = g(x_2 uv) = h(x_2 uv)$, injectiveness of $h$
implies $x_1 uv = x_2 uv$; hence $x_1 = x_2$, 
so $g$ is injective.
 \ \ \ $\Box$

\begin{defn} {\bf (maximum extension).}  \label{MaxExtMorph} 
 \ For any right-ideal morphism $f$: $A^* \to A^*$ we define

\smallskip

 \ \ \ $f_{\sf e, max} \ = \ \bigcup \, \{ g : g$ is a right-ideal morphism 
       with $g \equiv_{\sf end} f \}$.
\end{defn}
It can happen that $f_{\sf e, max} \not\in {\cal RM}^{\sf P}$
when $f \in {\cal RM}^{\sf P}$ (and that is in fact the ``usual'' case 
 -- see Prop.\ \ref{bMaxRec}).

\begin{pro} \label{MaxExtPROP} \hspace{-.14in} .

\noindent (1) For every right-ideal morphism $f$, $f_{\sf e, max}$ is a 
function, and a right-ideal morphism $A^* \to A^*$.  
It is the maximum extension of $f$ among all right-ideal morphisms that are
$\, \equiv_{\sf end} f$. So, $ f_{\sf e, max} \equiv_{\sf end} f$, and
$f_{\sf e, max}$ is the unique right-ideal morphism that is maximal (under 
$\subseteq$) in the set {\rm $\{g : g$ is a right-ideal morphism, 
and $f \equiv_{\sf end} g \}$}.

\smallskip

\noindent (2) For any right-ideal morphisms $h,k: A^* \to A^*$: 

\smallskip

 \ \ \  \ \ \ $h \equiv_{\sf end} k$ \ \ iff 
 \ \ $h_{\sf e, max} = k_{\sf e, max}$.
\end{pro}
{\bf Proof.} (1) If $f_{\sf e, max}$ were not a function there would exist
right-ideal morphisms $g_1, g_2$ such that
$g_1 \equiv_{\sf end} g_2 \equiv_{\sf end} f$, 
$f \subseteq g_1$, $f \subseteq g_2$, and for some
$x \in {\sf Dom}(g_1) \cap {\sf Dom}(g_2)$: $g_1(x) \neq g_2(x)$.
But $f \subseteq g_1 \cup g_2$, and by Prop.\ \ref{endequivUnionIntersLatt},
$f \equiv_{\sf end} g_i \equiv_{\sf end} g_1 \cup g_2$, and $g_1 \cup g_2$ 
is a function. Hence $g_1(x) = g_2(x)$ for all 
$x \in {\sf Dom}(g_1) \cap {\sf Dom}(g_2)$.
So, $f_{\sf e, max}(x)$ has at most one value for all $x$.
The facts that $f_{\sf e, max}$ is a right-ideal morphism, and that it is
maximum, are straightforward.

By Prop.\ \ref{topolInterior}, \ ${\sf Dom}(f)$ \ $\equiv_{\sf end}$
 \ $\bigcup_{g \equiv_{\sf end} f} {\sf Dom}(g)$. 
Since $f_{\sf e, max}$ agrees with $f$ on ${\sf Dom}(f)$ we conclude that 
$f \equiv_{\sf end} f_{\sf e, max}$.

\smallskip

\noindent (2) This follows immediately from 
$f \equiv_{\sf end} f_{\sf e, max}$. 
  \ \ \ $\Box$

\bigskip

A right-ideal morphism $f$: $A^* \to A^*$ can be extended to the partial 
function $f$: $A^{\omega} \to A^{\omega}$, defined as follows:

\hspace{0.9in}  $f(pw) = f(p) \, w$, \ for any $p \in {\sf domC}(f)$ and 
                $w \in A^{\omega}$.

\medskip 

\noindent Then ${\sf Dom}(f) \cap A^{\omega}$ $=$
${\sf domC}(f) \, A^{\omega}$ $=$ ${\sf Dom}(f) \, A^{\omega}$;  and
${\sf Im}(f) \cap A^{\omega} = {\sf imC}(f) \, A^{\omega}$ $=$
${\sf Im}(f) \, A^{\omega}$.
We use the same name $f$ for the extended function, and its restrictions
to $A^*$ or to $A^{\omega}$; the context will always make it clear which
function is being used.

Any right-ideal morphism $f$: $A^* \to A^*$ is continuous (with respect to 
the topology defined on $A^*$ by the right ideals); indeed, for every right 
ideal $R \subseteq A^*$, $f^{-1}(R)$ is a right ideal. Similarly, the 
extension of $f$ to $A^{\omega}$ is a continuous function in the Cantor 
space topology.

\begin{lem} \label{MaxCompNEQ} 
 \ There exist right-ideal morphisms $f, g: A^* \to A^*$ such that 

\medskip

 \ \ \  \ \ \ ${\sf Dom}(f_{\sf e, max}) \ A^{\omega}$ \ $\subsetneqq$
 \ ${\sf in}({\sf cl}({\sf Dom}(f) \ A^{\omega}))$, \ \ \ and

\medskip

 \ \ \  \ \ \    
$g_{\sf e, max} \circ f_{\sf e, max} \ \neq \ (g \circ f)_{\sf e, max}$ .
\end{lem}
{\bf Proof.} For example, let ${\sf domC}(f) = 0^*1$, and 
$f(0^{2n} 1) = 0^{2n+1} 1$, and $f(0^{2n+1} 1) = 0^{2n} 1$, for all 
$n \ge 0$.  Then $f = f_{\sf e, max}$; indeed, any strict extension of $f$ 
would need to make $f(0^m)$ defined for some $m \ge 0$; but such an extension 
to a right-ideal morphism would not agree with $f$ on ${\sf Dom}(f)$ (since
$f$ transposes $0^{2n} 1$ and $0^{2n+1} 1$).  
So, ${\sf Dom}(f_{\sf e, max}) \ \{0,1\}^{\omega}$
$=$ $0^*1 \, \{0,1\}^{\omega}$, whereas ${\sf cl}(0^*1 \, \{0,1\}^{\omega})$
$=$ $\{0,1\}^{\omega}$, and ${\sf in}(\{0,1\}^{\omega}) = \{0,1\}^{\omega}$.

Since $f$ in the above example is injective, we can let $g = f^{-1}$. 
We have $f = f_{\sf e, max}$ and $f^{-1} = (f^{-1})_{\sf e, max}$.  Then 
$\, g_{\sf e, max} \circ f_{\sf e, max} \, $ is the identity restricted to
$0^*1 \, \{0,1\}^*$, whose maximum end-equivalent extension is the full
identity ${\bf 1}_{A^*}$. So in this example, 
$g_{\sf e, max} \circ f_{\sf e, max}$ $\neq$ $(g \circ f)_{\sf e, max}$.  
 \ \ \ $\Box$

%
%

\begin{lem} \label{longEnough}
 \ Let $f_1, f_2$ be right-ideal morphisms with $f_1 \equiv_{\sf end} f_2$,
and let $x \in {\sf Dom}(f_1)$. Then there exists $v \in A^*$ such that 
$\, xv A^* \subseteq {\sf Dom}(f_2)$, and for all $w \in A^*$: 
$f_2(xv w) = f_1(xv w)$.
\end{lem}
{\bf Proof.} If $x \in {\sf Dom}(f_1)$ then the right ideal $xA^*$ 
intersects ${\sf Dom}(f_1)$, hence $xA^*$ intersects ${\sf Dom}(f_2)$
(since ${\sf Dom}(f_1) \equiv_{\sf end} {\sf Dom}(f_2)$).
Thus there exists $x v \in xA^*$ such that $x v \in {\sf Dom}(f_2)$. Hence 
$xv A^* \subseteq {\sf Dom}(f_2)$ (since ${\sf Dom}(f_2)$ is a right ideal). 
We have $f_2(xv w) = f_1(xv w)$ because
$f_1$ and $f_2$ agree where they are both defined.  
 \ \ \ $\Box$

\bigskip

\noindent Just as we saw for right-ideal morphisms in general, 
$f_{\sf e, max}$ can be extended to 
${\sf domC}(f_{\sf e, max}) \, A^{\omega}$.

%
%
%
%

\begin{pro} \label{endequivActionAomega} 
 \ For all right-ideal morphisms $f, g: A^* \to A^*$ we have:

\smallskip
 
 \ \ \  \ \ \ $g \equiv_{\sf end} f$ \ \ iff 
 \ \ $g_{\sf e, max} = f_{\sf e, max}$ on $A^{\omega}$. 
\end{pro}
{\bf Proof.} The implication ``$\Rightarrow$'' is clear from the definitions.
Conversely, if $g_{\sf e, max}$ and $f_{\sf e, max}$ act the same on 
$A^{\omega}$ then ${\sf Dom}(g_{\sf e, max}) \ A^{\omega}$  $=$ 
${\sf Dom}(f_{\sf e, max}) \ A^{\omega}$, hence 
 \ ${\sf cl}({\sf Dom}(g_{\sf e, max}) \ A^{\omega})$  $=$
${\sf cl}({\sf Dom}(f_{\sf e, max}) \ A^{\omega})$. Hence (by Prop.\
\ref{topolEndEqu}), 
${\sf Dom}(g_{\sf e, max}) \equiv_{\sf end} {\sf Dom}(f_{\sf e, max})$.
Since ${\sf Dom}(h) \equiv_{\sf end} {\sf Dom}(h_{\sf e, max})$ for any
righ-ideal morphism $h$, we conclude that 
${\sf Dom}(g) \equiv_{\sf end} {\sf Dom}(f)$.

Since for all $w \in A^{\omega}$, $g_{\sf e, max}(w) = f_{\sf e, max}(w)$, 
we have for all $x \in {\sf Dom}(g) \cap {\sf Dom}(f)$ and all 
$v \in A^{\omega}$: \ $g_{\sf e, max}(xv) = f_{\sf e, max}(xv)$.  Hence 
(since $g_{\sf e, max}(xv) = g_{\sf e, max}(x) \ v$, and similarly for $f$),
$g_{\sf e, max}(x) \ v = f_{\sf e, max}(x) \ v$, for all $v \in A^{\omega}$.
Taking $v = 1 \, 0^{\omega}$ (for example) then
implies  $g_{\sf e, max}(x) = f_{\sf e, max}(x)$,
for all $x \in {\sf Dom}(g) \cap {\sf Dom}(f)$. Since $g_{\sf e, max}$
agrees with $g$ on ${\sf Dom}(g)$ (and similarly for $f$), we conclude that
$g(x) = f(x)$, for all $x \in {\sf Dom}(g) \cap {\sf Dom}(f)$.
Hence, $g \equiv_{\sf end} f$.
 \ \ \ $\Box$

\medskip

\noindent {\bf Remark.} It is also true that \ $g \equiv_{\sf end} f$ \ is
equivalent to the following: \ ${\sf Dom}(g) \equiv_{\sf end} {\sf Dom}(f)$,
and $g$ and $f$ agree on 
${\sf ends}({\sf Dom}(g)) \cap {\sf ends}({\sf Dom}(f))$.  However, 
\ $g \equiv_{\sf end} f$ \ is {\em not} equivalent to the property that
$g$ and $f$ agree on $A^{\omega}$; 
indeed, the actions of $g$ and $f$ on $A^{\omega}$ could have different 
domains (even if $g \equiv_{\sf end} f$).  In Cor.\ \ref{end_boundActCOR}
we will see that $g$ and $f$ agree on $A^{\omega}$ iff 
$g \equiv_{\sf bd} f$ (which is a different congruence than 
$\equiv_{\sf end}$).

\begin{pro} \label{endequivCongr}
The relation $\equiv_{\sf end}$ is a {\em congruence} for right-ideal
morphisms; i.e., for all right-ideal morphisms $f_1, f_2, g$: if
$f_1 \equiv_{\sf end} f_2$, then $f_1 g \equiv_{\sf end} f_2 g$ and
$g f_1 \equiv_{\sf end} g f_2$. 
\end{pro}
{\bf Proof.} The result follows from the next four claims.

\smallskip

\noindent 
{\sf Claim 1.} \ ${\sf Dom}(f_1 g) \equiv_{\sf end} {\sf Dom}(f_2 g)$.

\smallskip

\noindent {\sf Proof.} Let $R$ be a right-ideal that intersects
${\sf Dom}(f_1 g)$, so there exists $x_1 \in R$ such that 
$x_1 \in {\sf Dom}(f_1 g)$; equivalently, $g(x_1) \in {\sf Dom}(f_1)$.
Thus, $g(R)$ intersects ${\sf Dom}(f_1)$, therefore (since 
${\sf Dom}(f_1) \equiv_{\sf end} {\sf Dom}(f_2)$, and $g(R)$ is a right 
ideal) $g(R)$ intersects ${\sf Dom}(f_2)$.  So, for some 
$g(x_2) \in g(R)$ with $x_2 \in R$, $g(x_2) \in {\sf Dom}(f_2)$. 
The latter is equivalent to $x_2 \in {\sf Dom}(f_2 g)$; so $R$ intersects 
${\sf Dom}(f_2 g)$. 
 \ {\small \sf [This proves Claim 1.]}

\smallskip

\noindent
{\sf Claim 2.} \ $f_1 g$ and $f_2 g$ agree on 
${\sf Dom}(f_1 g) \cap {\sf Dom}(f_2 g)$.

\smallskip

\noindent {\sf Proof.} 
Suppose $x \in {\sf Dom}(f_1 g) \cap {\sf Dom}(f_2 g)$.  Then $f_1 g(x)$ 
and $f_2 g(x)$ are defined, so 
$g(x) \in {\sf Dom}(f_1) \cap {\sf Dom}(f_2)$. Hence, since $f_1$ and $f_2$
agree on ${\sf Dom}(f_1) \cap {\sf Dom}(f_2)$ (because 
$f_1 \equiv_{\sf end} f_2$), we have $f_1 g(x) = f_2 g(x)$. 
 \ {\small \sf [This proves Claim 2.]}

\smallskip

\noindent 
{\sf Claim 3.} \ ${\sf Dom}(g f_1) \equiv_{\sf end} {\sf Dom}(g f_2)$.

\smallskip

\noindent {\sf Proof.} Let $R$ be a right ideal that intersects 
${\sf Dom}(g f_1)$, so there exists $x \in R$ such that $g f_1(x)$ is 
defined, 
hence $f_1(x)$ is defined.  Hence, by Lemma \ref{longEnough}, there exists 
$v \in A^*$ such that $xv \in {\sf Dom}(f_1) \cap {\sf Dom}(f_2)$, and 
$f_1(xv) = f_2(xv)$. Thus, $gf_1(xv) = gf_2(xv)$. 
So, $xv \in {\sf Dom}(gf_2)$. Since $R$ is a right ideal, $xv \in R$, hence
$R$ intersects ${\sf Dom}(g f_2)$.

In a similar way one proves that every right ideal that intersects 
${\sf Dom}(g f_2)$ also intersects ${\sf Dom}(g f_1)$.  
 \ {\small \sf [This proves Claim 3.]}

\smallskip

\noindent {\sf Claim 4.} \ $g f_1$ and $g f_2$ agree on
${\sf Dom}(g f_1) \cap {\sf Dom}(g f_2)$.

\smallskip

\noindent {\sf Proof.} Since $f_1 \equiv_{\sf end} f_2$, $f_1$ and
$f_2$ agree on ${\sf Dom}(f_1) \cap {\sf Dom}(f_2)$. 
Moreover, ${\sf Dom}(g f_i) \subseteq {\sf Dom}(f_i)$ (for $i = 1, 2$),
so the Claim holds.  
 \ {\small \sf [This proves Claim 4.]}
 \ \ \ \ \ \ $\Box$

\begin{cor} \label{CORcircMax}   
 \ For all right-ideal morphisms $f, g: A^* \to A^*$, 
 \ $(g_{\sf e, max} \circ f_{\sf e, max})_{\sf e, max}$ \ $=$ 
 \ $(g \circ f)_{\sf e, max}$.  
\end{cor}
{\bf Proof.} By Prop.\ \ref{MaxExtPROP}(1), 
$g_{\sf e, max} \equiv_{\sf end} g$ and $f_{\sf e, max} \equiv_{\sf end} f$. 
Hence by Prop.\ \ref{endequivCongr}, 
$g_{\sf e, max} \circ f_{\sf e, max} \equiv_{\sf end} g \circ f$.
 \ \ \ $\Box$

\smallskip

\begin{defn} {\bf (ends monoid ${\cal M}_{\sf end}^{\sf P}$).} 
The {\em ends monoid} consists of the $\, \equiv_{\sf end}$-classes of 
${\cal RM}^{\sf P}$; the multiplication is the multiplication of 
$\, \equiv_{\sf end}$-classes. It is denoted by 
${\cal M}_{\sf end}^{\sf P}$, or by 
${\cal RM}^{\sf P}\!/\!\equiv_{\sf end}$.
\end{defn}
As a set, ${\cal M}_{\sf end}^{\sf P}$ $=$ 
$\{[f]_{\sf end} : f \in {\cal RM}^{\sf P}\}$.
The multiplication in ${\cal M}_{\sf end}^{\sf P}$ is well-defined since
$\equiv_{\sf end}$ is a congruence, by Prop.\ \ref{endequivCongr}. Hence,
${\cal M}_{\sf end}^{\sf P}$ is a monoid which is a homomorphic image of 
${\cal RM}^{\sf P}$.
The monoid version $M_{2,1}$ of the Richard Thompson group $V$ (a.k.a.\ 
$G_{2,1}$) is a submonoid of ${\cal M}_{\sf end}^{\sf P}$; $M_{2,1}$ is 
defined in \cite{JCBmonThH}; see \cite{CFP} for more information on the 
Thompson group.
 
There is a one-to-one correspondence between $\, \equiv_{\sf end}$-classes 
and maximum end-extensions of elements of ${\cal RM}^{\sf P}$ (by Prop.\
\ref{endequivActionAomega}). So, ${\cal M}_{\sf end}^{\sf P}$ can also be 
defined as

\medskip

 \ \ \ ${\cal M}_{\sf end}^{\sf P}$ \ $=$       
 \ $(\{f_{\sf e, max} : f \in {\cal RM}^{\sf P}\}, \ \cdot)$

\medskip

\noindent with multiplication ``$\cdot$'' defined by

\medskip

 \ \ \ $g_{\sf e, max} \cdot f_{\sf e, max}$ \ $=$         
  \ $(g_{\sf e, max} \circ f_{\sf e, max})_{\sf e, max}$ 
 \ \ ($= \ (g \circ f)_{\sf e, max}$).
\medskip

\noindent Here we used Cor.\ \ref{CORcircMax}.

\bigskip

For the remainder of this section we need some definitions.

For any monoid $M$, the {\em ${\cal L}$-order} (denoted by $\le_{\cal L}$) 
and the {\em $\cal R$-order} (denoted by $\le_{\cal R}$), are defined (for 
any $s, t, u, v \in M$) by
$t \le_{\cal L} s$ iff there exists $m \in M$ such that $t = ms$; and
$v \le_{\cal R} u$ iff there exists $n \in M$ such that $v = un$.
The {\em ${\cal D}$-relation} on $M$ is defined (for $x,y \in M$) by 
$x \equiv_{\cal D} y$ iff there exists $z \in M$ such that 
$x \equiv_{\cal R} z \equiv_{\cal L} y$; equivalently, there exists 
$w \in M$ such that $x \equiv_{\cal L} w \equiv_{\cal R} y$.
A monoid $M$ is called {\em ${\cal D}^0$-simple} if $M$ has only one 
${\cal D}$-class, except possibly for a zero.
These are well-known concepts in semigroup theory; see e.g., 
\cite{Grillet}.

A right ideal $R \subseteq A^*$ is called {\em essential} iff 
$R \equiv_{\sf end} A^*$ (iff ${\sf cl}({\sf ends}(R)) = A^{\omega}$). 
Equivalently, the prefix code that generated $R$ (as a right ideal) is
a {\em maximal} prefix code. 

For any function $f$: $A^* \to A^*$, the relation ${\sf mod}f$ is the 
equivalence relation defined on ${\sf Dom}(f)$ by 
$x_1 \ {\sf mod}f \ x_2$ iff $f(x_1) = f(x_2)$.  The equivalence classes 
of ${\sf mod}f$ are $\{ f^{-1} f(x) : x \in {\sf Dom}(f)\}$.
For two partial functions $g,f: A^* \to A^*$,  we say 
${\sf mod}f \le {\sf mod}g$ (``the relation ${\sf mod}f$ is {\it coarser} 
than ${\sf mod}g$'', or ``${\sf mod}g$ is {\it finer} than ${\sf mod}f$'') 
 \ iff \ ${\sf Dom}(f) \subseteq {\sf Dom}(g)$, and for all 
$x \in {\sf Dom}(f)$: $\, g^{-1} g(x) \subseteq f^{-1} f(x)$. Equivalently, 
${\sf mod}f \le {\sf mod}g$ iff every ${\sf mod}f$-class is a union of 
${\sf mod}g$-classes.

A monoid $M$ is called {\em congruence-simple} iff $M$ is non-trivial and
the only congruences on $M$ are the equality relation and the one-class
congruence.

The {\em length-lexicographic order} on $\{0,1\}^*$ is a well-order, defined
as follows for any $x_1, x_2 \in \{0,1\}^*$:  $x_1 \leq_{\ell \ell} x_2$ 
 \ iff \ $|x_1| < |x_2|$, or $|x_1| = |x_2|$ and $x_1$ precedes $x_2$ in the 
dictionary order on $\{0,1\}^*$ (based on the alphabetic order $0 < 1$).  

The next lemma is the ${\cal RM}^{\sf P}$-version of Prop.\ 2.1 of 
\cite{s1f}.

\begin{lem} \label{regLR}
 \ If $f, r \in {\cal RM}^{\sf P}$ and $r$ is {\em regular} with an inverse
$r' \in {\cal RM}^{\sf P}$ then:

\smallskip

\noindent (1) \ \ $f \leq_{\cal R} r$ \ \ iff \ \ $f = r r' f$
 \ \ iff \ \ ${\sf Im}(f) \subseteq {\sf Im}(r)$.

\smallskip

\noindent (2) \ \ $f \leq_{\cal L} r$ \ \ iff \ \ $f = fr'r$
 \ \ iff \ \ ${\sf mod}f \leq {\sf mod}r$.
\end{lem}
{\bf Proof.} (1) \ $f \leq_{\cal R} r$ iff for some 
$u \in {\cal RM}^{\sf P}$: $f = ru$.  Then $f = r r' r u = r r' f$.
Also, it is straightforward that $f = ru$ implies
 ${\sf Im}(f) \subseteq {\sf Im}(r)$.

Conversely, if ${\sf Im}(f) \subseteq {\sf Im}(r)$ then
 \ ${\bf 1}_{{\sf Im}(f)}$ $=$
${\bf 1}_{{\sf Im}(r)} \circ {\bf 1}_{{\sf Im}(f)}$ $=$
$r \circ r'|_{{\sf Im}(r)} \circ {\bf 1}_{{\sf Im}(f)}$. Hence,
$f \ = \ {\bf 1}_{{\sf Im}(f)} \circ f \ = \ $
$r \circ r'|_{{\sf Im}(r)} \circ {\bf 1}_{{\sf Im}(f)} \circ f \ = \ $
$r \circ r'|_{{\sf Im}(r)} \circ f \ \leq_{\cal R} r$.

\smallskip

\noindent (2) \ $f \leq_{\cal L} r$ iff for some 
$v \in {\cal RM}^{\sf P}$: $f = vr$.  Then $f = v r r' r = f r' r$.
And it is straightforward that $f = vr$ implies
 ${\sf mod}f \leq {\sf mod}r$.

Conversely, if ${\sf mod}f \leq {\sf mod}r$ then for all
$x \in {\sf Dom}(f)$, \ $r^{-1}r(x) \subseteq f^{-1}f(x)$.  And for
every $x \in {\sf Dom}(f)$, \ $\{f(x)\} = f \circ f^{-1} \circ f(x)$.
Moreover,
$f \circ r^{-1} \circ r(x) \subseteq f \circ f^{-1} \circ f(x) = \{f(x)\}$,
and since $r^{-1} \circ r(x) \neq \varnothing$, it follows that
$f \circ r^{-1} \circ r(x) = \{f(x)\}$.  So, $f = f \circ r^{-1} \circ r$.
Moreover, $f \circ r' \circ r(x) \in f \circ r^{-1} \circ r(x) = \{f(x)\}$,
hence $f \circ r' \circ r(x) = f(x)$.  Hence, $f = f r' r \leq_{\cal L} r$.
 \ \ \ $\Box$

\begin{defn} \label{rankfunction} {\bf (rank function).}
 \ For any set $S \subseteq A^*$ the {\em rank function} of $S$ is defined
for all $x \in S$ by $\, {\sf rank}_S(x) = $
$|\{z \in S: z \leq_{\ell \ell} x\}| \, $ (where $\leq_{\ell \ell}$ denotes 
the length-lexicographic order).  
When $x \not\in S$, ${\sf rank}_S(x)$ is undefined.
\end{defn}

We will use padding with a fully time-constructible function in order to
turn any algorithm into a linear-time algorithms.
A ``Turing machine'' will always mean a multi-tape Turing machine.
By definition,
a function $t$: ${\mathbb N} \to {\mathbb N}$ is {\em fully
time-constructible} iff $t$ is total, and increasing, and there exists a 
deterministic Turing machine such that for some $n_0 \in {\mathbb N}$, and 
for all $n \ge n_0$, and for every input of length $n$, the machine runs 
for time exactly $t(n)$.

For example, any polynomial function $n \mapsto c \, (n^d + 1)$ (where
$c, d$ are positive integers), and any exponential function
$n \mapsto c \, d^n$ (where $c > 0$ and $d \ge 2$ are integers) are known
to be fully time-constructible. The sum $t_1(n) + t_2(n)$, and the product
$t_1(n) \cdot t_2(n)$ of two fully time-constructible functions are also
fully time-constructible.
See e.g.\ \cite{HU}, \cite{BDG1}, \cite{Rothe}, for information about
time-constructible functions.

\begin{lem} \label{recTimeConstr}
 \ Let $f$: $\{0,1\}^* \to \{0,1\}^*$ be a partial recursive right-ideal
morphism with decidable domain.  Then we have:

\smallskip

\noindent {\rm (1)} There exists a fully time-constructible function $t$
such that $f$ is computed by a Turing machine with time-complexity
$\le t$.

\smallskip

\noindent {\rm (2)} Let $t$ be any fully time-constructible function such
that $f$ is computed by a Turing machine with time-complexity $\le t$.
Let $F$ be the restriction of $f$ to
 \ $\bigcup_{x \in {\sf domC}(f)} x \, A^{|x| \cdot t(|x|)} \, A^*$.
In other words, $ \, {\sf domC}(F) \, = \,$
$\bigcup_{x \in {\sf domC}(f)} x \, A^{|x| \cdot t(|x|)}$; and
$ \, F(x u v) = f(x) \, u v \, $ for all
$x \in {\sf domC}(f)$, $u \in A^{|x| \cdot t(|x|)}$, $v \in A^*$.

Then $F \equiv_{\sf ends} f$, and $F$ has linear time-complexity and has
linear balance (both bounded from above by the function $n \mapsto 3n$).
\end{lem}
{\bf Proof.} (1) Let $M$ be a deterministic Turing machine that computes
$f$ and that eventually halts on every input. 
We construct a new Turing machine $M'$ for $f$ which has the same
running time for on all inputs of length $\le n$, for all $n$.
On input $x \in A^n$, $M'$ simulates $M$ on all inputs of length $n$, but
only outputs $f(x)$.
If $f(x)$ is not defined, $M'$ produces no output, but since ${\sf Dom}(f)$
is decidable, $M'$ nevertheless
has a time-complexity for every input.
Let $t$ be the time-complexity function of $M'$.
Then $t$ is fully time-constructible, and it is the running time of a
Turing machine that that halts for all inputs and that computes $f$.

(2) Since each set $A^{|x| \cdot t(|x|)}$ is a maximal prefix code, we have
${\sf domC}(F) \equiv_{\sf ends} {\sf domC}(f)$. Also, $F$ is a restriction
of $f$. Hence, $F \equiv_{\sf ends} f$.

To compute $F(w)$ in linear time, we first consider the case where
$w \in {\sf Dom}(F)$, i.e., $w = xuv$ for some $x \in {\sf domC}(f)$,
$u \in A^{|x| \cdot t(|x|)}$, $v \in A^*$.
We first run the Turing machine for $f$ on the prefixes of $xuv$ until $x$
(the smallest prefix on which $f$ is defined) is found.
In detail, each prefix of $w$ is considered in turn, and copied on a work
tape; when the next prefix is considered, one more letter is added on the
right, and the head of the work tape is moved back to the left end. So the
copying of prefixes takes time $\le |x|^2$ ($\le |x| \cdot t(|x|)$).
Checking that the prefix belongs to ${\sf Dom}(f)$ takes time $\le t(|x|)$,
so check all the prefixes takes time $\le |x| \cdot t(|x|)$.
In total, the time to find $x$ and to compute $f(x)$ takes time
$\le  2 \, |x| \cdot t(|x|)$.

Then we check that the rest of the input, namely $uv$, has length
$|x| \cdot t(|x|)$.  Since $|x| \cdot t(|x|)$ is time-constructible,
this can be done in time $\le |x| \cdot t(|x|)$. During this time, $uv$ is
copied to the output tape; this takes time $|uv| = |x| \cdot t(|x|) + |v|$. So
the total time is
$\le 2 \, |x| \, t(|x|) + |x| \, t(|x|) + |v| \le 3 \cdot |xuv| \, $
(since $|u| = |x| \, t(|x|)$).
Thus, $F(xuv)$ is computed in time $\le 3 \cdot |xuv|$.

The complexity bound implies that $|F(xuv)| \le 3 \cdot |xuv|$. For the
input balance we have: $|F(xuv)| = |f(x)| +  |u| + |v| \ge $
$\frac{1}{2} \cdot |x| \cdot t(|x|) + \frac{1}{2} \cdot |u| + |v| \ge$
$\frac{1}{2} \cdot |xuv|$.

To handle the case of an arbitrary input $w$ (not necessarily in
${\sf Dom}(F)$), we follow the same procedure as above, but we add a
counter that stops the computation after time $3 \, |w|$. The machine
rejects, and produces no output, if $w$ has not been found to be in
${\sf Dom}(F)$ by that time.
 \ \ \ $\Box$

\begin{thm} \label{M_is_reg}
 \ The monoid ${\cal M}_{\sf end}^{\sf P}$ is regular, ${\cal D}^0$-simple, 
and finitely generated.
Moreover, every $\equiv_{\sf end}$-class contains a regular element of
${\cal RM}^{\sf P}$.
\end{thm}
{\bf Proof.} An initial remark: Every ${\cal D}^0$-simple monoid is regular;
but we prove regularity separately first because it will be used in the 
proof of ${\cal D}^0$-simplicity. 

For every $f \in {\cal RM}^{\sf P}$ there exists an inverse function $f'$
that is balanced, but that is not necessarily polynomial-time computable;
balance is inherited from $f$ if we restrict the domain of $f'$ to
${\sf Im}(f)$.
Let $T(.)$ be a fully time-constructible upper bound on the
time-complexity
of $f'$ (see Lemma \ref{recTimeConstr}). We can restrict both $f$ and $f'$
in order to reduce the time-complexity (padding argument), while preserving
end-equivalence, as in Lemma \ref{recTimeConstr}: Namely,
we replace ${\sf domC}(f')$ by the prefix code
 \ $\bigcup_{y \in {\sf domC}(f')} \, y \, A^{|y| \cdot T(|y|)}$.
Let $F'$ be this restriction of $f'$, and let the restriction of $f$ be
$F = f \circ F' \circ f = {\sf id}_{{\sf Dom}(F')} \circ f$.
Then $F$ and $F'$ have polynomial (in fact, linear) time-complexity, so
$F, F' \in {\cal RM}^{\sf P}$. Moreover, $F'$ is an inverse of $F$, and
$f \equiv_{\sf end} F$, and $f' \equiv_{\sf end} F'$.  Therefore,
$f \equiv_{\sf end} F = F F' F$ $\in$
$[f]_{\sf end} \cdot [F']_{\sf end} \cdot [f]_{\sf end}$. Hence
$\, [f]_{\sf end} = [f]_{\sf end} \cdot [F']_{\sf end} \cdot [f]_{\sf end}$,
so $[f]_{\sf end}$ is regular in ${\cal M}_{\sf end}^{\sf P}$.
And $F$ ($\equiv_{\sf end} f$) is a regular element of ${\cal RM}^{\sf P}$.

\medskip

Proof of ${\cal D}^0$-simplicity:  
For every non-empty $f \in {\cal RM}^{\sf P}$ there exists 
$f_0 \in {\cal RM}^{\sf P}$ such that $f \equiv_{\sf end} f_0$ and such that
${\sf imC}(f_0)$ is infinite. Indeed, let us pick some 
$x_0 \in {\sf domC}(f)$ and define $f_0$ by 

\smallskip

 \ \ \ \ \    
 ${\sf domC}(f_0) \ = \ ({\sf domC}(f) - \{x_0\}) \ \ \cup \ \ x_0 \ 0^*1$;

\smallskip

 \ \ \ \ \    
$f_0(x_0 \, 0^n 1) = f(x_0) \ 0^n1$ \ \ for all $n \ge 0$, \ and 

\smallskip

 \ \ \ \ \
$f_0(x) = f(x)$ \ \ for all $x \in {\sf domC}(f) - \{x_0\}$. 

\smallskip

\noindent Then ${\sf imC}(f_0)$ contains $f(x_0) \ 0^* 1$, hence it is 
infinite. So, from now on we assume that for $f$ itself, ${\sf imC}(f)$ is 
infinite.

\smallskip

\noindent {\sf Claim:} \ If $f \in {\cal RM}^{\sf P}$ has infinite 
${\sf imC}(f)$, then there exists a right-ideal morphism $g$ with the 
following properties: $g$ is partial recursive with decidable domain, 
${\sf domC}(g)$ is a maximal prefix code, $g$ is injective, and 
${\sf Im}(g) = {\sf Im}(f)$.

\smallskip

\noindent {\sf Proof of the Claim:} We construct $g$ as follows. Let 
$\leq_{\ell \ell}$ denote the length-lexicographic order on $\{0,1\}^*$.  
For any $y \in {\sf imC}(f)$, we can compute
${\sf rank}(y) = |\{z \in {\sf imC}(f): z \leq_{\ell \ell} y\}|$;
computability follows from the fact that $f$ is polynomially balanced.
The function ${\sf rank}$ is injective; it is also onto ${\mathbb N}$ since
${\sf imC}(f)$ is infinite. We define $g(0^n 1)$ to be the element
$y \in {\sf imC}(f)$ such that ${\sf rank}(y) = n$ (for any $n \ge 0$).
So, $g$ is injective and $g^{-1}(y) = 0^{{\sf rank}(y)} 1$ for all 
$y \in {\sf imC}(f)$.
We have ${\sf domC}(g) = 0^* 1$, which is a maximal prefix code; obviously,
$0^* 1$ is a decidable language.
Then $g$ is partial recursive with decidable domain, and injective, and 
${\sf Im}(f) = {\sf Im}(g)$. This proves the Claim. 

\smallskip

Let $t(.)$ be a fully time-constructible upper bound on the time
complexities of $g$, $g^{-1}$, $f$ and $f'$.
Then, by padding $f$ and $g$ as in Lemma \ref{recTimeConstr} we obtain
functions $f_1, g_1 \in {\cal RM}^{\sf P}$ such that
$f \equiv_{\sf end} f_1$, $g \equiv_{\sf end} g_1$, both $f_1$ and $g_1$
are computable in linear time, and both $f_1$ and $g_1$ are regular in
${\cal RM}^{\sf P}$.
Moreover, ${\sf domC}(g_1)$ ($\equiv_{\sf end} 0^* 1$) is a maximal prefix
code (equivalently, ${\sf Dom}(g_1)$ is an essential right ideal),
$g_1$ is injective, and ${\sf Im}(f_1) = {\sf Im}(g_1)$.

Since ${\sf Im}(f_1) = {\sf Im}(g_1)$ and $f_1, g_1$ are regular elements of
${\cal RM}^{\sf P}$, Lemma \ref{regLR}(1) implies $f_1 \equiv_{\cal R} g_1$. 

Since $g_1$ is injective, the relation ${\sf mod}g_1$ is the equality 
relation on ${\sf Dom}(g_1)$. Hence, since $g_1$ is a regular element of 
${\cal RM}^{\sf P}$, Lemma \ref{regLR}(2) implies that 
$g_1 \equiv_{\cal L} {\bf 1}_{{\sf Dom}(g_1)}$ (the identity map restricted 
to ${\sf Dom}(g_1)$). Since ${\sf domC}(g_1)$ is a maximal prefix code, 
${\sf Dom}(g_1)$ is an essential right ideal; equivalently, 
${\bf 1}_{{\sf Dom}(g_1)} \equiv_{\sf end} {\bf 1}$.

Overall we now have $f \equiv_{\sf end} f_1 \equiv_{\cal R} g_1$
$\equiv_{\cal L}$ ${\bf 1}_{{\sf Dom}(g_1)} \equiv_{\sf end} {\bf 1}$.
Hence, in ${\cal M}_{\sf end}^{\sf P}$, 
$[f]_{\sf end} \equiv_{\cal D} [{\bf 1}]_{\sf end}$.

\medskip

Proof of finite generation: The proof is based on the fact that 
${\cal RM}^{\sf P}$ has evaluation maps for programs with bounded
balance and time-complexity. This was described in detail in
Section 4 of \cite{s1f} and Section 2 of \cite{BiInfGen}. We briefly give 
the definition here: For a polynomial $q_2$ such that 
$q_2(n) = a \, (n^2 +1) \,$ (for some fixed large constant $a$), we define 
an {\em evaluation map} ${\sf evR}_{q_2}^C \in {\cal RM}^{\sf P}$ by 

\smallskip

 \ \ \ \ \ ${\sf evR}_{q_2}^C({\sf code}(w) \ 11 \, u v)$ \ $=$
 \ ${\sf code}(w) \ 11 \ \phi_w(u) \ v$, 

\smallskip

\noindent for all Turing machine programs $w$ with balance and 
time-complexity $\le q_2$,  and all $u \in {\sf domC}(\phi_w)$, and 
$v \in A^*$.  Here, $\phi_w \in {\cal RM}^{\sf P}$ denotes the function
with program $w$.  Then we have

\smallskip

 \ \ \ \ \     
 $\phi_w \ = \ \pi'_{|{\sf code}(w) \, 11|} \circ {\sf evR}_{q_2}^C$
    $\circ$  $\pi_{{\sf code}(w) \, 11}$.

\smallskip

\noindent where $\pi'_n$ is defined by $\pi'_n(x_1 x_2) = x_2$ whenever
$x_1, x_2 \in A^*$ with $|x_1| = n$ (and $\pi'_n$ is undefined on other 
arguments); and $\pi_u$ is defined by $\pi_u(x) = u x$ for all 
$u, x \in A^*$. See \cite{BiInfGen} for the proof that such a function
${\sf evR}_{q_2}^C$ exists. 

In the proof of regularity of ${\cal M}_{\sf end}^{\sf P}$ above, we saw 
that every $\phi_v \in {\cal RM}^{\sf P}$ is end-equivalent to some
$\phi_w \in {\cal RM}^{\sf P}$ such that $\phi_w$ has linear
time-complexity (in fact, it is $\le 3 \, n$, by Lemma \ref{recTimeConstr}).
We can obtain $\phi_w$ as $\phi_w = \phi_v \circ {\bf 1}_{P_w}$, where

\smallskip

 \ \ \  \ \ \ $P_w \ = \ $
  $\bigcup_{x \in {\sf domC}(\phi_v)} x \, A^{T(|x|)^2}$;

\smallskip

\noindent here, $T(.)$ is the time-complexity of $\phi_w$.
Since $T(.)$ is a polynomial of the form $c \, (n^2 +1)$, the function
$n \mapsto n \cdot T(n)^2$ is fully time-constructible.
Since the time-complexity of $\phi_w$ is linear with coefficient $\le 3$
(by Lemma \ref{recTimeConstr}),
the evaluation map ${\sf evR}_{q_2}^C$ can evaluate $\phi_w$ without any
need for further padding; so we have

\smallskip

 \ \ \ \ \ $\phi_v$ \ $\equiv_{\sf end}$ 
 \ $\phi_w \ = \ \pi'_{|{\sf code}(w) \, 11|} \circ {\sf evR}_{q_2}^C$
    $\circ$  $\pi_{{\sf code}(w) \, 11}$ .

\smallskip

\noindent
So, ${\cal M}_{\sf end}^{\sf P}$ is generated by $\, \{[\pi_0]_{\sf end},$
$[\pi_1]_{\sf end}, [\pi_1']_{\sf end}, [{\sf evR}_{q_2}^C]_{\sf end} \}$.
 \ \ \ $\Box$

\bigskip

In the proof of Theorem \ref{Mcongr_simple} the concept of
${\cal J}^0$-simplicity is used. By definition, a monoid $M$ with a zero
is ${\cal J}^0$-simple iff for all non-zero elements
$a, b \in M$ there exist $x_1, x_2, x_3, x_4 \in M$ such that
$a = x_1 b x_2$ and $b = x_3 a x_4$. For more information on the
${\cal J}$-relation and the ${\cal J}$-preorder, see e.g.\ \cite{Grillet}.
Obviously, ${\cal D}^0$-simplicity implies ${\cal J}^0$-simplicity.

\begin{thm} \label{Mcongr_simple}
 \ ${\cal M}_{\sf end}^{\sf P}$ is congruence-simple. 
\end{thm}
{\bf Proof.} The proof is similar to the proof of congruence-simplicity of
the Thompson-Higman monoid $M_{2,1}$ in \cite{BiThompsMonV3}.
Let {\bf 0} be the $\, \equiv_{\sf end}$-class of the empty 
map; this class consists only of the empty map $0$. 
When $\cong$ is any congruence on ${\cal M}_{\sf end}^{\sf P}$ that is not 
the equality relation, we will show that the whole monoid 
${\cal M}_{\sf end}^{\sf P}$ is congruent to {\bf 0}.  We will make use of 
${\cal J}^0$-simplicity of ${\cal M}_{\sf end}^{\sf P}$, which follows 
from its ${\cal D}^0$-simplicity (and also from the ${\cal J}^0$-simplicity 
of ${\cal RM}^{\sf P}$, Prop.\ 2.7 in \cite{s1f}).

Since $\cong$ is a congruence on ${\cal M}_{\sf end}^{\sf P}$, and since 
$\equiv_{\sf end}$ is a congruence on ${\cal RM}^{\sf P}$, it follows that 
$\cong$ can also be defined as a congruence on ${\cal RM}^{\sf P}$ that 
is coarser than $\equiv_{\sf end}$.  We will show that if $\cong$ on 
${\cal RM}^{\sf P}$ is not $\equiv_{\sf end}$, then $\cong$ is the trivial 
(one-class) congruence.

\smallskip

\noindent {\sf Case (0):} \
Assume that $\Phi \cong {\bf 0}$ for some element $\Phi \neq {\bf 0}$ in 
${\cal M}_{\sf end}^{\sf P}$. Then for all $\alpha, \beta \in $
${\cal M}_{\sf end}^{\sf P}$ we have obviously 
$\alpha \, \Phi \, \beta \cong {\bf 0}$. Moreover, by 
${\cal J}^0$-simplicity, ${\cal M}_{\sf end}^{\sf P}$ $=$
$\{\alpha\, \Phi\, \beta : \alpha,\beta \in {\cal M}_{\sf end}^{\sf P}\}$,
since $\Phi \neq {\bf 0}$. Hence all elements are congruent to {\bf 0}.

\smallskip

For the remainder of the proof we let $\varphi, \psi$ $\in$ 
${\cal RM}^{\sf P} - \{0\}$ be representatives of two different 
$\, \equiv_{\sf end}$-classes (i.e., $\varphi \not\equiv_{\sf end} \psi$)
such that $[\varphi]_{\sf end} \cong [\psi]_{\sf end}$. 
Notation: For any $u, v \in A^*$, $(v \leftarrow u)$ denotes the right-ideal 
morphism $u w \mapsto v w$ (for all $w \in A^*$).

\smallskip

\noindent {\sf Case (1):}
 \ ${\sf Dom}(\varphi) \not\equiv_{\sf end} {\sf Dom}(\psi)$.

Then there exists $x_1 \in A^*$ such that $x_1 A^*$ intersects 
${\sf Dom}(\varphi)$ (e.g., at $x_0$), but $x_1 A^*$ does not intersect 
${\sf Dom}(\psi)$. Then $x_0A^* \subseteq {\sf Dom}(\varphi)$, but
 \ ${\sf Dom}(\psi) \cap x_0A^* = \varnothing$. Or, vice versa, there exists 
$x_0 \in A^*$ such that $x_0A^* \subseteq {\sf Dom}(\psi)$, but
 \ ${\sf Dom}(\varphi) \cap x_0A^* = \varnothing$.  Let us assume the former.

Letting $\beta = (x_0 \leftarrow x_0)$, we have
$\varphi \, \beta(.) = (\varphi(x_0) \leftarrow x_0)$. But
$\psi \, \beta(.) = 0$, since 
$x_0A^* \cap {\sf Dom}(\psi) = \varnothing$.  So, 
$[\varphi \, \beta]_{\sf end} \cong [\psi \, \beta]_{\sf end} = {\bf 0}$, 
but $[\varphi \, \beta]_{\sf end} \neq {\bf 0}$.
Hence, applying case (0) to $\Phi = [\varphi \, \beta]_{\sf end}$ we 
conclude that the entire monoid ${\cal M}_{\sf end}^{\sf P}$ is congruent 
to {\bf 0}.

\smallskip

\noindent {\sf Case (2.1):}
 \ ${\sf Dom}(\varphi) \equiv_{\sf end} {\sf Dom}(\psi)$ \ and
 \ ${\sf Im}(\varphi) \not\equiv_{\sf end} {\sf Im}(\psi)$.

Then there exists $y_0 \in A^*$ such that
$y_0A^* \subseteq {\sf Im}(\varphi)$, but
${\sf Im}(\psi) \cap y_0A^* = \varnothing$;  or, vice versa,
$y_0A^* \subseteq {\sf Im}(\psi)$, but
${\sf Im}(\varphi) \cap y_0A^* = \varnothing$.
Let us assume the former. 

Let $x_0 \in A^*$ be such that $y_0 = \varphi(x_0)$. Then
$(y_0 \leftarrow y_0) \circ \varphi \circ (x_0 \leftarrow x_0)$
$=$ $(y_0 \leftarrow x_0)$.
On the other hand,
$(y_0 \leftarrow y_0) \circ \psi \circ (x_0 \leftarrow x_0) = {\bf 0}$.
Indeed, if $x_0A^* \cap {\sf Dom}(\psi) = \varnothing$ then for all
$w \in A^*$: $\psi \circ (x_0 \leftarrow x_0)(x_0w) = \psi(x_0w)$ $=$ 
$\varnothing$. And if $x_0A^* \cap {\sf Dom}(\psi) \neq  \varnothing$ then
for those $w \in A^*$ such that $x_0w \in {\sf Dom}(\psi)$ we have
$(y_0 \leftarrow y_0) \circ \psi \circ (x_0 \leftarrow x_0)(x_0w)$ $=$
$(y_0 \leftarrow y_0)(\psi(x_0w)) = \varnothing$, since
${\sf Im}(\psi) \cap y_0A^* = \varnothing$.  Now case (0) applies to
${\bf 0} \neq \Phi$ $=$
$[(y_0 \leftarrow y_0) \circ \varphi \circ (x_0 \leftarrow x_0)]_{\sf end}$
$ \cong {\bf 0}$;
hence all elements of ${\cal M}_{\sf end}^{\sf P}$ are congruent to {\bf 0}.

\smallskip

\noindent {\sf Case (2.2):}
 \ ${\sf Dom}(\varphi) \equiv_{\sf end} {\sf Dom}(\psi)$ \ and 
  \ ${\sf Im}(\varphi) \equiv_{\sf end} {\sf Im}(\psi)$.

Then we can restrict $\varphi$ and $\psi$ to
${\sf Dom}(\varphi) \cap {\sf Dom}(\psi)$ ($\equiv_{\sf end}$
${\sf Dom}(\varphi) \equiv_{\sf end} {\sf Dom}(\psi)$), by choice of 
representatives in $[\varphi]_{\sf end}$, respectively $[\psi]_{\sf end}$; so 
now \ ${\sf domC}(\varphi) = {\sf domC}(\psi)$. Since $\varphi \neq \psi$,
there exist $x_0 \in {\sf domC}(\varphi) = {\sf domC}(\psi)$ and
$y_0 \in {\sf Im}(\varphi)$, $y_1 \in {\sf Im}(\psi)$, such that 
$\varphi(x_0) = y_0 \neq y_1 = \psi(x_0)$. We have two subcases.

\smallskip

\noindent {\sf Subcase (2.2.1):} \ $y_0$ and $y_1$ are not prefix-comparable.

Then
 \ $(y_0 \leftarrow y_0) \circ \varphi \circ (x_0 \leftarrow x_0)$ $=$ 
$(y_0 \leftarrow x_0)$.

On the other hand,
$(y_0 \leftarrow y_0) \circ \psi \circ (x_0 \leftarrow x_0)(x_0 w)$ $=$
$(y_0 \leftarrow y_0)(y_1 w) = \varnothing$ \ for all $w \in A^*$ (since 
$y_0$ and $y_1$ are not prefix-comparable).  So, 
$(y_0 \leftarrow y_0) \circ \psi \circ (x_0 \leftarrow x_0) = {\bf 0}$.
Hence case (0) applies to \ ${\bf 0} \neq\Phi$ $=$
$[(y_0 \leftarrow y_0) \circ \varphi \circ (x_0 \leftarrow x_0)]_{\sf end}$
$\cong {\bf 0}$.

\smallskip

\noindent {\sf Subcase (2.2.2):} \ $y_0$ is a prefix of $y_1$, and
$y_0 \neq y_1$. (The case where $y_0$ is a prefix of $y_1$ is similar.)

Then $y_1 = y_0 a u_1$ for some $a \in A$, $u_1 \in A^*$.
Letting $b \in A - \{a\}$, and $y_2 = y_0 b$, we obtain a string $y_2$ that
is not prefix-comparable with $y_1$.  Now,
$\, (y_2 \leftarrow y_2) \circ \varphi \circ (x_0 \leftarrow x_0)(x_0 b)$
$ \ = \ $ $(y_2 \leftarrow y_2)(y_0 b) \ = \ y_2$; so, $\Phi = $
$(y_2 \leftarrow y_2) \circ \varphi \circ (x_0 \leftarrow x_0) \ne {\bf 0}$.
But for all $w \in A^*$,
$ \, (y_2 \leftarrow y_2) \circ \psi \circ (x_0 \leftarrow x_0)(x_0w)$
$ \ = \ (y_2 \leftarrow y_2)(y_1 w) \ = \ \varnothing$, since
$y_2$ and $y_1$ are not prefix-comparable; so,
$(y_2 \leftarrow y_2) \circ \psi \circ (x_0 \leftarrow x_0) = {\bf 0}$.
Thus, case 0 applies to \ ${\bf 0} \neq \Phi = $
$(y_2 \leftarrow y_2) \circ \varphi \circ (x_0 \leftarrow x_0)$
$\cong (y_2 \leftarrow y_2) \circ \psi \circ (x_0 \leftarrow x_0)$
$= {\bf 0}$.
 \ \ \ $\Box$

\begin{pro} \label{Mgroupunits}
 \ The group of units of ${\cal M}_{\sf end}^{\sf P}$ is  

\smallskip

 \ \ \ \ \ {\rm \ $\{[f]_{\sf end} \ : \ f \in {\cal RM}^{\sf P}$ and 
 $f$ is a bijection between two essential right ideals of $A^*\}$.     }
\end{pro}
{\bf Proof.} If $f \in {\cal RM}^{\sf P}$ is a bijection between essential
right ideals, then $f$ is also a bijection from $R_1 = {\sf Dom}(f)$ onto
$R_2 = {\sf Im}(f)$; and $R_1$ and $R_2$ are decidable subsets of $A^*$
(since $R_1 \in {\sf P}$ and $R_2 \in {\sf NP}$). Hence
$f^{-1}$: $R_2 \to R_1$ is partial recursive, and has decidable domain and
image.  Also, $f \circ f^{-1} = {\sf id}_{R_2}$, and
$f^{-1} \circ f = {\sf id}_{R_1}$.  Since $R_1, R_2$ are essential right
ideals, ${\sf id}_{R_2}$ $\equiv_{\sf end}$ ${\sf id}$ $\equiv_{\sf end}$
${\sf id}_{R_1}$.
So, $[f]_{\sf end} \cdot [f^{-1}]_{\sf end} = [{\sf id}]_{\sf end}$ $=$
$[f^{-1}]_{\sf end} \cdot [f]_{\sf end}$. By Lemma \ref{recTimeConstr},
$f^{-1}$ is $\equiv_{\sf end}$-equivalent to an element of
${\cal RM}^{\sf P}$ (with linear time-complexity and linear balance).
Hence, $[f^{-1}]_{\sf end} \in$ ${\cal M}_{\sf end}^{\sf P}$; so
$[f]_{\sf end}$ belongs to the group of units.

Conversely, suppose $[F]_{\sf end} \equiv_{\cal H} [{\sf id}]_{\sf end}$ in
${\cal M}_{\sf end}^{\sf P}$, where $F \in {\cal RM}^{\sf P}$.  Then there
exists $f \in {\cal RM}^{\sf P}$ such that $f \equiv_{\sf end} F$, and $f$
regular in ${\cal RM}^{\sf P}$ (by Theorem \ref{M_is_reg}).
Since $[f]_{\sf end}$ $\equiv_{\cal L}$ $[{\sf id}]_{\sf end}$, there
exist $f_2, {\sf id}_{R_2} \in {\cal RM}^{\sf P}$ such that
$f_2 \circ f = {\sf id}_{R_2}$; and ${\sf id}_{R_2} \in {\cal RM}^{\sf P}$
implies $R_2 \in {\sf P}$.
Since $[f]_{\sf end}$ $\equiv_{\cal R}$ $[{\sf id}]_{\sf end}$, there
exist $f_1, {\sf id}_{R_1} \in {\cal RM}^{\sf P}$ such that
$f \circ f_1 = {\sf id}_{R_1}$; and ${\sf id}_{R_1} \in {\cal RM}^{\sf P}$
implies $R_1  \in {\sf P}$.  Since
${\sf id}_{R_2} \equiv_{\sf end} {\sf id} \equiv_{\sf end} {\sf id}_{R_1}$,
$R_1$ and $R_2$ are essential.
Since $f$ and ${\sf id}_{R_2}$ are regular, Lemma \ref{regLR}(1) implies
${\sf Im}(f) = R_2$. Since $f$ and ${\sf id}_{R_1}$ are regular, Lemma
\ref{regLR}(2) implies $f$ is injective and ${\sf Dom}(f) = R_1$.
Hence, $f$ is a bijection from $R_1$ onto $R_2$.
 \ \ \ $\Box$

%

\bigskip

\noindent We prove next that in the definition of
${\cal M}_{\sf end}^{\sf P}$ we can replace ${\cal RM}^{\sf P}$ by the
monoid ${\cal RM}^{\sf rec}$, defined by

\smallskip 

 \ \ \ \ \  ${\cal RM}^{\sf rec} \ = \ $
 $\{f \ : \ f$ is a right-ideal morphism on $A^*$ that is partial recursive,
 ${\sf Dom}(f)$ is

\hspace{1.4in} decidable, and $f$ has a total recursive input-output
    balance\}.

\medskip

\noindent Recall that $[f]_{\sf end} \in {\cal RM}^{\sf P}$ is defined by
$\, [f]_{\sf end}$ $=$
$\{F \in {\cal RM}^{\sf P} : \, F \equiv_{\sf end} f\}$.

\begin{pro} \label{RMRecPolyRMPpoly}
The monoid  ${\cal RM}^{\sf rec}\!/\!\equiv_{\sf end}$ is isomorphic to
${\cal RM}^{\sf P}\!/\!\equiv_{\sf end}$ \ ($= {\cal M}_{\sf end}^{\sf P}$).
\end{pro}
{\bf Proof.} Let us show that the map $H$: $[f]_{\sf end} \ \longmapsto$
$ \ \{F \in {\cal RM}^{\sf rec} : F \equiv_{\sf end} f\} \,$
(for all $f \in {\cal RM}^{\sf P}$) is a bijection from
${\cal M}_{\sf end}^{\sf P}$ onto
${\cal RM}^{\sf rec}\!/\!\equiv_{\sf end}$.
The map  $H$ is injective because different $\equiv_{\sf end}$-classes are
disjoint, in both ${\cal RM}^{\sf P}$ and ${\cal RM}^{\sf rec}$. The map is
also surjective because for every $g \in {\cal RM}^{\sf rec}$ there exists
$g_{\sf pad} \in {\cal RM}^{\sf P}$ such that
$g_{\sf pad} \equiv_{\sf end} g$.  We can take $g_{\sf pad}$ to be the
restriction of $g$ to
$\, \bigcup_{y \in {\sf domC}(g)} \, y \, A^{|y| \cdot t(|y|)} \, A^*$,
as in Lemma \ref{recTimeConstr}.
Moreover, $H$ is a homomorphism since $\equiv_{\sf end}$ is a congruence.
 \ \ \ $\Box$

\bigskip

\noindent {\bf Question:} 
We proved that ${\cal M}_{\sf end}^{\sf P}$ is a congruence-simple 
homomorphic image of ${\cal RM}^{\sf P}$. Does ${\cal RM}^{\sf P}$ have 
other congruence-simple homomorphic images?


\section{Bounded end-equivalence}

\begin{defn} \label{equiv_boundDEF} {\bf (bounded end-equivalence of sets).}

\noindent {\bf (1)} Two sets $L_1, L_2 \subset A^*$ are {\em boundedly 
end-equivalent} (denoted by $L_1 \equiv_{\sf bd} L_2$) \ iff 
 \ $L_1 \equiv_{\sf end} L_2$, and there exists a total function
$\beta$: ${\mathbb N} \to {\mathbb N}$ such that for all $x_1 \in L_1$ and 
$x_2 \in L_2$: $\, x_1 \parallel_{\sf pref} x_2 \, $ implies 
$|x_1| \le \beta(|x_2|)$ and $|x_2| \le \beta(|x_1|)$.  

\smallskip

\noindent {\bf (2)} More generally, let $\cal T$ be any non-empty family of 
total functions ${\mathbb N} \to {\mathbb N}$ such that:

\noindent $\bullet$ 
 \ $\cal T$ contains {\em upper bounds on sum and composition};
this means that for all $\tau_1, \tau_2 \in {\cal T}$ there exist      
$\tau_3, \tau_4 \in {\cal T}$ such that for all $n \in {\mathbb N}$:
$\, \tau_1(n) + \tau_2(n) \le \tau_3(n)$, and 
$\, \tau_1(\tau_2(n)) \le \tau_4(n)$.

\noindent $\bullet$ \ There exists $\tau \in {\cal T}$ such that $\tau$ is
an increasing function, and $n \le \tau(n)$ for all $n \in {\mathbb N}$.

\smallskip

Two sets $L_1, L_2 \subset A^*$ are {\em $\cal T$-end-equivalent} (denoted
by $L_1 \equiv_{\cal T} L_2$) \ iff \ $L_1 \equiv_{\sf end} L_2$, and
there exists a function $\tau \in {\cal T}$ such that for all $x_1 \in L_1$
and $x_2 \in L_2$: $x_1 \parallel_{\sf pref} x_2$ implies
$|x_1| \le \tau(|x_2|)$ and $|x_2| \le \tau(|x_1|)$.
\end{defn}
Note that the bounding function $\beta$ or $\tau$ for 
$L_1 \equiv_{\sf bd} L_2$  
depends on $L_1$ and $L_2$. The only assumption on the function 
$\beta$: ${\mathbb N} \to {\mathbb N}$ is that it is {\em total} on 
${\mathbb N}$; no computability assumptions are made.

\bigskip

\noindent {\bf Examples and counter-examples:}

\smallskip

For any prefix code $P \subset A^*$ and any total function $\beta$:
${\mathbb N} \to {\mathbb N}$ we have
$\, \bigcup_{x \in P} \, x \, A^{\beta(|x|)} \equiv_{\sf bd} P$.

For the prefix codes $\{0,1\}$ and $0^* 1$ we have
$\{0,1\} \equiv_{\sf end} 0^* 1$, but
$\{0,1\} \not\equiv_{\sf bd} 0^* 1$.

When $P$ is a prefix code, $P \not\equiv_{\sf bd} P \, A^*$; in this,
$\equiv_{\sf bd}$ differs from $\equiv_{\sf end}$.
And $A^*$ is not $\equiv_{\sf bd}$-equivalent to itself, so
$\equiv_{\sf bd}$ is not reflexive in general.
When $L$ is a prefix code, $L$ is boundedly end-equivalent to itself. If
$L$ is a union of two prefix codes, then $L$ might not be boundedly
end-equivalent to itself (e.g., $\{0,1\} \, \cup \, \{0^n 1: n \ge 0\}$).
 From here on we will use $\equiv_{\sf bd}$ only between prefix codes.

\bigskip

Closure of $\cal T$ under composition guarantees that $\equiv_{\cal T}$ is 
transitive.
Typical examples of families $\cal T$ as above are the following (where we
only take those functions that are increasing and satisfy $n \le \tau(n)$): 

\noindent $\bullet$
 \ ${\mathbb N}^{\mathbb N}$, i.e., the family of all  total functions on
${\mathbb N}$; then $\equiv_{\cal T}$ is $\equiv_{\sf bd}$.

\noindent $\bullet$ \ {\sf rec} = the family of all partial recursive 
functions ${\mathbb N} \to {\mathbb N}$ with decidable domain, i.e., the 
partial recursive functions that are extendable to total recursive
functions.

\noindent $\bullet$ 
 \ {\sf E3} = the family of all {\em elementary recursive functions}, i.e., 
level 3 of the Grzegorczyk hierarchy; these are the primitive recursive 
functions with size bounded by a constant iteration of exponentials. 


\noindent $\bullet$ \ {\sf poly} = the family of all polynomials with 
non-negative integer coefficients.

\noindent $\bullet$ \ {\sf lin} = the family of all affine functions of
the form $n \mapsto a \, n + b$ (where $a \ge 1, b \ge 0$). 

\bigskip

\noindent We have the following Cantor-space characterization of 
$\, \equiv_{\sf bd}$ between prefix codes:

\begin{pro} \label{equiv_boundEndsPROP} 
 \ For prefix codes $P_1, P_2 \subset A^*$ we have 
$P_1 \equiv_{\sf bd} P_2$ \ iff \ ${\sf ends}(P_1) = {\sf ends}(P_2)$.  
\end{pro}
{\bf Proof.} $[\Leftarrow]$ If ${\sf ends}(P_1) = {\sf ends}(P_2)$ then
by applying closure we obtain $P_1 \equiv_{\sf end} P_2$ (by Prop.\ 
\ref{topolEndEqu}). To prove boundedness of this end-equivalence, let
$x_1 \in P_1$, $x_2 \in P_2$ be such that $x_1 \parallel_{\sf pref} x_2$.
Let us assume $x_1$ is a prefix of $x_2$ (if $x_2$ is a prefix of $x_1$ the
reasoning is symmetric). The existence of a total function $\beta$: 
${\mathbb N} \to {\mathbb N}$ such that $|x_2| \le \beta(|x_1|)$ for
every $x_2$ that has $x_1$ as a prefix, is equivalent to the finiteness of 
$x_1 A^* \cap P_2$ for every $x_1 \in P_1$. Indeed, the lengths of the words
in a set $S \subseteq A^*$ (over a finite alphabet $A$) are bounded iff that 
set $S$ is finite.

Since ${\sf ends}(P_1) = {\sf ends}(P_2)$, every end that passes through
$x_1$ (i.e., that belongs to the subtree $x_1 A^*$ of $A^*$) intersects
$P_2$. Hence, the tree 
$x_1 A^* - P_2 A^+ = (x_1 A^* - P_2 A^*) \, \cup \, (P_2 \cap x_1 A^*)$ 
has no infinite path.  By the K\"onig Infinity Lemma, this implies that this 
tree is finite. Hence $x_1 A^* \cap P_2$, which is the set of leaves of this 
finite tree, is finite.

\smallskip

\noindent $[\Rightarrow]$ If $P_1 \equiv_{\sf bd} P_2$ with bounding
function $\beta$: ${\mathbb N} \to {\mathbb N}$, consider 
$x_1 w \in {\sf ends}(P_1)$, with $x_1 \in P_1$ and $w \in A^{\omega}$.
If $x_1$ has a prefix $x_2 \in P_2$ then obviously, 
$x_1 w \in {\sf ends}(P_2)$. 

Let us assume next that $x_1$ does not have a prefix in $P_2$; we want to
show that in this case too, $x_1 w \in {\sf ends}(P_2)$.
For every $n \in {\mathbb N}$, let $w_n$ be the prefix of length $n$ of $w$.
Since $P_1 \equiv_{\sf end} P_2$, every right ideal $x_1 w_n A^*$ intersects
$P_2 A^*$; so, $x_1 w_n u_n = x_{2,n} v_n$ for some $x_{2,n} \in P_2$ and
some $u_n, v_n \in A^*$. It follows that $x_1 \parallel_{\sf pref} x_{2,n}$,
hence $x_1$ is a prefix of $x_{2,n}$ (since we assumed that $x_1$ does not
have a prefix in $P_2$).  Hence, $|x_{2,n}| \le \beta(|x_1|)$.
It also follows from $x_1 w_n u_n = x_{2,n} v_n$ that
$x_1 w_n \parallel_{\sf pref} x_{2,n}$, so either $x_1 w_n$ is a prefix of
$x_{2,n}$, or $x_{2,n}$ is a prefix of $x_1 w_n$ (while $x_1$ is also a
prefix of $x_{2,n}$).

Case 1: $x_1 w_n$ is a prefix of $x_{2,n}$.
Then $|x_{2,n}| \ge |x_1| + n$.
If we choose $n$ so that $n > \beta(|x_1|)$, this case is ruled out.

Case 2: $x_{2,n}$ is a prefix of $x_1 w_n$.
Then $x_{2,n}$ is a vertex on the end $x_1 w$ (between $x_1$ and $x_1 w_n$),
hence $x_1 w$ is equal to an end through $x_{2,n}$, so
$x_1 w \in {\sf ends}(P_2)$.
 \ \ \ $\Box$

\bigskip

\noindent {\bf Remarks:} 

\noindent
{\bf (1)} Prop.\ \ref{equiv_boundEndsPROP}  was proved for prefix codes. 
When $P_1, P_2 \subseteq A^*$ are not prefix codes, the proposition does not
always hold. E.g., $A^* \not\equiv_{\sf bd} A^*$ (non-reflexivity, as we 
saw), but obviously ${\sf ends}(A^*) = {\sf ends}(A^*)$.

One could argue that ${\sf ends}(L_1) = {\sf ends}(L_2)$ is the more 
reasonable definition of ``$L_1 \equiv_{\sf bd} L_2$''. But our 
definition of $\equiv_{\sf bd}$ (Def.\ \ref{equiv_boundDEF}) has the 
advantage of generalizing to $\equiv_{\cal T}$.
 
In any case, since we will use $\equiv_{\sf bd}$ only with prefix codes, the
question doesn't matter in this paper.  
 
\smallskip

\noindent {\bf (2)}
The relation $\equiv_{\sf bd}$ can be generalized to a pre-order, denoted
by $\subseteq_{\sf bd}$: For prefix codes $P_1, P_2 \subset A^*$ we define
 \ $P_1 \subseteq_{\sf bd} P_2$ \ iff
 \ ${\sf ends}(P_1) \subseteq {\sf ends}(P_2)$. Equivalently,
$P_1 \subseteq_{\sf bd} P_2$  \ iff \ there exists $Q \subseteq P_2$ such
that $P_1 \equiv_{\sf bd} Q$. Similarly, $\equiv_{\cal T}$ can be generalized
by defining  \ $P_1 \subseteq_{\cal T} P_2$ \ iff \ there exists
$Q \subseteq P_2$ such that $P_1 \equiv_{\cal T} Q$.

\medskip

\noindent {\bf Notation:} For any right ideal $R \subseteq A^*$, the prefix
code that generates $R$ (as a right ideal) is denoted by ${\sf prefC}(P)$.

\begin{pro} \label{infUnionInterBd}
 \ For any prefix code $P \subset A^*$, we have:

\smallskip

 \ \ \ $\bigcap \, \{ {\sf ends}(Q): Q$ is a prefix code and 
    $Q \equiv_{\sf bd} P \} \ = \ {\sf ends}(P)$,

\medskip

 \ \ \ $\bigcap \, \{ Q A^* : Q$ is a prefix code and
    $Q \equiv_{\sf bd} P \} \ = \ \varnothing$, 

\medskip

 \ \ \ ${\sf prefC}\big(\bigcup \, \{ Q A^* : Q$ is a prefix code and 
$Q \equiv_{\sf bd} P\} \big)$ $ \ \equiv_{\sf bd} \ P$ .
\end{pro}
{\bf Proof.} The first intersection is ${\sf ends}(P)$ since
${\sf ends}(Q) = {\sf ends}(P)$ when $Q \equiv_{\sf bd} P$; so this result
is different than the corresponding result for $\equiv_{\sf end}$ (in
Prop.\ \ref{topolInterior}).
For the second intersection result this is similar to the
proof of Prop.\ \ref{topolInterior}. 

For the union of the $Q A^*$ we have by Prop.\ \ref{topolInterior}, 
$\, {\sf prefC}\big(\bigcup_{Q \equiv_{\sf bd} P} Q A^* \big)$ 
$ \, \equiv_{\sf end} P$. Also, by Prop.\ \ref{equiv_boundEndsPROP},
${\sf ends}(Q A^*) = {\sf ends}(P)$ for every prefix code $Q$ such that   
$Q \equiv_{\sf bd} P$; hence, 
${\sf ends}\big( \bigcup_{Q \equiv_{\sf bd} P} Q A^* \big) = {\sf ends}(P)$.
Then the result follows by Prop.\ \ref{equiv_boundEndsPROP}.
\ \ \ $\Box$

\begin{defn} \label{boundequiv_RM}
{\bf (bounded end-equivalence of functions).}
 \ Two right-ideal morphisms $f, g$ are {\em boundedly end-equivalent} 
(denoted by $f \equiv_{\sf bd} g$) \ iff
 \ ${\sf domC}(f) \equiv_{\sf bd} {\sf domC}(g)$, and $f(x) = g(x)$
for all $x \in {\sf Dom}(f) \cap {\sf Dom}(g)$.
Equivalently, $f \equiv_{\sf bd} g$ \ iff \ $f \equiv_{\sf end} g$
 \ and \ ${\sf domC}(f)$ $\equiv_{\sf bd}$ ${\sf domC}(g)$.

For any family ${\cal T}$ of total functions as in Def.\ 
\ref{equiv_boundDEF} we define: $f \equiv_{\cal T} g$ iff 
$f \equiv_{\sf end} g$ and ${\sf domC}(f) \equiv_{\cal T} {\sf domC}(g)$. 
\end{defn}  
{\bf Notation:} When $f \in {\cal RM}^{\sf P}$, the 
$\, \equiv_{\sf bd}$-class in ${\cal RM}^{\sf P}$ of $f$ is denoted by 
$[f]_{\sf bd}$. So, $[f]_{\sf bd} = $
$\{g \in {\cal RM}^{\sf P}: g \equiv_{\sf bd} f\}$.
More generally, $[f]_{\cal T}$ denotes the $\, \equiv_{\cal T}$-class in 
${\cal RM}^{\sf P}$ of $f$.  Note that we define $[f]_{\sf bd}$ and 
$[f]_{\cal T}$ to only contain elements of ${\cal RM}^{\sf P}$. 

\medskip

For the rest of this Section we study $\equiv_{\sf bd}$.  The relations 
$\equiv_{\sf poly}$ and $\equiv_{\sf E3}$ will be investigated 
in the next Section.

\begin{pro} \label{bdequivUnionIntersLatt}
 \ (1) Let $P_1, P_2 \subset A^*$ be prefix codes such that
$P_1 \equiv_{\sf bd} P_2$, let $P_{\cap}$ be the prefix code that
generates the right ideal $P_1 A^* \cap P_2 A^*$, and let $P_{\cup}$ be the
prefix code that generates the right ideal $P_1 A^* \cup P_2 A^*$.  Then
$P_1 \equiv_{\sf bd} P_2 \equiv_{\sf bd} P_{\cap}$ \ $\equiv_{\sf bd}$ $P_{\cup}$.

\noindent
(2) Let $f_1, f_2$ be right-ideal morphisms such that
$f_1 \equiv_{\sf bd} f_2$. Then $f_1 \cap f_2$ and $f_1 \cup f_2$ are
right-ideal morphisms, and $f_1 \equiv_{\sf bd} f_2$ \
$\equiv_{\sf bd}$ \ $f_1 \cap f_2 \ \equiv_{\sf bd} \ f_1 \cup f_2$.
\end{pro}
{\bf Proof.}
(1) This follows from Prop.\ \ref{equiv_boundEndsPROP}, since
${\sf ends}(P_1 A^* \cap P_2 A^*) = {\sf ends}(P_1) \cap {\sf ends}(P_2)$
and ${\sf ends}(P_1 A^* \cup P_2 A^*)$ $=$ 
${\sf ends}(P_1) \cup {\sf ends}(P_2)$. 

\noindent For (2) the proof is similar to the proof of Prop.\
\ref{endequivUnionIntersLatt}.  
 \ \ \ $\Box$

\bigskip

\noindent Just as for $\equiv_{\sf end}$ (see Def.\ \ref{MaxExtMorph}), we 
define a maximum extension within a $\, \equiv_{\sf bd}$-class or a 
$\, \equiv_{\cal T}$-class.

\begin{defn} \label{MaxExtBoundMorph}
 \ For any right-ideal morphism $f$: $A^* \to A^*$ we define

\smallskip

 \ \ \ $f_{\sf b, max} \ = \ \bigcup \, \{ g : g$ is a right-ideal morphism
       with $g \equiv_{\sf bd} f \}$.

\smallskip

\noindent For a family $\cal T$ of functions as in Def.\ 
\ref{equiv_boundDEF} we define

\smallskip

 \ \ \ $f_{{\cal T}, {\sf max}} \ = \ \bigcup \, \{ g : g$ is a right-ideal
 morphism with $g \equiv_{\cal T} f \}$. 
\end{defn}
Then, just as in Prop.\ \ref{MaxExtPROP}, we have:

\begin{pro} \label{MaxExtBoundPROP} \hspace{-.13in} {\bf .}

\noindent (1) For every right-ideal morphism $f$, $f_{\sf b, max}$ is a
function, and a right-ideal morphism $A^* \to A^*$.
Moreover, $f \equiv_{\sf bd} f_{\sf b, max}$.

\smallskip

\noindent (2) For any right-ideal morphisms $f,g$ we have:
 \ \ $g \equiv_{\sf bd} f$ \ iff \ $g_{\sf b, max} = f_{\sf b, max}$.
\end{pro}
{\bf Proof.}  The same proof as for Prop.\ \ref{MaxExtPROP} works here 
(using Prop.\ \ref{bdequivUnionIntersLatt} and Prop.\ 
\ref{infUnionInterBd}).
 \ \ \ $\Box$

\bigskip

\noindent Recall the action of a right-ideal morphism $f$: $A^* \to A^*$ on 
$A^{\omega}$: For any $p \in {\sf domC}(f)$ and $w \in A^{\omega}$, we 
define $f(pw) = f(p) \, w$. The domain of the action of $f$ on $A^{\omega}$ 
is $\, {\sf domC}(f) \ A^{\omega}$.
Accordingly, ${\cal RM}^{\sf P}$ acts on $A^{\omega}$ (non-faithfully).
We have the following characterization of $\, \equiv_{\sf bd} \, $ in terms 
of the Cantor space: 

\begin{cor} \label{end_boundActCOR} 
 \ Two right-deal morphisms $f, g: A^* \to A^*$ have the same action on
$A^{\omega}$ \ iff \ $g \equiv_{\sf bd} f$.
Hence on $A^{\omega}$:
 \ $g_{\sf b, max} \circ f_{\sf b, max} = (g \circ f)_{\sf b, max}$.

\smallskip

The relation $\, \equiv_{\sf bd}$ is a {\em congruence} on the monoid of all
right-ideal morphisms of $A^*$, and in particular on ${\cal RM}^{\sf P}$.
\end{cor} 
{\bf Proof.} If $g \equiv_{\sf bd} f$ then 
${\sf domC}(g) \equiv_{\sf bd} {\sf domC}(f)$, hence by Prop.\
\ref{equiv_boundEndsPROP}:  
${\sf ends}({\sf domC}(g)) = {\sf ends}({\sf domC}(f))$. 
So the actions of $f$ and $g$ on $A^{\omega}$ have the same domain. Since 
$g \equiv_{\sf bd} f$, the functions $f$ and $g$ agree on their common
domain in $A^{\omega}$, so they have the same action on $A^{\omega}$.

Conversely, if $f$ and $g$ act in the same way on $A^{\omega}$ then 
${\sf ends}({\sf domC}(g)) = {\sf ends}({\sf domC}(f))$, so 
${\sf domC}(g) \equiv_{\sf bd} {\sf domC}(f)$ (by Prop.\
\ref{equiv_boundEndsPROP}). Also, if $f$ and $g$ act in the same way on
$A^{\omega}$ then for all $x \in {\sf Dom}(f) \cap {\sf Dom}(g)$ and 
for all $w \in A^{\omega}$: $\, f(xw) = g(xw)$. Hence (since $x \in$
${\sf Dom}(f) \cap {\sf Dom}(g)$), $\, f(x) = g(x)$. So, $f$ and $g$ 
agree on ${\sf Dom}(f) \cap {\sf Dom}(g)$, hence $f \equiv_{\sf bd} g$.

The rest of the corollary follows now.
 \ \ \ $\Box$

\begin{pro} \label{bMaxRec}
 \ There exists $g \in {\cal RM}^{\sf rec}$ such that 
${\sf Dom}(g_{\sf b, max})$ is undecidable; so 
$g_{\sf b, max} \not\in {\cal RM}^{\sf rec}$. 
Moreover, this function $g$ can be chosen so that in addition we have 
$\, g_{\sf b, max} = g_{\sf e, max} $.
\end{pro}
{\bf Proof.} Let $L \subset 0^*$ (over the one-letter alphabet $\{0\}$) be 
an r.e.\ language that is undecidable. We assume in addition that for all
$0^i, 0^j \in L$ we have $|i -j| > 2$ (if $i \neq j$).
Let $M$ be a deterministic Turing machine that accepts $L$.
Let $T(0^n)$ be the running time of $M$ on input $0^n \in L$; $\, T(w)$ is
undefined for $w \in \{0,1\}^* - L$. 
We define the right-ideal morphism $g$ by

\smallskip

 \ \ \ ${\sf domC}(g) \ = \ $
$\bigcup_{0^n \in L} \ 0^n \ \{0,1\} \ 1 \ \{0,1\}^{T(0^n)}$,

\smallskip

 \ \ \ $g(0^n a 1 z) \ = \ 0^n \ \ov{a} \ 1 z$, 

\smallskip

\noindent if $0^n \in L$, $a \in \{0,1\}$, and 
$z \in \{0,1\}^{T(0^n)} \, \{0,1\}^*$.
Here $\ov{a}$ denotes the complement of $a$ (i.e., $\ov{0} = 1$, 
$\ov{1} = 0$). The set ${\sf domC}(g)$ is a prefix code because of the
$|i -j| > 2$ condition on $L$.

Membership in ${\sf Dom}(g)$ is decidable in linear time:
for an input $0^n a 1 z$ with $a \in \{0,1\}$ and $z \in \{0,1\}^*$, it
suffices to run the machine $M$ for $\le |z|$ steps on $0^n$.
And $|g(0^n a 1 z)| = |0^n a 1 z|$, hence $g \in {\cal RM}^{\sf P}$.
 
We consider the following right-ideal morphism $h$, which extends $g$:

\smallskip

 \ \ \ ${\sf domC}(h) \ = \ L \ \{0,1\} \ 1$, 

\smallskip

 \ \ \ $h(0^n a 1) \ = \ 0^n \ \ov{a} \ 1$, 

\smallskip

\noindent if $0^n \in L$, $a \in \{0,1\}$.

We claim that $h = g_{\sf b, max} = g_{\sf e, max}$. Indeed, we have
$h \equiv_{\sf bd} g$, since ${\sf ends}({\sf domC}(g))$ $=$
$L \, \{0,1\} \, 1 \, \{0,1\}^{\omega}$ $=$ ${\sf ends}({\sf domC}(h)))$
(using Prop.\ \ref{equiv_boundEndsPROP}).
And $h$ cannot be further extended to a right-ideal morphism that is 
$\, \equiv_{\sf bd} h$ (because $h$ permutes $0^n 0 1, \, 0^n 11$).

Also, ${\sf cl}({\sf ends}({\sf domC}(h))) \, = \, $
$L \, \{0,1\} \, 1 \, \{0,1\}^{\omega} \ \cup \ \{ 0^{\omega}\}$.
But $h$ cannot be extended to any prefix $0^i$ ($i \in {\mathbb N}$) 
of $0^{\omega}$ (again because $h$ permutes $0^n 0 1, \, 0^n 11$).
So, $h$ is also the maximal $\equiv_{\sf end}$-equivalent extension of $g$.
 \ \ \ $\Box$

\begin{lem} \label{pointwiseExt}
 \ For every right-ideal morphisms $f$ and $g$ such that $f \subseteq g$
and $f \equiv_{\sf bd} g$ we have: 

\smallskip

\noindent (1) \ For all $x \in {\sf domC}(g)$:
 \ $xA^* \cap {\sf domC}(f)$ is finite.

\smallskip

\noindent (2) \ For any $x \in {\sf domC}(g)$, $f$ can be 
extended to a right-ideal morphism whose domain code is 
$\, ({\sf domC}(f) - xA^*) \cup \{x\}$ \ (this is called a {\em one-point
extension} of $f$).

 \ If $f \in {\cal RM}^{\sf P}$ then this extension is also in 
${\cal RM}^{\sf P}$.

\smallskip

\noindent (3) \ For any {\em finite} subset $C \subset {\sf domC}(g)$, 
$f$ can be extended to a right-ideal morphism whose domain code includes 
$C$ \ (this is called a {\em finite extension} of $f$).

This extension belongs to ${\cal RM}^{\sf P}$ if $f \in {\cal RM}^{\sf P}$.

\smallskip

\noindent (4) \ Items (1), (2), (3) hold in particular when 
$g = f_{\sf b, max}$. 

\smallskip

\noindent (5) \ There exist $f \in {\cal RM}^{\sf P}$ such that for some 
$x \in {\sf domC}(f_{\sf e, max})$ the set $xA^* \cap {\sf domC}(f)$ is 
infinite.  
\end{lem}
{\bf Proof.} (1) Let $\beta(.)$ be the bounding function that corresponds to
$f \equiv_{\sf bd} g$. Then the tree  

\smallskip

 \ \ \ \ \     
$T_{x,f} \ = \ \{z \in A^*: x$ is a prefix of $z$ and $z$ is a prefix of 
       some word in ${\sf domC}(f)\}$ 

\smallskip

\noindent has $x$ as root and $xA^* \cap {\sf domC}(f)$ as set of leaves, 
and has depth $\le \beta(|x|)$. Moreover, the degree of each vertex is 
$\le 2$. Hence the tree $T_{x,f}$ and its set of leaves 
$xA^* \cap {\sf domC}(f)$ are finite.

\noindent (2) Since $f \equiv_{\sf bd} g$ and $f \subseteq g$, every end that 
starts at $x$ intersects ${\sf domC}(f)$. Since $g$ and $f$ agree on 
${\sf domC}(f)$, $f$ can be extended from ${\sf Dom}(f)$ to 
$(xA^* \cap {\sf domC}(f)) \cup {\sf Dom}(f)$. For domain codes, the effect
of this extension is to replace $xA^* \cap {\sf domC}(f)$ by $\{x\}$.

If $f \in {\cal RM}^{\sf P}$ then the one-point extension is also in 
${\cal RM}^{\sf P}$, since $xA^* \cap {\sf domC}(f)$ is finite.

\noindent (3) follows from (2). 

\noindent (4) Recall that $f \equiv_{\sf bd} f_{\sf b, max}$, by Prop.\ 
\ref{MaxExtBoundPROP}.(1); hence the Lemma applies to the maximum 
$\equiv_{\sf bd}$-equivalent extension $f_{\sf b, max}$. 

\noindent (5)
For $\equiv_{\sf end}$ and $f_{\sf e, max}$ the situation is different. 
E.g., when $f = {\bf 1}|_{0^* 1}$, we have $f_{\sf e, max} = {\bf 1}$ and 
${\sf domC}(f_{\sf e, max}) = \{\varepsilon\}$. Then for $x = \varepsilon$ 
we have $x \, \{0,1\}^* \cap {\sf domC}(f) = 0^* 1$, which is infinite.
 \ \ \ $\Box$

\medskip

\noindent The precise definitions of one-point and finite extensions of a
right-ideal morphism are as follows.

\begin{defn} \label{onepointExt}
 \ Let $f, g$ be right-ideal morphisms. We call $g$ a {\em one-point
extension} of $f$ \ iff
 \ {\rm (1)} $g$ is an extension of $f$,
 \ {\rm (2)} $g \equiv_{\sf bd} f$,
 \ {\rm (3)} there exists $x_0 \in {\sf domC}(g)$ such that
$ \, {\sf Dom}(g) = x_0 A^* \, \cup \, {\sf Dom}(f)$.

We call $g$ a {\em finite extension} of $f$ \ iff
 \ {\rm (1)} $f \subseteq g$, and
{\rm (2)} $g \equiv_{\sf bd} f$,  as above, and {\rm (3)} there
exists a finite prefix code $F \subseteq {\sf domC}(g)$
such that $ \, {\sf Dom}(g) = F A^* \, \cup \, {\sf Dom}(f)$.
\end{defn}
A one-point extension is a special case of a finite extension (when
$F = \{x_0\}$), and a finite extension can be constructed by a finite
sequence of one-point extensions.

\smallskip

It follows from the definition that for a finite extension,
$ \, {\sf domC}(g) =  F \, \cup \, ({\sf domC}(f) - F A^*)$.

Moreover, ${\sf domC}(f) \, \cap \, F A^* \, $ is  finite (by Lemma
\ref{pointwiseExt}(1)). Since
${\sf domC}(g)$ $=$ $F \, \cup \, $
$({\sf domC}(f) - ({\sf domC}(f) \, \cap \, F A^*))$, we
conclude that the symmetric difference
${\sf domC}(g) \vartriangle {\sf domC}(f)$ is finite.
Conversely, suppose that $f \subseteq g$,
$g \equiv_{\sf bd} f$, and
${\sf domC}(g) \vartriangle {\sf domC}(f)$ is finite; then
$g$ is a finite extension of $f$ (in the sense of Def.\
\ref{onepointExt}). Indeed, $F = {\sf domC}(g) - {\sf domC}(f)$ is
finite, and satisfies
$ \, {\sf Dom}(g) = F A^* \, \cup \, {\sf Dom}(f)$.
Thus, {\em for right ideal morphisms $g, f$ that satisfy
$f \subseteq g$ and $g \equiv_{\sf bd} f$ we have:
$g$ is a finite extension of $f \, $ iff
$\, {\sf domC}(g) \vartriangle {\sf domC}(f)$ is finite.}

\medskip

By Lemma \ref{pointwiseExt}, $f_{\sf b, max}$ can be constructed
from $f$ by an $\omega$-sequence of one-point extensions.
On the other hand, the example $f = {\sf id}|_{0^*1}$ shows that a
right-ideal morphism $f$ might be extendable
(to ${\sf id} \equiv_{\sf end} f$ in this example),
without having any finite extension. So, in general,
$f_{\sf e, max}$ cannot be obtained from $f$ by an $\omega$-sequence of
finite extensions.

\medskip

\noindent The following consequence of the Lemma is a little surprising.

\begin{pro} \label{Phi_is_poly}
 \ Let $\cal T$ be any family of functions as in Def.\ \ref{equiv_boundDEF} 
such that, in addition, ${\sf poly} \subseteq {\cal T}$.
Then for every right-ideal morphism $f$,
 \ $f_{{\cal T}, {\sf max}} = f_{\sf b, max}$.

In particular, $f_{\sf poly, max} = f_{\sf b, max}$.
\end{pro}
{\bf Proof.} Since $\equiv_{\cal T}$ implies $\equiv_{\sf bd}$, 
$f_{{\cal T}, {\sf max}} \subseteq f_{\sf b, max}$.

On the other hand let $x \in {\sf domC}(f_{\sf b, max})$.
Since $f \equiv_{\sf bd} f_{{\cal T}, {\sf max}} \equiv_{\sf bd}$
$f_{\sf b, max}$, if we let $F = f_{{\cal T}, {\sf max}}$ we have
$f_{\sf b,max} = F_{\sf b, max}$.  For every
$x \in {\sf domC}(F_{\sf b, max})$, $x A^* \cap {\sf domC}(F)$ is finite
(by Lemma \ref{pointwiseExt}); and $F$ can be extended to a
right-ideal morphism $F_0$ which is defined on $x$. Since
$x A^* \cap {\sf domC}(F)$ is finite, $F_0 \equiv_{\cal T} F$. But since
$F$ is already $\equiv_{\cal T}$-maximum, the extension $F_0$ of $F$ is
$F$ itself. So, $F(x) = F_0(x) = f_{\sf b, max}(x)$. Since this holds
for every $x \in {\sf domC}(f_{\sf b, max})$, it follows that
$F = f_{\sf b, max}$.
 \ \ \ $\Box$

\begin{pro} \label{finExtThompson}
 \ For a right-ideal morphism $f$, the following are equivalent:

\noindent (1) \ $f$ is finitely extendable, to a strictly larger domain;

\noindent (2) \ there exist $x_0, y_0 \in A^*$ such that 
$(x_0 0, \, y_0 0), \, (x_0 1, \, y_0 1) \in f$, and 
$(x_0, y_0) \not\in f$; 

\noindent (3) \ $f \, \neq \, f_{\sf b, max}$.
\end{pro}
{\bf Proof.} The implication (2)$\Rightarrow$(1) is clear, since 
$(x_0 0, \, y_0 0), \, (x_0 1, \, y_0 1) \in f$ implies that $f$ can be
extended to $f \cup \{(x_0, y_0)\}$.   And (1) implies that 
$\, f \, \subsetneqq \, f \cup \{(x_0, y_0)\} \, \subseteq f_{\sf b, max}$, 
so (1) implies (3).

Let us prove that that (3) implies (2). By (3) there exists
$z \in {\sf domC}(f_{\sf b, max})$ with $z \not\in {\sf Dom}(f)$.  
Consider the rooted tree with root $z$ and vertex set and edges set
respectively 
 
\smallskip

 \ \ \ $V = \{z w \, : \, w \in A^*$ and $zw$ is a prefix of a word in
$z A^* \cap {\sf domC}(f)\}$, 

\smallskip

 \ \ \ $E = \{(v, va) \in V \times V \, : \, a \in A\}$. 

\smallskip

\noindent This is a binary tree (every vertex has $\le 2$ children), and it 
is {\em saturated} (i.e., every vertex has either 2 children or none). 
The set of leaves is $z A^* \cap {\sf domC}(f)$, and this set is finite 
(by Lemma \ref{pointwiseExt}); hence since the tree is saturated, it is 
finite. Also, since $z \not\in {\sf Dom}(f)$ and since the tree is saturated, 
the tree has at least 3 vertices. 
Let $d$ be the depth of the tree (number of edges in a longest path from the 
root). Let $x_0$ be any non-leaf vertex at distance $d - 1$ from the root; 
since $x_0$ is not a leaf, it has two children, namely $x_0 0$ and $x_0 1$. 
Then $x_0 0$ and $x_0 1$ are at distance $d$ from the root, so they are 
leaves, i.e., $x_0 0, x_0 1 \in z A^* \cap {\sf domC}(f)$.
Let $y_0 = f_{\sf b, max}(x_0)$; then for $a \in \{0,1\}$ we have
$y_0 a = f_{\sf b, max}(x_0 a)$, and the latter is equal to $f(x_0 a)$ 
(since $x_0 0, x_0 1 \in {\sf domC}(f)$).
Hence, $(x_0, y_0)$ satisfies (2). 
 \ \ \ $\Box$

\medskip

\noindent We mention the following, which will however not be used in
this paper:

\medskip

\noindent {\bf Fact.} {\it 
For every $f \in {\sf fP}$, the encoded right-ideal morphism
$f^C$ is not finitely extendable, and not infinitely
extendable. In other words, $f^C =$
$(f^C)_{\sf b, max} = (f^C)_{\sf e, max}$.
}

\smallskip

\noindent
{\bf Proof.} The encoding $f^C$ was defined in the Introduction. It follows
from that definition that ${\sf domC}(f^C) \subseteq \{0,1\}^* \, 11$, so we 
never have $z0, z1 \in {\sf domC}(f^C)$ for any $z \in A^*$. Hence, by Prop.\
\ref{finExtThompson}, $f^C$ is not finitely extendable, so
$f^C = (f^C)_{\sf b, max}$.

Proof that $f^C$ is not infinitely extendable: Let $P \subset A^*$ be any
prefix code such that ${\sf domC}(f^C) \subseteq P A^*$ and
${\sf domC}(f^C) \equiv_{\sf end} P$; we want to show that
$P = {\sf domC}(f^C)$, which means that $f^C$ cannot be extended to a larger
domain.  Let us abbreviate ${\sf domC}(f^C)$ by $D$.

Since $P \equiv_{\sf end} D$ we have (by Lemma \ref{equivtoEndEquiv}): 
$(\forall p \in P)(\exists x11 \in D)[ \, p \parallel_{\sf pref} x11 \, ]$;
moreover, since $D \subset PA^*$, this $p$ is a prefix of $x11$. 
Hence, since $D \subset \{00, 01\}^* \, 11$, we conclude that 
$p \in \{00, 01\}^* \cdot \{ \varepsilon, 1, 11\}$; i.e., for each 
$p \in P$ we have three possibilities: $p \in \{00, 01\}^*$, \,
$p \in \{00, 01\}^* \, 1$, \, $p \in \{00, 01\}^* \, 11$.

\smallskip

\noindent {\sf Claim 1.} If $p \in \{00, 01\}^* \, 11$ then $p \in D$.

\smallskip

\noindent {\sf Proof.} Since $p$ is a prefix of a word $x11 \in D$, and
$p \in \{00, 01\}^* \, 11$, we conclude $p = x11$ (since $\{00, 01\}^* \, 11$
is a prefix code). \ \ \ {\sf [End, proof of Claim 1.]}

\smallskip

\noindent {\sf Claim 2.} \ $P \cap \{00, 01\}^* \, 1 = \varnothing$; 
 \ i.e., $p$ cannot be in $\{00, 01\}^* \, 1$.

\smallskip

\noindent {\sf Proof.}  If there exists $p \in P \cap \{00, 01\}^*\, 1$
and $p$ is a prefix of some $x 11 \in \{00, 01\}^* \, 11$, then
$p = x 1$. But then $P \not\equiv_{\sf end} D$, since the right ideal
$p0 A^* = x 10 A^*$ $\subseteq \{00, 01\}^* \, 10 \, A^*$
intersects $PA^*$ but not $\{00, 01\}^* \, 11 \, A^*$.
 \ \ \ {\sf [End, proof of Claim 2.]}

\smallskip

\noindent {\sf Claim 3.} \ $P \cap \{00, 01\}^* = \varnothing$; \ i.e.,
$p$ cannot be in $\{00, 01\}^*$.

\smallskip

\noindent {\sf Proof.} Assume there exists $p \in P \cap \{00, 01\}^*$
such that $p$ is a prefix of some $x 11 \in \{0,1\}^* \, 11$. But then
$P \not\equiv_{\sf end} D$, since the right ideal
$p 10 A^* = x 10 A^* \subseteq \{00, 01\}^* \, 10 \, A^*$ intersects $PA^*$
but not $\{00, 01\}^* \, 11 \, A^*$.
 \ \ \ {\sf [End, proof of Claim 3.]}

\smallskip

We are left with only case 1, i.e., $P \subseteq D$. Hence $P = D$, since
$D \equiv_{\sf end} P$.
 \ \ \ $\Box$

\bigskip

\noindent {\bf Notation:} The quotient monoid of ${\cal RM}^{\sf P}$ under
the congruence $\equiv_{\sf bd}$ of ${\cal RM}^{\sf P}$ is denoted by
${\cal RM}^{\sf P}\!/\! \equiv_{\sf bd}$; this is also the quotient monoid
for the action of ${\cal RM}^{\sf P}$ on $A^{\omega}$.

Moreover, ${\cal RM}^{\sf P}\!/\! \equiv_{\sf bd}$ will also be denoted by
${\cal M}_{\sf bd}^{\sf P}$.

\medskip

Recall that ${\cal RM}^{\sf rec}$ consists of all right-ideal morphisms 
that are partial recursive with decidable domain, and that have a total 
recursive input-output balance.

\begin{pro} \label{RMRecPolyRMPpolybd}
The monoid  ${\cal RM}^{\sf rec}\!/\!\equiv_{\sf bd}$ is isomorphic to
${\cal RM}^{\sf P}\!/\!\equiv_{\sf bd}$ \ ($= {\cal M}_{\sf bd}^{\sf P}$).
\end{pro}
{\bf Proof.} This is proved in the same way as Prop.\ 
\ref{RMRecPolyRMPpoly}.
 \ \ \ $\Box$

\begin{lem} \label{bdequivFinite} 
 \ If $P, Q$ are prefix codes such that $P \equiv_{\sf bd} Q$,
and if $P$ is finite, then $Q$ is finite.
\end{lem}
{\bf Proof.} Let $\beta$ be the length-bounding (total) function associated 
with $P \equiv_{\sf bd} Q$. Since $P \equiv_{\sf bd} Q$, every element 
$x_2 \in Q$ is prefix-comparable to some $x_1 \in P$. Since $P$ is finite, 
the elements of $P$ can have only finitely many prefixes, hence the set of 
elements of $Q$ that are prefixes of some element(s) of $P$ is finite. 
Moreover, each $x_1 \in P$ can be the prefix of only finitely many 
$x_2 \in Q$, since every such $x_2$ has length $\le \beta(|x_1|)$. Since $P$ 
is finite, the set of length bounds $\{\beta(|x_1|) : x_1 \in P\}$ is 
finite.  Hence $Q$ is finite. 
 \ \ \ $\Box$

\bigskip

The next theorem refers to the well-known Richard Thompson group $V$
(a.k.a.\ $G_{2,1}$). In order to make the paper self-contained we define
$V$ next.  First, let

\medskip

${\sf riAut}^{\sf fin} \ = \ \{f \; : \, f$ is a right-ideal morphism of
$A^*$, such that

\hspace{1.1in} (1) \ $f$ is injective,

\hspace{1.1in} (2) \ ${\sf domC}(f)$ and ${\sf imC}(f)$ are maximal
  prefix codes,

\hspace{1.1in} (3) \ ${\sf domC}(f)$ (and hence ${\sf imC}(f)$) is finite\}.

\medskip

\noindent The notation ``${\sf riAut}^{\sf fin}$'' stands for right-ideal 
automorphism with finite domain code.
Every element of ${\sf riAut}^{\sf fin}$ can be given by a
bijection between two finite maximal prefix codes, and it is
straightforward to prove that ${\sf riAut}^{\sf fin}$ is a submonoid of
${\cal RM}^{\sf P}$.
For every finite maximal prefix code $P$, ${\sf id}_{PA^*}$ is an idempotent
of ${\sf riAut}^{\sf fin}$, hence ${\sf riAut}^{\sf fin}$ is not a group.
The Thompson group $V$ is defined by

\medskip

$V \ = \ {\sf riAut}^{\sf fin}/\!\!\equiv_{\sf bd}$
 \ \ ($\le {\cal M}_{\sf bd}^{\sf P}$).

\medskip

\noindent See \cite{BiThomps} for details, and a proof that this is a group;
${\sf riAut}^{\sf fin}$ is studied in \cite{BiMonMap}.
This group has remarkable properties (e.g., it is finitely presented and
simple), and can be defined in several ways. It was introduced by Richard
J.\ Thompson in the 1960s along with two other remarkable groups; see
\cite{CFP} for more background.

\begin{thm} \label{boundReg} \hspace{-.13in} {\bf .}

\noindent (1) \ The monoid ${\cal M}_{\sf bd}^{\sf P}$ is regular, 
${\cal J}^0$-simple, and finitely generated.
Every element of ${\cal RM}^{\sf P}$ is $\equiv_{\sf bd}$-equivalent
to a regular element of ${\cal RM}^{\sf P}$.

\smallskip

\noindent (2) \ For any family $\cal T$ of total functions as in Def.\
\ref{equiv_boundDEF}, the group of units of ${\cal M}_{\cal T}^{\sf P}$
is the Richard Thompson group  $V$ (a.k.a.\ $G_{2,1}$).
In particular, the group of units of ${\cal M}_{\sf bd}^{\sf P}$ is $V$.

\smallskip

\noindent (3) \ ${\cal M}_{\cal T}^{\sf P}$, and in particular
${\cal M}_{\sf bd}^{\sf P}$,  is not congruence-simple. 
\end{thm}
{\bf Proof.} (1) Regularity and finite generation are proved in the same
way as for ${\cal M}_{\sf end}^{\sf P}$ (Theorem \ref{M_is_reg}, and Lemma
\ref{recTimeConstr}). Since ${\cal RM}^{\sf P}$ is ${\cal J}^0$-simple, so
is its homomorphic image ${\cal RM}^{\sf P}\!/\!\equiv_{\sf bd}$
($= {\cal M}_{\sf bd}^{\sf P}$).

\smallskip

\noindent (2) By Lemma \ref{bdequivFinite} (and since $\equiv_{\cal T}$
implies $\equiv_{\sf bd}$), if $f \in {\cal RM}^{\sf P}$
has a finite ${\sf domC}(f)$ then every element of $[f]_{\cal T}$ has a
finite domain code. Hence when ${\sf domC}(f)$ is finite,
$[f]_{\cal T} = [f]_{\sf bd}$.
It now follows immediately from the definition so every element of $V$
belongs to ${\cal M}_{\cal T}^{\sf P}$. Thus, $V$ is a subgroup of the
group of units.

Conversely, suppose $[F]_{\cal T} \equiv_{\cal H} [{\sf id}]_{\cal T}$ in
${\cal M}_{\cal T}^{\sf P}$, where $F \in {\cal RM}^{\sf P}$.
We note first that if  $e \in {\cal RM}^{\sf P}$ satisfies
$e \equiv_{\cal T} {\sf id}$, then $e = {\sf id}_{P A^*}$, for some finite
maximal prefix code $P \subset A^*$ (finiteness follows from Lemma
\ref{bdequivFinite}).
Now $[F]_{\cal T} \equiv_{\cal H} [{\sf id}]_{\cal T}$ implies that there
is a maximal prefix code $P_1 \subset A^*$ in {\sf P} such that
${\sf id}_{P_1 A^*} \equiv_{\cal T} {\sf id}$, and there exists
$g_2 \in {\cal RM}^{\sf P}$ such that
$g_2 \circ F = {\sf id}_{P_1 A^*}$ (since
 $[F]_{\cal T} \equiv_{\cal L} [{\sf id}]_{\cal T}$).
Then, ${\sf id}_{P_1 A^*} \equiv_{\cal T} {\sf id}$ implies
$F \circ {\sf id}_{P_1 A^*} \equiv_{\cal T} F$.
Letting $f = F \circ {\sf id}_{P_1 A^*}$ ($ \equiv_{\cal T} F$), we now have
$[f]_{\cal T} \equiv_{\cal H} [{\sf id}]_{\cal T}$ and
$g_2 \circ f = {\sf id}_{P_1 A^*}$ with
$f \circ g_2 \circ f = f$. So $f$ is regular.

We also have $f \circ g_1 = {\sf id}_{P_2 A^*}$ (since
  $[f]_{\cal T} \equiv_{\cal R} [{\sf id}]_{\cal T}$), where
$g_1 \in {\cal RM}^{\sf P}$ and $P_2$ is a finite maximal prefix code.
Since ${\sf id}_{P_2 A^*} \equiv_{\cal T} {\sf id}$ $\equiv_{\cal T}$
${\sf id}_{P_2 A^*}$, we have $P_2 \equiv_{\cal T}$
$\{\varepsilon\} \equiv_{\cal T} P_2$.  Since $f$ and
${\sf id}_{P_2 A^*}$ are regular, Lemma \ref{regLR}(1) implies
${\sf Im}(f) = P_2 A^*$.
Since $f$ and ${\sf id}_{P_1 A^*}$ are regular, Lemma \ref{regLR}(2) implies
that $f$ is injective and ${\sf Dom}(f) = P_1 A^*$.
Hence, $f$ is a bijection from $P_1 A^*$ onto $P_2 A^*$, with $P_1, P_2$
finite. Thus, $[f]_{\cal T}$ ($= [F]_{\cal T}$) belongs to the Thompson
group $V$.

\smallskip

\noindent
(3) Obviously, the congruence $\equiv_{\cal T}$ is a refinement of 
$\equiv_{\sf end}$, so there is a surjective homomorphism  
${\cal M}_{\cal T}^{\sf P} \ \twoheadrightarrow  \ $ 
${\cal M}_{\sf end}^{\sf P}$.  By (2) and Prop.\ 
\ref{Mgroupunits}, these two monoids have different groups of units, so they 
are not isomorphic, hence the above surjective homomorphism is not injective. 
So, $\equiv_{\cal T}$ is a strict refinement of $\equiv_{\sf end}$, so 
${\cal M}_{\cal T}^{\sf P}$ is not congruence-simple.
  \ \ \ $\Box$

\begin{pro} \label{RMpoly_Dclasses} \hspace{-.13in} {\bf .}

\noindent (1) \ In ${\cal RM}^{\sf P}$ we have: If 
 \ ${\bf 1} \equiv_{\cal D} f$, \ then ${\sf imC}(f)$ is finite.

\smallskip

\noindent (2) \ The monoid ${\cal RM}^{\sf P}$ is not ${\cal D}^0$-simple.
\end{pro}
{\bf Proof.} (1) If ${\bf 1} \equiv_{\cal D} f$ then $f$ is obviously 
regular (since {\bf 1} is regular, and the whole ${\cal D}$-class is 
regular if it contains a regular element). By definition of 
$\equiv_{\cal D}$, ${\bf 1} \equiv_{\cal D} f$ \ iff 
 \ ${\bf 1} \equiv_{\cal L} g \equiv_{\cal R} f$ for some 
$g \in {\cal RM}^{\sf P}$. By Lemma \ref{regLR} this implies that 
${\sf domC}(g)$ is finite (and, moreover, $g$ is injective, and 
${\sf domC}(g)$ is a maximal prefix code), and that 
${\sf Im}(g) = {\sf Im}(f)$; hence ${\sf imC}(g) = {\sf imC}(f)$. 
Since ${\sf domC}(g)$ is finite, ${\sf imC}(g)$ is finite. Hence 
${\sf imC}(f)$ ($ = {\sf imC}(g)$) is finite.

\smallskip

\noindent (2) Consider $f \in {\cal RM}^{\sf P}$ defined (for all
$n \ge 0$) by $f(0^{2n} 1) = 0^{2n+1} 1$ and $f(0^{2n+1} 1) = 0^{2n} 1$;
so, ${\sf domC}(f) = {\sf imC}(f) = 0^* 1$.  Then ${\sf imC}(f)$ is
infinite, so by (1), $f \not\equiv_{\cal D} {\sf id}$ in
${\cal RM}^{\sf P}$.
 \ \ \ $\Box$

\bigskip

Let $M_{2,1}$ denote the monoid generalization of the Thompson
group $V$ ($= G_{2,1}$).
To define $M_{2,1}$, consider first
$\, {\cal RM}^{\sf fin} \, = \, \{f \in {\cal RM}^{\sf P} : $
    ${\sf domC}(f)$ is finite\}.
Then, $M_{2,1} \, = \, {\cal RM}^{\sf fin}/\!\!\equiv_{\sf bd}$;
 \ see \cite{JCBmonThH}.

\begin{pro} \label{bdDiffDclasses} \hspace{-.13in} {\bf .}

\noindent (1) \ In ${\cal M}_{\sf bd}^{\sf P}$, the ${\cal D}$-class of the 
identity contains $M_{2,1} - \{{\bf 0}\}$. 

\smallskip

\noindent (2) \ In ${\cal M}_{\sf bd}^{\sf P}$, the ${\cal R}$-class of the 
identity contains some elements that are not in $M_{2,1}$.  
\end{pro} 
{\bf Proof.} (1) The monoid $M_{2,1}$ is the submonoid 
$\{ [f]_{\sf bd} \in {\cal M}_{\sf bd}^{\sf P}: {\sf domC}(f)$ is finite\}. 
In Theorem 2.5 of \cite{JCBmonThH} it was proved 
that $M_{2,1}$ is ${\cal D}^0$-simple. Therefore, $M_{2,1} - \{{\bf 0}\}$ is 
contained in the ${\cal D}$-class of {\bf 1} in ${\cal M}_{\sf bd}^{\sf P}$.

\noindent
(2) By Lemma 2.11 in \cite{s1f}, the ${\cal R}$-class of {\bf 1} in 
${\cal RM}^{\sf P}$ is 
$\{f \in {\cal RM}^{\sf P}: \varepsilon \in {\sf Im}(f)\}$.
Consider $f \in {\cal RM}^{\sf P}$ defined by $f(0^n 1) = 0^n$ for all 
$n \ge 0$; so, ${\sf domC}(f) = 0^* 1$ and ${\sf imC}(f) = \{\varepsilon\}$. 
Hence $f$ is in the ${\cal R}$-class of {\bf 1}. But $[f]_{\sf bd}$ 
($\in {\cal M}_{\sf bd}^{\sf P}$) does not belong to $M_{2,1}$, since the 
infinite prefix code $0^* 1$ is not $\, \equiv_{\sf bd}$-equivalent to a 
finite prefix code, by Lemma \ref{bdequivFinite}. 
 \ \ \ $\Box$

\begin{lem} \label{fingenessRI}
If $R_2 \subseteq R_1$ are right ideals of $A^*$ and $R_2$ is 
{\em essential}, then $R_1$ is also essential. 

If $R_2$ is {\em essential and finitely generated} (as a right ideal), then 
$R_1$ is finitely generated.
\end{lem}
{\bf Proof.} Every right ideal $x A^*$ intersects $R_2$, hence $x A^*$
obviously intersects $R_1$ (since $R_2 \subseteq R_1$). So $R_1$ is 
essential. 
If $R_2$ is a finitely generated essential right ideal then $R_2 = P_2 A^*$ 
for a finite maximal prefix code $P_2$.
It follows that $A^* - P_2 A^*$ is a finite set. Moreover, $R_1$ is 
generated by a subset of $P_2 \cup (R_1 - P_2 A^*)$. This is a subset of
$P_2 \cup (A^* - P_2 A^*)$, which is finite.
 \ \ \ $\Box$

\begin{lem} \label{bdequ_Im} {\rm (See also Lemma \ref{Im_endequiv}.)}
 \ If $f_1, f_2 \in {\cal RM}^{\sf P}$ satisfy $f_1 \equiv_{\sf bd} f_2$,
then ${\sf imC}(f_1) \equiv_{\sf bd} {\sf imC}(f_2)$.
\end{lem}
{\bf Proof.} By Prop.\ \ref{equiv_boundEndsPROP} and Cor.\
\ref{end_boundActCOR}: $\, f_1 \equiv_{\sf bd} f_2$ iff
${\sf domC}(f_1) \ A^{\omega} = {\sf domC}(f_2) \ A^{\omega}$ and 
$f_1, f_2$ agree on ${\sf Dom}(f_1) \cap {\sf Dom}(f_2)$.  Thus 
$f_1 \equiv_{\sf bd} f_2$ implies $\, {\sf imC}(f_1) \ A^{\omega}$ $=$ 
$f_1({\sf domC}(f_1) \ A^{\omega})$ $=$
$f_1({\sf domC}(f_1) \ A^{\omega} \cap {\sf domC}(f_2) \ A^{\omega})$ $=$
$f_2({\sf domC}(f_1) \ A^{\omega} \cap {\sf domC}(f_2) \ A^{\omega})$ $=$ 
$f_2({\sf domC}(f_2) \ A^{\omega})$ $=$ ${\sf imC}(f_2) \ A^{\omega}$. 
So ${\sf imC}(f_1) \ A^{\omega}$ $=$ ${\sf imC}(f_2) \ A^{\omega}$, hence 
(by Prop.\ \ref{equiv_boundEndsPROP} again),
${\sf imC}(f_1) \equiv_{\sf bd} {\sf imC}(f_2)$.
\ \ \ $\Box$

\bigskip

\noindent {\bf Notation:} For any right ideal $R \subseteq A^*$, the
(unique) prefix code that generates $R$ as a right ideal is denoted by
${\sf prefC}(R)$.

\begin{lem} \label{RMbdNotD0simple} \hspace{-0.13in} {\bf .} 

\noindent (1) Let $R \subseteq A^*$ be an essential right ideal such that 
$R \in {\sf P}$. If 
$\, [{\bf 1}]_{\sf bd} \, \equiv_{\cal D} \, [{\bf 1}_R]_{\sf bd} \, $
in ${\cal M}_{\sf bd}^{\sf P}$, then $R$ is finitely generated (as a right 
ideal).

\noindent (2) The monoid ${\cal M}_{\sf bd}^{\sf P}$ is not 
${\cal D}^0$-simple.
\end{lem}
{\bf Proof.} (1) We have $[{\bf 1}]_{\sf bd} \equiv_{\cal D} [f]_{\sf bd}$ 
 \ iff \ $[{\bf 1}]_{\sf bd} \equiv_{\cal L} [g]_{\sf bd} \equiv_{\cal R}$ 
$[f]_{\sf bd}$ for some $g \in {\cal RM}^{\sf P}$.
The relation $[{\bf 1}]_{\sf bd} \equiv_{\cal L} [g]_{\sf bd}$ is
equivalent to $[{\bf 1}]_{\sf bd} = [m]_{\sf bd} \, [g]_{\sf bd}$ for some
$m \in {\cal RM}^{\sf P}$, hence ${\bf 1} \equiv_{\sf bd} m g$.
So (by Lemma \ref{bdequivFinite}), ${\bf 1}_{P A^*} = m g$ for some 
finite maximal prefix code $P \subset A^*$.
It follows from ${\bf 1}_{P A^*} = m g$ that $P A^* \subseteq {\sf Dom}(g)$.
Since $P A^*$ is a finitely generated essential right ideal, it follows (by 
Lemma \ref{fingenessRI}) that ${\sf Dom}(g)$ is also a finitely generated 
essential right ideal.  
In summary, so far we have shown that
$\, [{\bf 1}]_{\sf bd} \equiv_{\cal L} [g]_{\sf bd} \equiv_{\cal R}$ 
$[f]_{\sf bd}$, where $g$ is such that 

\smallskip

 \ \ \ \ \ ${\sf domC}(g)$ is a finite maximal prefix code. 

\smallskip

We are interested in the case when $f = {\bf 1}_R$, where $R \subset A^*$ 
is an essential right ideal with 
$[{\bf 1}]_{\sf bd} \, \equiv_{\cal D} \, [{\bf 1}_R]_{\sf bd}$. 
Then $\, [{\bf 1}]_{\sf bd} \equiv_{\cal L} [g]_{\sf bd}$ $\equiv_{\cal R}$ 
$[f]_{\sf bd}$, with $g$ as above.
The relation $[g]_{\sf bd} \equiv_{\cal R} [{\bf 1}_R]_{\sf bd}$ implies 
that ${\bf 1}_R \equiv_{\sf bd} g h$ for some $h \in {\cal RM}^{\sf P}$. 
And this implies that ${\bf 1}_{R_1} = g h$ for some right ideal $R_1$ such 
that ${\sf prefC}(R_1) \equiv_{\sf bd} {\sf prefC}(R)$;  hence, $R_1$ is
essential (since $R$ is essential).  From ${\bf 1}_{R_1} = g h$ it follows 
that $R_1 \subseteq {\sf Im}(g)$; this implies that ${\sf Im}(g)$ is 
essential (by Lemma \ref{fingenessRI}). 
Moreover, since ${\sf Dom}(g)$ is finitely generated (as we saw above), 
$g({\sf Dom}(g)) = {\sf Im}(g)$ is finitely generated.  So, ${\sf Im}(g)$ 
is a finitely generated essential right ideal, i.e.,

\smallskip

 \ \ \ \ \ ${\sf imC}(g)$ is a finite maximal prefix code. 

\smallskip

The relation $[g]_{\sf bd} \equiv_{\cal R} [{\bf 1}_R]_{\sf bd}$
also implies that $g \equiv_{\sf bd} {\bf 1}_R \cdot k' \, $ for some 
$k' \in {\cal RM}^{\sf P}$; let $g_1' = {\bf 1}_R \cdot k'$. 

\smallskip

\noindent Let $g_1 = g_1' \cdot {\bf 1}_{{\sf Dom}(g)}$; hence
${\bf 1}_{{\sf Dom}(g_1)} = $
${\bf 1}_{{\sf Dom}(g_1')} \cdot {\bf 1}_{{\sf Dom}(g)}$ 
$ = {\bf 1}_{{\sf Dom}(g)} \cdot {\bf 1}_{{\sf Dom}(g_1')}$.
And let  $k = k' \cdot {\bf 1}_{{\sf Dom}(g_1)}$.
Multiplying $g \equiv_{\sf bd} g_1' = {\bf 1}_R \cdot k'$ on the
right by ${\bf 1}_{{\sf Dom}(g)}$ yields:
 
\smallskip

 \ \ \ \ \ $g \, \equiv_{\sf bd}$
$ \, g_1 \, = \, {\bf 1}_R \cdot k' \cdot {\bf 1}_{{\sf Dom}(g)}$,

\smallskip

\noindent and then multiplying this on the right by 
${\bf 1}_{{\sf Dom}(g_1')}$ yields

\smallskip

 \ \ \ \ \ ($g \, \equiv_{\sf bd}$)
$ \, g_1 \, = \, {\sf id}_R \cdot k'$
$\cdot$ ${\sf id}_{{\sf Dom}(g)} \cdot {\sf id}_{{\sf Dom}(g_1')}$
$ \, = \, {\sf id}_R \cdot k'$
$\cdot$ ${\sf id}_{{\sf Dom}(g_1')} \cdot {\sf id}_{{\sf Dom}(g)}$,

\smallskip

\noindent since ${\bf 1}_{{\sf Dom}(g_1')}$ and ${\bf 1}_{{\sf Dom}(g)}$
commute, and since ${\sf Dom}(g_1) \subseteq {\sf Dom}(g_1')$. Moreover,

\smallskip

 \ \ \ \ \ 
$k' \cdot {\bf 1}_{{\sf Dom}(g_1')} \cdot {\bf 1}_{{\sf Dom}(g)}$
$ \, = \, $
$k' \cdot {\bf 1}_{{\sf Dom}(g_1)} \, = k$.
  

\smallskip

\noindent Thus we have $\, g_1 = {\bf 1}_R \cdot k$; this implies 
${\sf Im}(g_1) \subseteq R$. Since $g \equiv_{\sf bd} g_1$, and 
${\sf Im}(g)$ is essential (as we saw above), Lemma \ref{bdequ_Im}
implies that ${\sf Im}(g_1)$ is also essential. And since 
${\sf Dom}(g)$ is finitely generated (as we saw above), and 
$g \equiv_{\sf bd} g_1$, it follows that ${\sf Dom}(g_1)$ is finitely
generated (by Lemma \ref{bdequivFinite}). Hence, ${\sf Im}(g_1) =$
$g_1({\sf Dom}(g_1))$ is finitely generated. So now we have
${\sf Im}(g_1) \subseteq R$ (seen above), where ${\sf Im}(g_1)$ is an 
essential right ideal that is finitely generated. By Lemma 
\ref{fingenessRI}, it follows that $R$ is finitely generated.

\smallskip

\noindent (2) Let $R$ be an essential right ideal in {\sf P} such that
${\sf prefC}(R)$ is infinite. Such right ideals exist; examples are 
$0^*1 A^*$ and $0^* 1 0^* 1 A^*$.  Then by (1),  
$\, [{\bf 1}]_{\sf bd} \, \not\equiv_{\cal D} \, [{\bf 1}_R]_{\sf bd}$.
 \ \ \ $\Box$

\medskip

\noindent We can now completely characterize the ${\cal D}$-relation in
${\cal M}_{\sf bd}^{\sf P}$:

\begin{thm} \label{twoDclasses}
 \ The monoid ${\cal M}_{\sf bd}^{\sf P}$ has exactly two non-zero 
${\cal D}$-classes, namely

\smallskip

$D_1 \ = \ \{[f]_{\sf bd} : \, [f]_{\sf bd}$ {\rm contains} $f \neq 0$ 
{\rm such that} ${\sf imC}(f)$ {\rm is finite}$\}$ \ \ and  
 
\smallskip

$D_2 \ = \ \{[f]_{\sf bd} : \, [f]_{\sf bd}$ {\rm contains} $f$ {\rm such 
that} ${\sf imC}(f)$ {\rm is infinite}$\}$.
\end{thm}
{\bf Proof.} The sets $D_1$ and $D_2$ are disjoint and they form a
bipartition of ${\cal M}_{\sf bd}^{\sf P} \, - \, \{[{\bf 0}]_{\sf bd}\}$.
Indeed, if $g \in [f]_{\sf bd}$ then
${\sf imC}(g) \equiv_{\sf bd} {\sf imC}(f)$ (by Lemma \ref{bdequ_Im}); and
finiteness of ${\sf imC}(f)$ implies finiteness of ${\sf imC}(g) \,$
(by Prop.\ \ref{bdequivFinite}). Thus, if $[f]_{\sf bd} \in D_1$ then all
elements of $[f]_{\sf bd}$ have finite image code, so
$[f]_{\sf bd} \not\in D_2$.
Recall that by definition,
$[f]_{\sf bd} \, = \, \{\varphi \in {\cal RM}^{\sf P}:$
$\, \varphi \equiv_{\sf bd} f\}$. \ Let us also define

\smallskip

 \ \ \  \ \ \  $[[f]]_{\sf bd} \, = \, \{\xi \in {\cal RM}^{\sf rec}: $
$\, \xi \equiv_{\sf bd} f\}$.

\smallskip

\noindent By Prop.\ \ref{RMRecPolyRMPpolybd}, ${\cal M}_{\sf bd}^{\sf P}$
($ \, = \, {\cal RM}^{\sf P}\!/\!\equiv_{\sf bd}$) is isomorphic to
${\cal RM}^{\sf rec}\!/\!\equiv_{\sf bd}$.
For the remainder of this proof we represent ${\cal M}_{\sf bd}^{\sf P}$
by ${\cal RM}^{\sf rec}\!/\!\equiv_{\sf bd}$.
Let us reformulate the description of the two non-zero ${\cal D}$-classes
in terms of ${\cal RM}^{\sf rec}$:

\smallskip

$\Delta_1 \ = \ \{\, [[f]]_{\sf bd} : \, [[f]]_{\sf bd}$ contains
$f \neq 0$ such that ${\sf imC}(f)$ is finite\},

\smallskip

$\Delta_2 \ = \ \{\, [[f]]_{\sf bd} : \, [[f]]_{\sf bd}$ contains $f$ such
that ${\sf imC}(f)$ is infinite\}.

\smallskip

\noindent The sets $\Delta_1$ and $\Delta_2$ are disjoint and form a
bipartition of
${\cal RM}^{\sf rec}\!/\!\equiv_{\sf bd} \, - \, \{[[{\bf 0}]]_{\sf bd}\}$;
this holds for the same reason as in the case of $D_1$ and $D_2$ above,
since Lemma \ref{bdequ_Im} and Prop.\ \ref{bdequivFinite} apply to all
right-ideal morphisms.

By Lemma \ref{RMbdNotD0simple}(2)), ${\cal M}_{\sf bd}^{\sf P}$ is not
${\cal D}^0$-simple, hence $\Delta_1 \cup \Delta_2$ is not one
${\cal D}$-class. So to prove the Theorem it suffices to prove that all
elements in $\Delta_1$ are ${\cal D}$-related, and all elements in
$\Delta_2$ are ${\cal D}$-related in ${\cal M}_{\sf bd}^{\sf P}$.

\smallskip

\noindent (1) Let us prove that all elements of $\Delta_1$ are
${\cal D}$-related.

Every element of $[[{\sf id}]]_{\sf bd}$ is of the form ${\sf id}_{PA^*}$,
where $P$ is a maximal prefix code, and by Lemma \ref{bdequivFinite}, $P$ is
finite. If $g \in {\cal RM}^{\sf rec}$ is such that ${\sf domC}(g)$ is
finite, then $g \in {\cal RM}^{\sf P}$.
By Lemma \ref{regLR}, if $g \in {\cal RM}^{\sf P}$ is injective and
${\sf domC}(g)$ is a finite maximal prefix code, then $g \equiv_{\cal L}$
${\sf id}_{{\sf domC}(g) \, A^*}$; hence  $[[g]]_{\sf bd}$ $\equiv_{\cal L}$
$[[{\sf id}]]_{\sf bd}$ (in ${\cal M}_{\sf bd}^{\sf P}$).
For any $f \in {\cal RM}^{\sf rec}$ with finite ${\sf imC}(f)$, we want to
show that $f \equiv_{\cal R} g$ in ${\cal RM}^{\sf rec}$ for some
$g$ of the above type. Then we will have
$\, [[f]]_{\sf bd} \equiv_{\cal R}$
$[[g]]_{\sf bd}$ $\equiv_{\cal L} [[{\sf id}]]_{\sf bd} \, $ in
${\cal M}_{\sf bd}^{\sf P}$; so every $[[f]]_{\sf bd} \in \Delta_1$ will
be in the $\cal D$-class of $[[{\sf id}]]_{\sf bd}$.

Let ${\sf imC}(f) = \{y_i : i = 1, \ldots, N\}$, where
$N = |{\sf imC}(f)| \, $ (finite).
Let $X = \{x_i : i = 1, \ldots, N\} \subseteq {\sf domC}(f)$ be such that
$x_i \in f^{-1}(y_i)$; i.e., $X$ is a choice set for the restriction of
$f^{-1}$ to ${\sf imC}(f)$. Since $X \subseteq {\sf domC}(f)$, $X$ is a
finite prefix code.
Let $P \subset A^*$ be any finite maximal prefix code of size $N$, and let
$(p_1, \ldots, p_N)$ be any total ordering of $P$.
We define an injective right-ideal morphism $\alpha \in {\cal RM}^{\sf P}$
($\subseteq {\cal RM}^{\sf rec}$)
by $\, {\sf domC}(\alpha) = P$, \ ${\sf imC}(\alpha) = X$, and
$\alpha(p_i w) = x_i w$ \ (for all $i = 1, \ldots, N$, and $w \in A^*$).

Let $g = f \circ \alpha$. So, ${\sf domC}(g) = P \,$ (which is a finite
maximal prefix code), and ${\sf imC}(g) = {\sf imC}(f)$.
And $g$, restricted to ${\sf domC}(g)$, is a bijection from $P$ to
${\sf imC}(f)$, so $g$ is injective on ${\sf Dom}(g) = P A^*$.
(In fact, $[g]_{\sf bd} \in M_{2,1}$, the Thompson-Higman monoid defined in
\cite{JCBmonThH}).)
Let $\beta = g^{-1} \circ f$; this is well defined since $g$ is injective.
Note that $g \circ g^{-1}$ is the restriction of the identity map to
$\, {\sf imC}(f) \, A^* = {\sf Im}(f) = {\sf Im}(g)$; thus,
$g \circ \beta = g \circ g^{-1} \circ f = {\sf id}_{{\sf Im}(f)} \circ f = f$.
Thus, $g = f \circ \alpha$ and $g \circ \beta = f$, so $g \equiv_{\cal R} f$
in ${\cal RM}^{\sf rec}$.
 
\smallskip

\noindent 
(2) Let us prove that all elements of $\Delta_2$ are ${\cal D}$-related.

\smallskip

By Prop.\ \ref{RMRecPolyRMPpolybd}, ${\cal M}_{\sf bd}^{\sf P}$
($ = {\cal RM}^{\sf P}\!/\!\equiv_{\sf bd}$) is isomorphic to
${\cal RM}^{\sf rec}\!/\!\equiv_{\sf bd}$.

\smallskip

\noindent {\sf Claim 1:} \ Let $P, Q \subset A^*$ be prefix codes that
are infinite and decidable. Then there exists a bijection
$\alpha \in {\cal RM}^{\sf rec}$ from $P$ onto $Q$.

\smallskip

\noindent {\sf Proof:} For any infinite set $S \subseteq A^*$, the rank
function of $S$ is a bijection from $S$ onto ${\mathbb N}$, defined for
$x \in S$ by \ ${\sf rank}_S(x) \ = \ |\{w \in S: w <_{\ell \ell} x\}|$.
To make ${\sf rank}_S(.)$ a function between words, we represent a natural
integer $n \in {\mathbb N}$ by $0^n 1$, so ${\sf Im}({\sf rank}_S) = 0^*1$.

If $S$ is decidable, ${\sf rank}_S$ is partial recursive with decidable
domain $S$. And ${\sf rank}_S$ has a computable input-output balance when
$S$ is infinite and decidable. Then
$\, \alpha = {\sf rank}_Q^{-1} \circ {\sf rank}_P \, $ is a bijection from
$P$ onto $Q$ with the claimed properties. Finally, $\alpha$ can be extended
to a bijective right-ideal morphism from $P A^*$ onto $Q A^*$;
thus, $\alpha \in {\cal RM}^{\sf rec}$.
This proves Claim 1.

\smallskip

As a consequence of the proof of Claim 1,
${\sf id}_{P A^*} = \alpha \circ {\sf id}_{Q A^*}$ and
${\sf id}_{Q A^*} = \alpha^{-1} \circ {\sf id}_{P A^*}$. Hence for all
infinite decidable prefix codes $P$ and $Q$ we have:
$\, {\sf id}_{P A^*} \equiv_{\cal L} {\sf id}_{Q A^*} \,$ in
${\cal RM}^{\sf rec}$.

\smallskip

\noindent {\sf Claim 2:} \ For any $f \in {\cal RM}^{\sf rec}$,
 \ $f \, \equiv_{\cal R} \, {\sf id}_{{\sf Im}(f)}$.

\smallskip

\noindent {\sf Proof:} For $f \in {\cal RM}^{\sf rec}$, ${\sf Im}(f)$ is
a decidable set, because of the computable I/O-balance. Every
$f \in {\cal RM}^{\sf rec}$ has an inverse in ${\cal RM}^{\sf rec}$, and
since ${\sf Im}(f)$ is decidable, such an inverse can be restricted to
${\sf Im}(f)$. If $f'$ is such an inverse with domain ${\sf Im}(f)$, we
have: $f \circ f' = {\sf id}_{{\sf Im}(f)}$.
Moreover, ${\sf id}_{{\sf Im}(f)} \circ f = f$. This proves Claim 2.

\smallskip

By Claim 2 and the consequence of Claim 1 we now have: If
$f, g \in {\cal RM}^{\sf rec}$ and if ${\sf imC}(f)$, ${\sf imC}(g)$ are
infinite, then $f \, \equiv_{\cal R} \, {\sf id}_{{\sf Im}(f)} \, $
$\equiv_{\cal L} \, {\sf id}_{{\sf Im}(g)} \, \equiv_{\cal R} \, g$. 
 \ \ \ $\Box$


\section{Polynomial and exponential end-equivalences }

The relations $\equiv_{\cal T}$, in particular $\equiv_{\sf poly}$ and 
$\equiv_{\sf E3}$, were defined in Def.\ \ref{equiv_boundDEF} for prefix
codes, and in Def.\ \ref{boundequiv_RM} for right-ideal morphisms. 
We call $\equiv_{\sf poly}$ the {\em polynomial end-equivalence} relation, 
and $\equiv_{\sf E3}$ the {\em exponential end-equivalence} 
(or elementary recursive end-equivalence) relation.
Throughout this section, $\cal T$ denotes a family of functions as in Def.\
\ref{equiv_boundDEF}, possibly with additional properties.
When $\equiv_{\cal T}$ is applied between prefix codes, it is an
equivalence relation. Transitivity follows from the fact that $\cal T$ is
closed under composition.
For two prefix codes $P_1, P_2 \subset A^*$ and $\tau \in {\cal T}$, we say
that lengths in $P_1$ and $P_2$ are $\tau$-related iff
$|x_1| \le \tau(|x_2|)$ and $|x_2| \le \tau(|x_1|)$ for all $x_1 \in P_1$,
$x_2 \in P_2$ with $x_1 \|_{\sf pref} x_2$.
In particular, when $\cal T$ is the set of polynomials we say 
{\em ``polynomially related''}.
The latter is the most interesting, due to its connections with {\sf NP}. 

If $L_1, L_2, M_1, M_2$ are prefix codes such that 
$\, L_1 \equiv_{\cal T} L_2$, $\, M_1 \subseteq L_1$, 
$\, M_2 \subseteq L_2$, and $\, M_1 \equiv_{\sf end} M_2$, 
then $\, M_1 \equiv_{\cal T} M_2$.
Indeed, the bounding function $\tau \in {\cal T}$ that appear in
the definition of $L_1 \equiv_{\cal T} L_2$, also works for
$M_1 \equiv_{\cal T} M_2$.

There exists $f \in {\cal RM}^{\sf P}$ such that
$f({\sf domC}(f))$ $\not\equiv_{\cal T}$ ${\sf imC}(f)$, and 
$\equiv_{\cal T}$ is not reflexive on $f({\sf domC}(f))$. Moreover,
$f$ can be chosen so that there exist prefix codes 
$P_1, P_2 \subset {\sf Dom}(f)$ with 
$P_1 \equiv_{\cal T} P_2$, such that $f(P_1) \not\equiv_{\cal T} f(P_2)$.
As an example, let $f \in {\cal RM}^{\sf P}$ be defined by
$f(0^n1w) = 0^n w$ for all $n \ge 0$ and $w \in \{0,1\}^*$; so
${\sf domC}(f) = 0^*1$, and ${\sf imC}(f) = \{ \varepsilon\}$. 
Then, $f({\sf domC}(f)) = 0^* \not\equiv_{\sf bd} \{ \varepsilon\}$ $=$
${\sf imC}(f)$, and $f({\sf domC}(f)) = 0^* \not\equiv_{\sf bd} 0^*$ $=$
$f({\sf domC}(f))$ (non-reflexive). Note that $\not\equiv_{\sf bd}$ 
implies $\not\equiv_{\cal T}$. 
To show the possibility of $f(P_1) \not\equiv_{\cal T} f(P_2)$ when 
$P_1 \equiv_{\cal T} P_2$, let $f$ be as in the example above, and
let $P_1 = P_2 = 0^*1$. Then $f(P_1) = f(P_2) = 0^*$; but 
$0^* \not\equiv_{\sf bd} 0^*$ (non-reflexivity in this case).
 
There exist $f \in {\cal RM}^{\sf P}$ and prefix codes
$P_1, P_2 \subset {\sf Im}(f)$ such that $P_1 \equiv_{\sf poly} P_2$ but
$f^{-1}(P_1) \not\equiv_{\sf end} f^{-1}(P_2)$.
For example, let 
${\sf domC}(f) = \{00, 01, 1\}$, and $f(00) = 00$, $f(01) = 0$, 
$f(1) = \varepsilon$; so $f \in {\cal RM}^{\sf P}$.
Then $f^{-1}(\{ \varepsilon\}) = \{ 1\}$, $f^{-1}(\{0\}) = \{01, 10\}$, 
$f^{-1}(\{1\}) = \{11\}$, $f^{-1}(\{01\}) = \{101, 011\}$,
$f^{-1}(\{00\}) = \{100, 010, 00\}$. 
Let $P_1 = \{ \varepsilon\}$, $P_2 = \{0,1\}$, and $P_3 = \{00, 01, 1\}$.
Then $P_1 \equiv_{\sf poly} P_2 \equiv_{\sf poly} P_3$, but $f^{-1}(P_1)$,
$f^{-1}(P_2)$, and $f^{-1}(P_3)$ are all $\not\equiv_{\sf end}$,
since $f^{-1}(P_1) = \{11\}$, $f^{-1}(P_2) = \{01,10,11\}$, and
$f^{-1}(P_3) = \{11, 100, 010, 00, 101, 011\}$.

\bigskip

\noindent The following is the $\equiv_{\cal T}$ version of Propositions
\ref{endequivUnionIntersLatt} and \ref{bdequivUnionIntersLatt}.

\begin{pro} \label{polyequivIntersection}
For any family of functions $\cal T$, as in Def.\ \ref{equiv_boundDEF} we
have the following.

\smallskip

\noindent (1) \ Let $P_1, P_2 \subset A^*$ be prefix codes such that
$P_1 \equiv_{\cal T} P_2$, let $P_{\cap}$ be the prefix code that 
generates the right ideal $P_1 A^* \cap P_2 A^*$, and let $P_{\cup}$ be the
prefix code that generates the right ideal $P_1 A^* \cup P_2 A^*$.  Then 
$P_1 \equiv_{\cal T} P_2$
$\equiv_{\cal T} P_{\cap} \ \equiv_{\cal T} P_{\cup}$.

\smallskip

\noindent (2) \ Let $f_1, f_2$ be right-ideal morphisms such that
$f_1 \equiv_{\cal T} f_2$. Then $f_1 \cap f_2$ and $f_1 \cup f_2$ are 
right-ideal morphisms, and $f_1 \equiv_{\cal T} f_2$
$\equiv_{\cal T} f_1 \cap f_2 \equiv_{\cal T} f_1 \cup f_2$.
\end{pro}
{\bf Proof.} 
(1) By Prop.\ \ref{endequivUnionIntersLatt}, 
$P_{\cap} \equiv_{\sf end} P_{\cup} \equiv_{\sf end} P_1$ 
$\equiv_{\sf end}$ $P_2$.

Length bounds: It is known (and easily proved) that $P_{\cup}$ 
$\subseteq P_1 \cup P_2$.  Hence, if $x \in P_{\cup}$, $p_1 \in P_1$, and 
$x \parallel_{\sf pref} p_1$, then either $x \in P_1$ (and then $x = p_1$), 
or $x \in P_2$ (and then $|x|$, $|p_1|$ are $\tau$-related since 
$P_1 \equiv_{\cal T} P_2$). 
Similarly, if $x \in P_{\cup}$, $p_2 \in P_2$, and $x \parallel_{\sf pref}$
$p_2$, then $|x|$, $|p_2|$ are $\tau$-related. 

It is also the case that $P_{\cap} \subseteq P_1 \cup P_2$, and a similar 
reasoning applies here.

\noindent (2) We conclude from (1) that ${\sf Dom}(f_i)$ $\equiv_{\cal T}$
${\sf Dom}(f_1) \cap {\sf Dom}(f_2)$ $\equiv_{\cal T}$ 
${\sf Dom}(f_1) \cup {\sf Dom}(f_2)$, for $i = 1, 2$. 

Also, ${\sf Dom}(f_1 \cap f_2)$  $=$ ${\sf Dom}(f_1) \cap {\sf Dom}(f_2)$.
And since $f_1 = f_2$ on ${\sf Dom}(f_1 \cap f_2)$ we have
$f_1 = f_1 \cap f_2$ on ${\sf Dom}(f_1 \cap f_2)$. Hence 
$f_1 \equiv_{\cal T} f_1 \cap f_2$, and similarly for $f_2$. 

Also, ${\sf Dom}(f_1 \cup f_2)$  $=$ ${\sf Dom}(f_1) \cup {\sf Dom}(f_2)$.
And $f_1 = f_1 \cup f_2$ on ${\sf Dom}(f_1)$, and 
$f_2 = f_1 \cup f_2$ on ${\sf Dom}(f_2)$, hence 
$f_1 \equiv_{\cal T} f_1 \cup f_2 \equiv_{\cal T} f_2$.
 \ \ \ $\Box$

\begin{cor} \label{equiviPolyclassLattice} 
 \ Every $\, \equiv_{\cal T}$-class in ${\cal RM}^{\sf P}$ is a lattice 
under $\subseteq$, $\cup$ and $\cap$.   \ \ \ $\Box$
\end{cor}

\begin{thm} \label{polyequivCongr} 
 \ Let $\cal T$ be any family of functions as in Def.\ \ref{equiv_boundDEF}.
and let $\cal M$ be any monoid of right-ideal morphisms with I/O-balance 
function in $\cal T$. Then the relation $\, \equiv_{\cal T}$ is a congruence 
on $\cal M$.
\end{thm}
{\bf Proof.} Clearly, $\equiv_{\cal T}$ is an equivalence relation (for
transitivity we use the fact that $\cal T$ has upper bounds for composition,
see Def.\ \ref{equiv_boundDEF}(2)).
For the multiplicative property, let $f_1, f_2, g \in {\cal M}$, and
suppose $f_1 \equiv_{\cal T} f_2$;  we want to prove that
$f_1 \, g \equiv_{\cal T} f_2 \, g$  and
$g \, f_1 \equiv_{\cal T} g \, f_2$.
Since $\equiv_{\cal T}$ implies $\equiv_{\sf bd}$, the actions of $f_1$ and
$f_2$ on $A^{\omega}$ are the same, hence
$f_1 g \equiv_{\sf bd} f_2 g$, and $g f_1 \equiv_{\sf bd} g f_2$ (using
Cor.\ \ref{end_boundActCOR}). It now suffices to check the
$\cal T$-relation for lengths in the domain codes. 

\smallskip

\noindent $\bullet$  
 \ Proof that \, ${\sf domC}(f_1 g) \equiv_{\cal T} {\sf domC}(f_2 g)$: 

e want to show that lengths in ${\sf domC}(f_1 g)$ and ${\sf domC}(f_2 g)$
are $\cal T$-related.
Let $x_1 \in {\sf domC}(f_1g) \,$ and $ \, x_2 \in {\sf domC}(f_2g) \,$ be
prefix-comparable. By Prop.\ \ref{polyequivIntersection}(2) we can assume
that $f_2 \subseteq f_1$, hence ${\sf Dom}(f_2) \subseteq {\sf Dom}(f_1)$;
hence, $x_2 \ge_{\sf pref} x_1$ (i.e., $x_1$ is a prefix of $x_2$).

Since $x_2 \ge_{\sf pref} x_1$, we have $g(x_2) \ge_{\sf pref} g(x_1)$.
Since $x_i \in {\sf domC}(f_i g)$, we have $g(x_i) \in {\sf Dom}(f_i)$;
hence there exists $z_i \in {\sf domC}(f_i)$ with
$g(x_i) \ge_{\sf pref} z_i$ (for $i = 1, 2$).  Since
$z_2 \le_{\sf pref} g(x_2) \ge_{\sf pref} g(x_1) \ge_{\sf pref} z_1$, we
have $z_2 \parallel_{\sf pref} g(x_1)$ and $z_2 \parallel_{\sf pref} z_1$.
Since $f_1 \equiv_{\cal T} f_2$, $|z_1|$ and $|z_2|$ are
$\tau _{12}$-related for some $\tau_{12}  \in {\cal T}$ (depending only on
$f_1, f_2$). We can assume $\tau_{12}(n) \ge n$ for all
$n \in {\mathbb N}$ and that $\tau_{12}$ is increasing, since ${\cal T}$
contains such a function and ${\cal T}$ has upper bounds for sum (see
Def.\ \ref{equiv_boundDEF}(2)).

Since $f_2 \subseteq f_1$, ${\sf Dom}(f_2) \subseteq {\sf Dom}(f_1)$, hence
for all $d_2 \in {\sf domC}(f_2)$ and all $d_1 \in {\sf domC}(f_1)$: if
$d_2 \parallel_{\sf pref} d_1$ then $d_2 \ge_{\sf pref} d_1$; therefore,
$z_2 \ge_{\sf pref} z_1$.
We now have two cases:

 \ \ \ \ \ $g(x_2) \ \ge_{\sf pref} \ g(x_1) \ \ge_{\sf pref} \ z_2$
    $ \ge_{\sf pref} \ z_1$, \ \ or

 \ \ \ \ \ $g(x_2) \ \ge_{\sf pref} \ z_2 \ \ge_{\sf pref} \ g(x_1)$
    $ \ge_{\sf pref} \ z_1$.

\smallskip

\noindent Case 1: \ $g(x_1) \ge_{\sf pref} z_2$.

Now $g(x_1) \in {\sf Dom}(f_2)$ (since $g(x_1) \ge_{\sf pref} z_2$), hence
$x_1 \in g^{-1}({\sf Dom}(f_2))$, so $x_1 \in {\sf Dom}(f_2 g)$. Since
$x_2 \in {\sf domC}(f_2g)$ and $x_2 \ge_{\sf pref} x_1$ it follows that 
$x_2 = x_1$. 

This implies obviously that $|x_2| = |x_1|$, hence $|x_2|$ and $|x_1|$ are
$\tau_{12}$-related (since $\tau_{12}(n) \ge n$).

\smallskip

\noindent Case 2: \ $z_2 \ge_{\sf pref} g(x_1)$.

Being a right-ideal morphism of $A^*$, $g$ is a one-to-one correspondence 
between the $\ge_{\sf pref}$-chains 

\smallskip

 \ \ \ \ \ $x_2 \ge_{\sf pref} \ \ldots \ \ge_{\sf pref} x_1$ \ \ and  

 \smallskip

 \ \ \ \ \ $g(x_2) \ge_{\sf pref} \ \ldots \ \ge_{\sf pref} g(x_1)$ 

\smallskip

\noindent in $A^*$. (It is easy to see that any right-ideal morphism $f$ 
of $A^*$ is injective on any chain 
$\, x >_{\sf pref} xm_1 >_{\sf pref} xm_1 m_2 >_{\sf pref} \  \ldots \ $
$>_{\sf pref} xm_1 m_2 \ldots m_k$, if $x \in {\sf Dom}(f)$; this holds 
even if $f$ is not injective on all of ${\sf Dom}(f)$.)
Since $g(x_2) \ge_{\sf pref} z_2 \ge_{\sf pref} g(x_1)$, let $t_2$ be the
(unique) inverse image of $z_2$ in the upper chain; so 
$x_2 \ge_{\sf pref} t_2 \ge_{\sf pref} x_1$, and $g(t_2) = z_2$ 
($\in {\sf domC}(f_2)$). 
Then $t_2 \in g^{-1}({\sf domC}(f_2))$, hence $t_2 \in {\sf Dom}(f_2 g)$.
Moreover, $t_2 \le_{\sf pref} x_2 \in {\sf domC}(f_2 g)$ implies that
$t_2 = x_2$ (since ${\sf domC}(f_2 g)$ is a prefix code, and in a prefix 
code, prefix comparable elements are equal).
Therefore, $g(t_2) = g(x_2)$, hence (since $g(t_2) = z_2$), $z_2 = g(x_2)$.
Thus, $g(x_2) \in {\sf domC}(f_2)$.
Since we also have $z_1 \in {\sf domC}(f_1)$ we conclude that
$\, |g(x_2)| \le \tau_{12}(|z_1|) \, $ (since lengths in
${\sf domC}(f_2)$ and ${\sf domC}(f_1)$ are $\tau_{12}$-related).
And $\, |z_1| \le |g(x_1)| \,$ (since $g(x_1) \ge_{\sf pref} z_1$), hence
$\, |g(x_2)| \le \tau_{12}(|g(x_1)|) \,$ (since $\tau_{12}$ is increasing).
Letting $\tau_g$ denote the balance polynomial of $g$, we obtain:
$ \ |x_2| \, \le \, \tau_g(|g(x_2)|) \, $
$\le \, \tau_g(\tau_{12}(|g(x_1)|)) \, $
$\le \, \tau_g \circ \tau_{12} \circ \tau_g(|x_1|)$.

We also have $\, |x_2| \ge |x_1| \, $ (since $x_2 \ge_{\sf pref} x_1$).

In summary,  ${\sf domC}(f_1 g) \equiv_{\cal T} {\sf domC}(f_2 g)$ for any
function in $\cal T$ that bounds $\tau_g \circ \tau_{12} \circ \tau_g$
from above.

\medskip 

\noindent $\bullet$  
 \ Proof that \, ${\sf domC}(g f_1) \equiv_{\cal T} {\sf domC}(g f_2)$: 

Let $x_1 \in {\sf domC}(g f_1), \ x_2 \in {\sf domC}(g f_2)$ be
prefix-comparable. We want to show that $|x_1|, \, |x_2|$ are 
$\cal T$-related. As before, we can assume that $f_2 \subseteq f_1$, hence 
$x_2 \ge_{\sf pref} x_1$ (i.e., $x_1$ is a prefix of $x_2$).

Since  $g f_1(x_1)$ is defined, $f_1(x)$ is defined for all 
$x \ge_{\sf pref} x_1$. 

Since $f_2 \subseteq f_1$ and since $x_2 \in {\sf Dom}(g f_2)$
$\subseteq {\sf Dom}(f_2) \subseteq {\sf Dom}(f_1)$, we have 
$f_2(x_2) = f_1(x_2)$.  
Let $z_2 \in {\sf domC}(f_2)$ be such that $x_2 \ge_{\sf pref} z_2$. Then, 
$z_2 \le_{\sf pref} x_2 \ge_{\sf pref} x_1$,
hence $z_2 \parallel_{\sf pref} x_1$; so we have two cases: 
 \ $x_1 \ge_{\sf pref} z_2$, or 
 \ $x_2 \ge_{\sf pref} z_2 \ge_{\sf pref} x_1$. 

\smallskip

\noindent Case 1: \ $x_1 \ge_{\sf pref} z_2$ ($\in {\sf domC}(f_2)$):

Then $x_1 \in {\sf Dom}(f_2)$, thus $f_2(x_1) = f_1(x_1)$ (since 
$f_2 \subseteq f_1$, and $f_2(x_1)$ is defined). So $gf_2(x_1) = gf_1(x_1)$,
and $gf_2(x_1)$ is defined, i.e., $x_1 \in {\sf Dom}(g f_2)$. 
Since $x_2 \in {\sf domC}(f_2)$ and  $x_2 \parallel_{\sf pref} x_1$, it
follows that $x_1  \ge_{\sf pref} x_2$.
So $x_1 = x_2$ (since we also have $x_2 \ge_{\sf pref} x_1$). 
So $|x_1|, \, |x_2|$ are $\cal T$-related.
 
\smallskip

\noindent Case 2: \ $x_2 \ge_{\sf pref} z_2 \ge_{\sf pref} x_1$.

Since $z_2 \in {\sf domC}(f_2)$, $f_2(z_2)$ is defined, hence 
$f_2(z_2) = f_1(z_2)$ (since $f_2 \subseteq f_1$). Hence $g f_2(z_2)$ is
defined (since $g f_1(x)$ is defined for all $x \ge_{\sf pref} x_1$, and
$f_2(z_2) = f_1(z_2)$), i.e., $z_2 \in {\sf Dom}(g f_2)$. Therefore, 
$z_2 \ge_{\sf pref} x_2$, since $z_2 \parallel_{\sf pref} x_2$ and 
$x_2 \in {\sf domC}(g f_2)$.
But since we also have $x_2 \ge_{\sf pref} z_2$, it follows that 
$z_2 = x_2$.

Let $z_1 \in {\sf domC}(f_1)$ be such that $x_1 \ge_{\sf pref} z_1$.
Then we have the $\ge_{\sf pref}$-chain 
 \ $x_2 = z_2 \ge_{\sf pref} x_1 \ge_{\sf pref} z_1$.
Since $|z_2| \le \tau_{12}(|z_1|)$ it follows (since $x_2 = z_2$ and
$|z_1| \le |x_1|$) that $|x_2| \le \tau_{12}(|x_1|)$. Also,
$|x_1| \le |x_2|$.
So, $|x_1|$ and $|x_2|$ are $\tau_{12}$-related.

In summary, ${\sf domC}(g f_1) \equiv_{\cal T} {\sf domC}(g f_2)$ for the
function $\tau_{12}$. 
 \ \ \ $\Box$

\bigskip

\noindent {\bf Notation:} Let ${\cal RM}^{\sf P}\!/\!\! \equiv_{\cal T}$
denote the set of $\, \equiv_{\cal T}$-congruence classes in 
${\cal RM}^{\sf P}$; we will abbreviate 
${\cal RM}^{\sf P}\!/\! \equiv_{\cal T}$ by ${\cal M}^{\sf P}_{\cal T}$.
In particular, we will consider ${\cal M}^{\sf P}_{\sf poly}$, 
${\cal M}^{\sf P}_{\sf E3}$, and ${\cal M}^{\sf P}_{\sf lin}$.
 
Similarly, ${\cal RM}^{\bf NP}\!/\! \equiv_{\cal T}$ denotes the set of
$\, \equiv_{\cal T}$-congruence classes in ${\cal RM}^{\bf NP}$.  Here, 
${\cal RM}^{\bf NP}$ (also called ${\cal RM}^{\Sigma_1}$) is the monoid 

\medskip

${\cal RM}^{\bf NP} \ = \ $
$\{f : f$ is a polynomially balanced right-ideal morphism $A^* \to A^*$
that is computable 

\hspace{1.1in}  by a polynomial-time deterministic Turing machine with an
oracle in {\sf NP}\}.

\medskip

\noindent See \cite{s1f}, section 6, for more details on the similarly 
defined ${\sf fP}^{\bf NP}$.

\medskip

By Theorem \ref{polyequivCongr} we have (if
${\sf poly} \subseteq {\cal T}$):

\begin{cor} \label{monoidEmbedding}
 \ Let ${\cal T}$ be as in Def.\ \ref{equiv_boundDEF} with, in addition,
${\sf poly} \subseteq {\cal T}$. Then for every monoid $\cal M$ of 
right-ideal morphisms with polynomial I/O-balance, 
${\cal M}/\! \equiv_{\cal T}$ is a monoid. 
In particular, ${\cal RM}^{\bf P}\!/\! \equiv_{\sf poly}$ 
($ = {\cal M}^{\sf P}_{\sf poly}$) and 
$\, {\cal RM}^{\bf NP}\!/\! \equiv_{\sf poly}$ are monoids, 
and there is a homomorphic embedding

\smallskip 

\hspace{.6in}
${\cal M}^{\sf P}_{\sf poly}$ \ \ $\hookrightarrow$
 \ \ ${\cal RM}^{\bf NP}\!/\! \equiv_{\sf poly}$ .
\end{cor}
{\bf Proof.} The first statements follow from the fact that 
$\equiv_{\cal T}$ is a congruence. 

Every $\, \equiv_{\sf poly}$-class of ${\cal RM}^{\bf NP}$ contains at most 
one $\, \equiv_{\sf poly}$-class of ${\cal RM}^{\sf P}$, since 
$\, \equiv_{\sf poly}$ is transitive; hence we have the embedding.
  \ \ \ $\Box$

\bigskip

\noindent But if ${\sf P} \neq {\sf NP}$ then a $\, \equiv_{\sf poly}$-class
of ${\cal RM}^{\bf NP}$ that contains elements of ${\cal RM}^{\sf P}$ could
also contain functions that are not in ${\cal RM}^{\sf P}$; i.e., a
$\, \equiv_{\sf poly}$-class of ${\cal RM}^{\sf P}$ could be a strict subset
of the corresponding $\, \equiv_{\sf poly}$-class if
${\sf P} \neq {\sf NP}$.
So if ${\sf P} \neq {\sf NP}$, the embedding above is not an inclusion.

\bigskip

Let ${\cal T}_1, {\cal T}_2$ be families of functions as in Def.\
\ref{equiv_boundDEF}. If ${\cal T}_1 \subseteq {\cal T}_2$ then
$\, \equiv_{{\cal T}_1} \ \subseteq \ \equiv_{{\cal T}_2}$; hence there 
exists a surjective monoid morphism 
$\, {\cal RM}^{\sf P}\!/\! \equiv_{{\cal T}_1} \ $
$\twoheadrightarrow$ $ \ {\cal RM}^{\sf P}\!/\! \equiv_{{\cal T}_2}$.
In particular we have surjective monoid morphisms

\medskip

\hspace{.8in}  ${\cal M}^{\sf P}_{\sf poly}$ $ \ \twoheadrightarrow \ $
${\cal M}^{\sf P}_{\sf E3}$ $ \ \twoheadrightarrow \ $
${\cal M}^{\sf P}_{\sf bd}$ $ \ \twoheadrightarrow \ $
${\cal M}^{\sf P}_{\sf end}$ $ \ \twoheadrightarrow \ $ \{{\bf 1}\}.

\medskip

\noindent Since ${\cal M}^{\sf P}_{\sf end}$ is congruence-simple, the 
right-most arrow (onto the one-element monoid) cannot be factored (except 
by using automorphisms as factors).

%

\begin{pro} \label{equitoLin} 
 \ Let $\cal T$ be as in Def.\ \ref{equiv_boundDEF}, with the additional
condition that ${\sf poly} \subseteq {\cal T}$. 
Then every $\, \equiv_{\cal T}$-class of ${\cal RM}^{\sf P}$ contains 
functions whose I/O-balance and time-complexity are linear (bounded from
above by the function $n \mapsto 3n$). 
\end{pro}
{\bf Proof.} We proceed as in Lemma \ref{recTimeConstr}(2).
For any $f \in {\cal RM}^{\sf P}$ and any function $\tau \in {\cal T}$, we
define the right-ideal morphism $F_{f,\tau}$ by

\medskip

 \ \ \ \ \ ${\sf domC}(F_{f,\tau}) \ = \ $
$\bigcup_{x \in {\sf domC}(f)} \ x \ \{0,1\}^{|x| \cdot \tau(|x|)}$, 
 \ \ and

\medskip

 \ \ \ \ \ $F_{f,\tau}(x z w) \ = \ f(x) \, z  w$, 

\medskip

\noindent for all $x \in {\sf domC}(f)$, $z \in$ 
$\{0,1\}^{|x| \cdot \tau(|x|)}$, and $w \in \{0,1\}^*$.  Then 
$f \equiv_{\cal T} F_{f,\tau} \, $ (since $\{0,1\}^{|x| \cdot \tau(|x|)}$ 
is a maximal prefix code, see Lemma \ref{recTimeConstr}(2)). 

Let $\tau$ be a polynomial upper bound on the I/O-balance and the
time-complexity of $f$. We can choose $\tau$ to be a polynomial
of the form $n \mapsto a \cdot (n^d + 1)$; then $\tau \in {\cal T}$ and 
$\tau$ is fully time-constructible.
Then $F_{f,\tau}$ has linear I/O-balance and time-complexity with 
coefficient $le 3$, by Lemma \ref{recTimeConstr}(2).  
 \ \ \ $\Box$

\bigskip

We saw in \cite{s1f} that {\sf fP} and ${\cal RM}^{\sf P}$ do not contain
evaluation maps (contrary to the monoid {\sf fR} of all partial recursive 
functions with partial recursive balance function).
By definition, a (coded) evaluation map for {\sf fR} is a partial function 
${\sf eval} \in {\sf fR}$ such that for every $f \in {\sf fR}$ there exists 
$w \in \{0,1\}^*$ (called a program for $f$) such that for all 
$x \in {\sf Dom}(f)$: ${\sf eval}({\sf code}(w) \ 11 \, x) = f(x)$. 
We saw that {\sf fP} and ${\cal RM}^{\sf P}$ contain partial functions 
that play the role of evaluation maps in a limited way: For every polynomial
$q$ of degree $> 1$ there exists ${\sf eval}_q$ that works as an evaluation
map for functions whose time-complexity and I/O-balance are less than $q$;
for {\sf fP}, see section 4 of \cite{s1f}, for ${\cal RM}^{\sf P}$, see
section 2 of \cite{BiInfGen}.

An interesting property of ${\cal M}^{\sf P}_{\sf poly}$ is that it has 
``evaluation elements'' that play the same role as evaluation maps; of 
course, elements of ${\cal M}^{\sf P}_{\sf poly}$ are not maps but 
equivalence classes of maps. We will also see that 
${\cal M}^{\sf P}_{\sf poly}$, contrary to ${\cal RM}^{\sf P}$, is finitely 
generated.

\begin{defn} \label{evalElement}
 \ A class $[e_0] \in {\cal M}^{\sf P}_{\sf poly}$ is called an 
{\em evaluation element} iff 
there exists $e \in [e_0]$ such that for every 
$[f_0] \in {\cal M}^{\sf P}_{\sf poly}$, there exists $u \in \{0,1\}^*$ 
such that $[e({\sf code}(u) \ 11 \, (\cdot))] = [f_0(\cdot)]$.

Here ${\sf code}(u) \ 11 \, (\cdot)$ denotes the function
$\, x \in \{0,1\}^* \longmapsto {\sf code}(u) \ 11 \, x$.
\end{defn}
Equivalently, 

$(\exists e \in [e_0])$
$(\forall [f_0] \in {\cal M}^{\sf P}_{\sf poly})$
$(\exists f \in [f_0])$ $(\exists u \in \{0,1\}^*)$
$(\forall x \in {\sf Dom}(f)) : \ \  e({\sf code}(u) \ 11 \, x) = f(x)$.
 
The function ${\sf code}(.)$ was defined in the Introduction (just before
the definition of ${\cal RM}^{\sf P}$).

\begin{thm} \label{evalfuncRMpolyequ} 
 \ The monoid ${\cal M}^{\sf P}_{\sf poly}$ has evaluation elements and is 
finitely generated.
\end{thm}
{\bf Proof.} For any polynomial $q$ of the form $q(n) = a \cdot (n^d +1)$
with $d >1$ and $a \ge 3$, we consider the evaluation function
${\sf evalR}_q^C$ defined by
$\, {\sf evalR}_q^C({\sf code}(u) \ 11 \, x z)$ $=$ $\phi_u(x) \ z$;
here $u$ is any program with linear time-complexity and I/O-balance (with
coefficient $\le 3$), $x \in {\sf domC}(\phi_u)$, and $z \in A^*$.
By Prop.\ \ref{equitoLin}, every $\phi_v \in {\cal RM}^{\sf P}$ is
$\equiv_{\sf poly}$-equivalent to some $\phi_u \in {\cal RM}^{\sf P}$ such
that $\phi_u$ has time-complexity and I/O-balance less than the function
$n \mapsto 3n$; thus, $[{\sf evalR}_q^C]$ is an evaluation element.

Defining ${\sf evR}_q^C$ by
$\, {\sf evR}_q^C({\sf code}(u) \ 11 \, x z)$ $=$
${\sf code}(u) \ 11 \ \phi_u(x) \ z$,  we also have

\smallskip

\hspace{.5in}
    $\phi_u \ = \ \pi'_{|{\sf code}(u) \, 11|} \circ {\sf evR}_q^C$
    $\circ$  $\pi_{{\sf code}(u) \, 11}$.

\smallskip

\noindent So, ${\cal M}^{\sf P}_{\sf poly}$ is generated by
$\{[\pi_0], \, [\pi_1], \, [\pi_1'], \, [{\sf evR}_q^C] \}$.
 \ \ \ $\Box$

\bigskip

\noindent We saw that when ${\sf poly} \subseteq {\cal T}$ then 
${\cal M}^{\sf P}_{\cal T}$ is a homomorphic image of 
${\cal M}^{\sf P}_{\sf poly}$. Hence we have:

\begin{cor}
 \ If ${\sf poly} \subseteq {\cal T}$ then ${\cal M}^{\sf P}_{\cal T}$ is
finitely generated.   \ \ \ \ \ $\Box$
\end{cor}

We do not know whether ${\cal M}^{\sf P}_{\sf poly}$ is regular (and this
is equivalent to ${\sf P} = {\sf NP}$ by Theorem \ref{regofRMequ}), but for
${\cal M}^{\sf P}_{\sf E3}$ we can prove:

\begin{pro} \label{M_Ereg} 
 \ The monoid ${\cal M}^{\sf P}_{\sf E3}$ is regular.
\end{pro}
{\bf Proof.} Consider $[f] \in {\cal M}^{\sf P}_{\sf E3}$, i.e., an 
$\equiv_{E3}$-class in ${\cal RM}^{\sf P}$ for some 
$f \in {\cal RM}^{\sf P}$. Suppose $f$ has 
I/O balance and time-complexity $\le T$ for some polynomial $T$.
To show that $f$ has an inverse with elementary recursive I/O balance and
time-complexity, let $y \in {\sf Im}(f)$ and consider all words $x$ of 
length $|x| \le T(|y|)$; for each such $x$ we test whether 
$x \in {\sf Dom}(f)$, and (if so) we compute $f(x)$, in time $\le T(|x|)$.
On input $y$ we output the first $x$ in length-lexicographic order such that 
$f(x) = y$.
All this takes time $\le |A|^{\ell} \cdot T(\ell)$, where $\ell$ is the
minimum length of $x \in f^{-1}(y)$; so $\ell \le T(|y|)$ (by I/O-balance).  
The bound $\, \tau(|y|) = |A|^{T(|y|)} \cdot T(T(|y|)) \,$ is elementary
recursive, and testing whether $y \in {\sf Im}(f)$ is also elementary
recursive, since ${\sf Im}(f) \in {\sf NP} \subset {\sf E}_{\sf 3}$.
So $f$ has an inverse $f'$ with elementary recursive I/O balance and
time-complexity.

Let $\tau(n)$ be a fully time-constructible elementary recursive upper bound
on $|A|^{T(n)} \cdot T(T(n))$ and on the time it takes to test whether
$y \in {\sf Im}(f)$ (when $|y| = n$). The function $n \mapsto 2^n$ is fully
time-constructible, and if $T(n) = a \cdot (n^d +1)$ then $T(T(n))$ has an
upper bound that has that form too. Moreover, the product of fully
time-constructible functions is fully time-constructible.

So $f'$ has balance and time-complexity bounded by $\tau$.
We use Lemma \ref{recTimeConstr}(2) in the same way as in the proof of
regularity of ${\cal M}_{\sf end}^{\sf P}$ (Theorem \ref{M_is_reg}); we
pad $f'$ by taking the restriction $F'$ of $f'$ to

\smallskip

 \ \ \  \ \ \ ${\sf Dom}(F') \ = \ $
$\bigcup_{y \in {\sf domC}(f')} y \, A^{|y| \cdot \tau(|y|)} \, A^*$.

\smallskip

\noindent And we restrict $f$ to $\, F = f \circ F' \circ f$
$= {\sf id}_{{\sf Dom}(F')} \circ f$.
Then $F F' F = F$, and $F$ and $F'$ have linear time-complexity and balance
(by Lemma \ref{recTimeConstr}(2)), hence $F, F' \in {\cal RM}^{\sf P}$.
Moreover, $f \equiv_{E3} F$ and $f' \equiv_{E3} F'$, since
$A^{|y| \cdot \tau(|y|)}$ is a maximal prefix code and
$n \to n \cdot \tau(n)$ is elementary recursive.
So $[F']$ is an inverse of $[f]$ ($ = [F]$) in ${\cal M}^{\sf P}_{\sf E3}$.  
 \ \ \ $\Box$


\section{Inverses in ${\cal M}^{\sf P}_{\sf poly}$}

In this section we study the regular elements of ${\cal RM}^{\sf P}$ and
of ${\cal M}^{\sf P}_{\sf poly}$, and we eventually show that 
${\cal RM}^{\sf P}$ is regular iff ${\cal M}^{\sf P}_{\sf poly}$ is 
regular.

\subsection{Properties of inverses in ${\cal RM}^{\sf P}$ }

Here are a few useful facts about inverses that were not proved in 
\cite{s1f}, \cite{BiInfGen}.

\begin{lem} \label{f_inv_imc}
 \ For any right-ideal morphism $g$: \,
$g^{-1}({\sf imC}(g)) \subseteq {\sf domC}(g)$.
\end{lem}
{\bf Proof.}  If $x \in g^{-1}({\sf imC}(g))$ ($\subseteq {\sf Dom}(g)$),
then $x = p w$ for some $p \in {\sf domC}(g)$ and $w \in A^*$; hence,
$g(x) = g(p) \ w$, and  $g(x) \in {\sf imC}(g)$.
Since $g(p) \ w \in {\sf imC}(g)$ and $g(p) \in {\sf Im}(g)$, we have
$g(p) \ w = g(p)$ (since ${\sf imC}(g)$ is a prefix code, and
$g(p) \, w \in {\sf imC}(g)$ cannot have a strict prefix in the right ideal
generated by ${\sf imC}(g)$).
So, $w = \varepsilon$, hence $x = pw = p \in {\sf domC}(g)$.
 \ \ \ $\Box$

\begin{lem} \label{injimC}
 \ For every right-ideal morphism $g$ we have: 

\noindent (1) \ $\, {\sf imC}(g) \subseteq g({\sf domC}(g))$.

\noindent (2) \ If $g$ is injective then 
$\, {\sf imC}(g) = g({\sf domC}(g))$.

\noindent (3) \ If $g'$ is an inverse of $g$ and if
$\, {\sf Dom}(g') = {\sf Im}(g)$, then the inverse $g'$ is injective. 
\end{lem}
{\bf Proof.} (1) By applying $g$ to the inclusion in Lemma \ref{f_inv_imc}
we obtain: $g g^{-1}({\sf imC}(g)) \subseteq g({\sf domC}(g))$; since
$g g^{-1} = {\sf id}_{{\sf Im}(g)}$, the result follows. 

\noindent (2) When $g$ is injective, let $g(x) \in g({\sf domC}(g))$ with
$x \in {\sf domC}(g)$. Then $g(x) \in {\sf Im}(g) = {\sf imC}(g) \, A^*$;
so $g(x) = u v$ for some $u \in {\sf imC}(g), \, v \in A^*$.
Let $z \in {\sf Dom}(g)$ be such that $g(z) = u$; then $z = st$ for some
$s \in {\sf domC}(g), \, t \in A^*$. Hence,
$g(x) = uv = g(z) \ v = g(z v) = g(stv)$. Since $g$ is injective, this
implies that $x = stv$, so $s$ and $x$ are prefix-comparable. But then
$s = x$, since $x$ and $s$ belong to the prefix code ${\sf domC}(g)$.
Therefore $t = v = \varepsilon$. It follows that
$g(x) = u v = u \in {\sf imC}(g)$, so $g(x) \in {\sf imC}(g)$.
Thus, $g({\sf domC}(g)) \subseteq {\sf imC}(g)$ when $g$ is injective.
 
\noindent (3) For all $y_1, y_2 \in {\sf Dom}(g') = {\sf Im}(g)$, 
$g'(y_i) \in g^{-1}(y_i)$ ($i = 1 , 2$). If $y_1 \neq y_2$ then 
$g^{-1}(y_1)$ is disjoint from $g^{-1}(y_2)$, so $g'(y_1) \neq g'(y_2)$. 
 \ \ \ $\Box$

\begin{pro} \label{injInverse}
 \ If $f \in {\cal RM}^{\sf P}$ (or $\in {\sf fP}$) is regular then $f$ has 
an {\em injective} inverse $f' \in {\cal RM}^{\sf P}$ (respectively $\in$ 
${\sf fP}$) with the additional property that 
$\, {\sf Dom}(f') = {\sf Im}(f)$.
\end{pro}
{\bf Proof.}  Let $F' \in {\cal RM}^{\sf P}$ (or $\in {\sf fP}$) be an
inverse of $f$, so ${\sf Im}(f) \subseteq {\sf Dom}(F')$.
Since $f$ is regular we know (by Prop.\ 1.9 in \cite{s1f}) that
${\sf Im}(f)$ is in {\sf P}. Hence the restriction $f' = F'|_{{\sf Im}(f)}$
belongs to ${\cal RM}^{\sf P}$ (respectively to {\sf fP}). Moreover, the
restriction of an inverse of $f$ to ${\sf Im}(f)$ is always an injective
inverse of $f$ (by Lemma \ref{injimC}(3)).
 \ \ \ $\Box$

\medskip

\noindent As a consequence of Prop.\ \ref{injInverse} we have:

\begin{cor} \label{injinvDclass}
 \ Every regular ${\cal D}$-class of {\sf fP} and of ${\cal RM}^{\sf P}$
contains injective partial functions.          \ \ \ $\Box$
\end{cor}
Question: Do non-regular ${\cal D}$-classes also contain injective partial
functions (that are thus not regular)?  This is equivalent to the existence 
of injective one-way functions.

\begin{lem} \label{invofPrefcode} \hspace{-.14in} {\bf .}

\noindent (1) \ For every right-ideal morphism $f$: $A^* \to A^*$ and every
prefix code $P \subset A^*$ we have:  $f^{-1}(P)$ is a prefix code, and
 $f^{-1}(P) \, A^* \ \subseteq \ f^{-1}(P A^*)$.

\smallskip

\noindent (2) \ There exists  $f \in {\cal RM}^{\sf P}$ and a prefix code
$P$ such that \ $f^{-1}(P) \ A^* \neq f^{-1}(P \,  A^*)$.

\smallskip

\noindent (3) \ There exists  $f \in {\cal RM}^{\sf P}$ and a prefix code
$P$ such that $f(P)$ is not a prefix code.
\end{lem}
{\bf Proof.} (1) If $x_1$ is a prefix of $x_2 = x_1 u$, with
$x_1, x_2 \in f^{-1}(P)$, then $f(x_1)$ and $f(x_2) = f(x_1) \, u$ both
belong to the prefix code $P$, hence $f(x_1) = f(x_1) \, u$, hence
$u = \varepsilon$. Now $x_2 = x_1 u$ implies $x_2 = x_1$.
So, in $f^{-1}(P)$, prefix-related words are equal, hence $f^{-1}(P)$ is a
prefix code.

Obviously, $f^{-1}(P) \subseteq f^{-1}(P A^*)$. Moreover, $PA^*$ is a
right ideal, hence $f^{-1}(P A^*)$ is a right ideal. Therefore,
$f^{-1}(P) \, A^* \subseteq f^{-1}(P A^*) \ A^* = f^{-1}(P A^*)$.

\noindent (2) Example: Let $f(0^n 1) = 0^n$ for all $n \ge 0$, with
${\sf domC}(f) = 0^* 1$, and ${\sf imC}(f) = \{ \varepsilon\}$.
Let $P = \{ \varepsilon\}$. Then
$f^{-1}(P) \ A^* = f^{-1}(\{\varepsilon\}) \ \{0,1\}^* = 1 \, \{0,1\}^*$,
and $f^{-1}(P \, A^*) = f^{-1}(\{0,1\}^*) = 0^* 1 \, \{0,1\}^*$.

\noindent (3) Example: For $f$ as in (2), let $P = 0^* 1$.
We obtain $f(P) = 0^*$, which is not a prefix code.
  \ \ \ $\Box$

\begin{defn} {\bf (normal morphism).}  \label{normal}
 \ A right-ideal morphism $f$ is called {\em normal} \, iff
 \, $f({\sf domC}(f)) = {\sf imC}(f)$.
\end{defn}
Thus, $f$ is normal iff its restriction to ${\sf domC}(f)$ maps
into (hence onto) ${\sf imC}(f)$; in other words, $f$ is entirely defined
by the way it relates ${\sf domC}(f)$ to ${\sf imC}(f)$. On the other 
hand, a non-normal right-ideal morphism $g$ will map ${\sf domC}(g)$ to a 
larger set than  ${\sf imC}(g)$, i.e.,  ${\sf imC}(g) \subsetneqq$
$g({\sf domC}(g))$.  

\medskip

\noindent {\bf Examples} of normal and non-normal right-ideal morphisms:

Every injective right-ideal morphism is normal (by Lemma \ref{injimC}).

The encodings of the elements of {\sf fP} are normal; we saw near the 
beginning of the Introduction that for all $f \in {\sf fP}$, $f^C$ is normal.

The following is a {\it non-normal regular} element of ${\cal RM}^{\sf P}$: 
Let ${\sf domC}(g) = 0 \, 1^*$; and let \ $g(0^n 1 \, w) = 0^n \, w \, $  
for all $n \ge 0$ and all $w \in \{0,1\}^*$.
Then ${\sf imC}(g) = \{ \varepsilon\} \, \neq \, 0^* = g({\sf domC}(g)$.
So in this example, ${\sf imC}(g)$ and $g({\sf domC}(g)$ are extremely
different.

\begin{lem} \label{equivNormal}
 \ A right-ideal morphism $f$ is normal \ iff
 \ ${\sf domC}(f)$ $=$ $f^{-1}({\sf imC}(f))$.
\end{lem}
{\bf Proof.} The right-to-left implication is trivial since
$f f^{-1} = {\bf 1}_{{\sf Im}(f)}$.

Conversely, let us assume normality, i.e.,
$f({\sf domC}(f)) = {\sf imC}(f)$. Then $f^{-1}({\sf imC}(f))$ $\subseteq$
${\sf domC}(f)$ by Lemma \ref{f_inv_imc}. To prove that
${\sf domC}(f) \subseteq f^{-1}({\sf imC}(f))$, let $x \in {\sf domC}(f)$.
Then $f(x) \in f({\sf domC}(f)) = {\sf imC}(f)$; the latter equality holds
by the assumption of normality.  So, $f(x) \in {\sf imC}(f)$, hence
$x \in f^{-1}({\sf imC}(f))$. Thus, ${\sf domC}(f)$ $\subseteq$
$f^{-1}({\sf imC}(f))$.
 \ \ \ $\Box$

\begin{pro} \label{notClosedComp}
 \ (1) The set $\, \{f \in {\cal RM}^{\sf P}: f$ {\rm is normal}$\} \,$  is 
{\em not} closed under composition, i.e., it is not a submonoid.
In fact, there exist regular normal elements of ${\cal RM}^{\sf P}$ whose 
composite is regular but not normal.

\noindent (2) There exist a regular  normal $g \in {\cal RM}^{\sf P}$ and a 
prefix code $P \subset {\sf Dom}(g)$ with $P \equiv_{\sf poly} {\sf domC}(g)$, 
such that $g(P)$ is not a prefix code.
\end{pro}
{\bf Proof.} (1) This is shown by the following example. Let
$f, g \in {\cal RM}^{\sf P}$ be defined by
$\, {\sf domC}(f) = \{0,1\} \, $ and $\, f(0) = 0$, $f(1) = 10$;
$\, {\sf domC}(g) = \{0,1\} \, $ and $\, g(0) = g(1) = 0$.
Then $f({\sf domC}(f)) = {\sf imC}(f) = \{0, 10\}$, and
$g({\sf domC}(g)) = {\sf imC}(g) = \{0\}$; so, $f$ and $g$ are normal.
Now, ${\sf domC}(gf) = \{0, 1\}$ and $gf(0) = 0$, $\, gf(1) = 00$.  
So, $gf({\sf domC}(gf)) = \{0, 00\}$, which is not a prefix code; and
${\sf imC}(gf) = \{0\}$. Thus, $gf$ is not normal.

\noindent (2) Take $g$ as above, and $P = \{0, 10, 11\}$. Then 
$g(P) = \{0, 00, 01\}$, which is not a prefix code.
 \ \ \ $\Box$.


\bigskip

\noindent {\bf Miscellaneous}

\medskip

\noindent The remaining Definition and Facts of this subsection will not 
be used in the rest of the paper.

\medskip


\noindent {\bf Definition 5.M1}     
{\it \ For any right-ideal morphism $f$: $A^* \to A^*$, the {\em normalization}
$f_N$ of $f$ is the restriction of $f$ to $f^{-1}({\sf imC}(f)) \ A^*$. So,
${\sf domC}(f_N) = f^{-1}({\sf imC}(f))$ ($\subseteq {\sf domC}(f)$).
}   

\medskip

\noindent
Then $f_N$ is normal: $f_N({\sf domC}(f_N) = f({\sf domC}(f_N))$ $=$
$f(f^{-1}({\sf imC}(f)) = {\sf imC}(f)$, and this is equal to
${\sf imC}(f_N)$, by Prop.\ 5.M2 (next). 
Moreover, $f^{-1}({\sf imC}(f))$ is a prefix code (by Lemma
\ref{invofPrefcode}(1)).

Note that $f_N \in {\cal RM}^{\sf P}$ iff $f^{-1}({\sf imC}(f)) \in {\sf P}$.
Indeed, if $f_N \in {\cal RM}^{\sf P}$ then ${\sf domC}(f_N) \in {\sf P}$;
and if ${\sf domC}(f_N) \in {\sf P}$ then the restriction of $f$
($\in {\cal RM}^{\sf P}$) is in ${\cal RM}^{\sf P}$.

We conjecture that $f_N$ is not always in ${\cal RM}^{\sf P}$ when
$f \in {\cal RM}^{\sf P}$; this conjecture is motivated by Prop.\
5.M3 below.  

\medskip

\noindent {\bf Proposition 5.M2}   
{\it 
\ For any right-ideal morphism $f$ and its normalization $f_N$, we have:
 \ ${\sf Im}(f_N) = {\sf Im}(f)$, and ${\sf imC}(f_N) = {\sf imC}(f)$.
}       

\smallskip

\noindent {\bf Proof.} Obviously, ${\sf Im}(f_N) \subseteq {\sf Im}(f)$.
Conversely, let $y \in {\sf Im}(f)$, so $y = q w$ for some
$q \in {\sf imC}(f)$ and $w \in A^*$. Then $q = f(p)$ for some
$p \in {\sf domC}(f) \cap f^{-1}({\sf imC}(f)) = {\sf domC}(f_N)$, hence
$y = f(p) \, w = f(pw) \in {\sf Im}(f_N)$.
 \ \ \ $\Box$

\bigskip

\noindent We know that for all $f \in {\sf fP}$, ${\sf Dom}(f)$ is in
{\sf P} and ${\sf Im}(f)$ is in {\sf NP} (Prop.\ 1.9 in \cite{s1f}). What
can be said about the complexity of ${\sf imC}(f)$?
The complexity class ${\sf DP}$ ($\subseteq \Delta_2^{\sf P}$ $=$
${\sf P}^{\sf NP}$) is defined by

\smallskip

 \ \ \  \ \ \ ${\sf DP} \ = \ \{L_1 - L_0 \ : \ L_1, L_0 \in {\sf NP}\}$.

\smallskip

\noindent Obviously, ${\sf NP} \cup {\sf coNP} \subseteq {\sf DP}$.
There exist {\sf DP}-complete problems, e.g., the following:
 \ ${\sf critical3SAT} \ = \ \{\beta : \beta$ is a boolean
formula in {\sf 3cnf} that is not satisfiable, but for every clause $c$ in 
$\beta$, the removal of $c$ results in a boolean formula $\beta - \{c\}$
that is satisfiable\}. See e.g.\ \cite{Papadim}.

\medskip

\noindent {\bf Proposition 5.M3}   
{\it
 \ For all $f \in {\sf fP}$, \ ${\sf imC}(f)$ is in {\sf DP}, and
when $f$ is normal then ${\sf imC}(f)$ is in {\sf NP}.
}   

\smallskip

\noindent
{\bf Proof.} We have $x \in {\sf imC}(f)$ iff the following hold:
(1) $x \in {\sf Im}(f)$, and (2) for every strict prefix $p$ of $x$:
$p \not\in {\sf Im}(f)$.  The second condition is equivalent to
${\sf not}(\exists p)[ \, p <_{\sf pref} x$ \ and \
$p \in {\sf Im}(f) \, ]$.  Hence, ${\sf imC}(f) \in {\sf DP}$.

However since ${\sf domC}(f)$ is in {\sf P}, $f({\sf domC}(f))$ is in
{\sf NP}; so when $f$ is normal then ${\sf imC}(f)$ ($= f({\sf domC}(f))$)
is in {\sf NP}.
 \ \ \ $\Box$

\medskip

Normalization works well with inverses:

\medskip

\noindent {\bf Proposition 5.M4}   
{\it
 \ For any right-ideal morphism $f$ and its normalization $f_N$ we have:

\noindent (1) \ If $f'_1$ is an inverse of $f_N$ then $f'_1$ is also an
inverse of $f$.

\noindent (2) \ If $f'$ is any inverse of $f$ then the restriction $F'$ of
$f'$ to ${\sf Im}(f)$ is an injective (hence normal) inverse of $f$.
Moreover, if $f, f' \in {\cal RM}^{\sf P}$ then $F' \in {\cal RM}^{\sf P}$.
}   

\smallskip

\noindent 
{\bf Proof.} (1) Let us show that $f_1'$ is an inverse of $f$. Since
${\sf Im}(f) = {\sf Im}(f_N)$, for all $x \in {\sf Dom}(f)$ there exists
$x_1 \in {\sf Dom}(f_N)$ such that $f(x) = f_N(x_1)$.
Now, $f(x) = f_N(x_1) = f_N f_1' f_N(x_1) = f_N f_1' f(x) \subseteq $
$f f_1' f(x)$; the latter holds since $f_N$ is a restriction of $f$. But
since $f f_1' f$ is a function, and $f_N f_1' f(x)$ is defined, we have
$f_N f_1' f(x) = f f_1' f(x)$.  Thus, $f(x) = f f_1' f(x)$.

\noindent
(2) The only part of the domain of an inverse of $f$ that matters (in the
relation $f f' f = f$) is ${\sf Im}(f)$ (which is always a subset of
${\sf Dom}(f')$). So the restriction $F'$ of $f'$ to ${\sf Im}(f)$ is an
inverse of $f$. For any inverse $f'$ of $f$ we have:
$f'({\sf imC}(f)) \subseteq f^{-1}({\sf imC}(f))$, which is a prefix code
by Lemma \ref{invofPrefcode}(1); so, $F'$ is normal.
For any inverse, the restriction to ${\sf Im}(f)$ is injective, since for
all $y_1 \neq y_2$ in ${\sf Im}(f)$, $f^{-1}(y_1)$ and $f^{-1}(y_2)$ are
disjoint.

If $f, f' \in {\cal RM}_2^{\sf P}$ then $f$ is regular, so by Prop.\ 1.9 in
\cite{s1f}, ${\sf Im}(f)$ is in {\sf P}. Hence the restriction of $f'$ to
${\sf Im}(f)$ is in ${\cal RM}_2^{\sf P}$.
 \ \ \ $\Box$.

\bigskip


We know (Prop.\ 6.1 in \cite{s1f}) that every element of {\sf fP}, and in
particular, every element of ${\cal RM}^{\sf P}$, has an inverse in
${\sf fP}^{\sf NP}$.
We show next that every element of ${\cal RM}^{\sf P}$ has an inverse in
${\cal RM}^{\bf NP}$. We first extend Prop.\ 2.6 of \cite{s1f} to
${\sf fP}^{\sf NP}$ and to ${\cal RM}^{\bf NP}$.

\medskip

\noindent {\bf Lemma 5.M5}   
{\it
 \ If an element $f \in {\cal RM}^{\sf P}$ has an inverse in {\sf fP} (or
in ${\sf fP}^{\sf NP}$), then $f$ also has an inverse in
${\cal RM}^{\sf P}$ (respectively in ${\cal RM}^{\bf NP}$).

Moreover, this inverse in $ {\cal RM}^{\sf P}$ (resp.\ ${\sf fP}^{\sf NP}$)
can be chosen to be injective (and hence normal).
}   

\smallskip

\noindent 
{\bf Proof.} If $f$ has an inverse in {\sf fP} then the result was proved
in Prop.\ 2.6 of \cite{s1f}.

Let $f_0' \in {\sf fP}^{\sf NP}$ be an inverse of $f$; we want to
construct an inverse $f'$ of $f$ that belongs to ${\cal RM}^{\bf NP}$.
We know (Prop.\ 1.9 of \cite{s1f}) that ${\sf Im}(f)$ is in {\sf NP}.
Hence we can restrict $f_0'$ to ${\sf Im}(f)$, i.e.,
${\sf Dom}(f_0') = {\sf Im}(f)$.
We proceed to define $f'(y)$ for $y \in {\sf Im}(f)$.

First, we compute the shortest prefix $p$ of $y$ that satisfies
$p \in {\sf Dom}(f_0') = {\sf Im}(f)$. Since  ${\sf Im}(f) \in {\sf NP}$,
this can be done in polynomial time with calls to an {\sf NP} oracle.
Now, $y = p \, z$ for some  string $z$.

Second, we define $f'(y) \ = \ f'_0(p) \ z$, \ where $p$ and $z$ are as
above.  Thus, $f'$ is a right-ideal morphism.

\smallskip

Let us verify that $f'$ has the claimed properties.
Clearly, $f'$ is computable in polynomial time with calls to an {\sf NP}
oracle, and is polynomially balanced (the latter following from the fact
that $f'$ is an inverse of $f$, which we prove next); thus,
$f'$ is a right-ideal morphism in ${\sf fP}^{\sf NP}$, so
$f' \in {\cal RM}^{\bf NP}$.
To prove that $f'$ is an inverse of $f$, let $x \in {\sf Dom}(f)$.
Then $f (f' (f(x))) = f (f'(p \, z))$, where $y = f(x) = p \, z$, and $p$
is the shortest prefix of $y$ such that $p \in {\sf Im}(f)$.
Then, $f'(p \, z) = f_0'(p) \, z$, by the definition of $f'$.
Then, since $f$ is a right-ideal morphism, $f(f_0'(p) \, z)$
$= f(f_0'(p)) \, z = p \, z$ (the latter since $f'_0$ is an inverse of $f$,
and since $p \in {\sf Im}(f)$).
Hence, $ff'|_{{\sf Im}(f)} = {\bf 1}_{{\sf Im}(f)}$. Thus, $f'$ is an
inverse of $f$.

Note that since ${\sf Dom}(f') = {\sf Im}(f)$, the inverse $f'$ described
above is injective. Indeed, if ${\sf Dom}(f') = {\sf Im}(f)$ then
$f' = f'|_{{\sf Im}(f)}$, so $ff' = {\bf 1}_{{\sf Dom}(f')}$,
which implies that $f'$ is injective (hence normal by Lemma \ref{injimC}).
 \ \ \ $\Box$

\medskip

\noindent {\bf Proposition 5.M6}   
{\it \ Every element of ${\cal RM}^{\sf P}$ has an inverse in
${\cal RM}^{\bf NP}$, and this inverse can be chosen to be injective
(and hence normal).
}       

\smallskip

\noindent {\bf Proof.} By Prop.\ 6.1 in \cite{s1f}, every element of
${\cal RM}^{\sf P}$ has an inverse in ${\sf fP}^{\sf NP}$.
The result then follows from Lemma 5.M5.  
 \ \ \ $\Box$


\bigskip

\noindent {\bf Proposition 5.M7}   
{\it
 \ Let $f_0 \in {\cal RM}^{\sf P}$, and let $f$ be any right-ideal morphism
such that $f \equiv_{\sf poly} f_0$.
Then $f \in {\cal RM}^{\sf P}$ \, iff \, ${\sf Dom}(f) \in {\sf P}$.

Hence, if $f_0 \in {\cal RM}^{\sf P}$, and ${\sf Dom}(f) \in {\sf P}$,
and $f \not\in {\cal RM}^{\sf P}$, then $f \not\equiv_{\sf poly} f_0$.
}   

\smallskip

\noindent {\bf Proof.} We know that for all $f \in {\cal RM}^{\sf P}$,
${\sf Dom}(f) \in {\sf P}$.
For the converse, if $x \in {\sf Dom}(f)$ (which can be checked in polynomial
time), then either $x \in {\sf Dom}(f_0)$ or $x$ is a prefix of a word
$xu \in$ ${\sf domC}(f_0)$. If $x \in {\sf Dom}(f_0)$ we can immediately
compute $f_0(x)$ ($= f(x)$), using the polynomial-time algorithm of $f_0$.

If $xu \in {\sf domC}(f_0)$ for some $u \in A^*$, we can compute $f_0(xu)$
in polynomial time (as a function of $|xu|$). Here, $u$ is the shortest word
such that $xu \in {\sf domC}(f_0)$. So, $|u|$ is polynomially bounded in terms
of $|x|$ (because ${\sf domC}(f) \equiv_{\sf poly} {\sf domC}(f_0)$).
Therefore, the computation of $f_0(xu)$ takes polynomial time as a function
of $|x|$.

Also, $f_0(xu) = f(xu) = f(x) \ u$; so we obtain $f(x)$ by removing the
suffix $u$ from $f_0(xu)$; we know $u$, since it is the shortest word such
that $xu \in {\sf domC}(f_0)$ (and ${\sf domC}(f_0) \in {\sf P}$ when
$f \in {\cal RM}^{\sf P}$).
 \ \ \ $\Box$


\bigskip

\noindent {\bf Proposition 5.M8}   

{\it
\noindent (1) There exist prefix codes $P_1, P_0 \subset A^*$ such that
$P_1 \equiv_{\sf poly} P_0$, and $P_0 \in {\sf P}$, but
$P_1 \not\in {\sf P}$. The prefix code $P_1$ can be chosen to have any
complexity above polynomial, or to be undecidable;
if ${\sf P} \neq {\sf NP}$ then $P_1$ can be chosen in {\sf DP}.

\smallskip

\noindent
(2) There exist right ideal morphisms $f_1, f_0$ such that
$f_1 \equiv_{\sf poly} f_0$, and $f_0 \in {\cal RM}^{\sf P}$, but
$f_1  \not\in {\cal RM}^{\sf P}$.
If ${\sf P} \neq {\sf NP}$ then $f_1$ can be chosen in
${\cal RM}^{\bf NP}$.
}   

\smallskip

\noindent 
{\bf Proof.} (1) We construct a family of examples. Let $L \subset A^*$ be
any set that is not in {\sf P}.  Let

\smallskip

 \ \ \ \ \ $P_0 \ = \ \{00, 01\}^* \, 11$, \ \ and

\smallskip

 \ \ \ \ \
$P_1 \ = \{ {\sf code}(x) \ 11 \ : \ x \in L\}$ $\cup$
    $\{ {\sf code}(x) \ 11 0 \ : \ x \not\in L \}$ $\cup$
    $\{ {\sf code}(x) \ 11 1 \ : \ x \not\in L \}$.

\smallskip

\noindent Then $P_1, P_0$ are prefix codes, $P_1 \equiv_{\sf poly} P_0$,
and $P_0 \in {\sf P}$. But $P_1 \not\in {\sf P}$ since $L$ is
polynomial-time reducible to $P_1$.

\smallskip

\noindent
(2) Let $f_0, f_1$ be the identity map restricted to $P_0 A^*$, respectively
$P_1 A^*$ (with $P_0, P_1$ as above). Then $f_1 \equiv_{\sf poly} f_0$, and
$f_0 \in {\cal RM}^{\sf P}$; but $f_1  \not\in {\cal RM}^{\sf P}$ since
${\sf domC}(f_1) = P_1 \not\in {\sf P}$.  If ${\sf P} \neq {\sf NP}$ then
$L$ can be chosen in ${\sf NP} - {\sf P}$, and then $f_1 \not\in$
${\cal RM}^{\sf P}$. 
 \ \ \ $\Box$

\bigskip

\noindent {\bf Question:} Assuming ${\sf P} \neq {\sf NP}$, is there
$F \in {\cal M}^{\sf P}_{\sf poly}$ such that for all
$f \in F:$ \ ${\sf Dom}(f) \not\in {\sf P}$?

\subsection{${\cal M}^{\sf P}_{\sf poly}$ vs.\ ${\cal RM}^{\sf P}$, 
            regarding regularity }

It is obvious that if ${\cal RM}^{\sf P}$ is regular then
${\cal M}^{\sf P}_{\sf poly}$ (= ${\cal RM}^{\sf P}\!/\! \equiv_{\sf poly}$)
is regular, being a homomorphic image of ${\cal RM}^{\sf P}$. The converse 
is also true, but the proof is not obvious, mainly because of the existence 
on non-normal functions in ${\cal RM}^{\sf P}$.  Many of the results of this 
sub-section hold for ${\cal M}^{\sf P}_{\cal T}$ (where $\cal T$ is any 
family of functions as in  Def.\ \ref{equiv_boundDEF}).

\begin{lem} \label{Im_endequiv}
 \ If $f_0, f$ are right-ideal morphisms with $f_0 \equiv_{\sf end} f$ 
and $f_0 \subseteq f$, then 
$\, {\sf imC}(f_0) \equiv_{\sf end} {\sf imC}(f)$.

If ${\sf poly} \subseteq {\cal T}$ and $f_0, f \in {\cal RM}^{\sf P}$ 
satisfy $f_0 \equiv_{\cal T} f$ and $f_0 \subseteq f$, then 
$\, {\sf imC}(f_0) \equiv_{\cal T} {\sf imC}(f)$.
\end{lem}
Compare with Lemma \ref{bdequ_Im}.

\smallskip

\noindent {\bf Proof.}  Since $f_0 \subseteq f$ we also have  
${\sf Im}(f_0)$ $\subseteq$ ${\sf Im}(f)$.
Suppose a right ideal $R$ intersects ${\sf Im}(f)$; so there exists
$x \in {\sf Dom}(f)$ such that $f(x) \in R \cap {\sf Im}(f)$. Hence,
$x \in f^{-1}(R \cap {\sf Im}(f)) = f^{-1}(R) \cap  {\sf Dom}(f)$; so
the right ideal $f^{-1}(R)$ intersects ${\sf Dom}(f)$, hence by
end-equivalence, $f^{-1}(R)$ also intersects ${\sf Dom}(f_0)$. So there
exists $x_0 \in f^{-1}(R) \cap {\sf Dom}(f_0)$, and this implies that
$f(x_0) \in ff^{-1}(R) \cap f({\sf Dom}(f_0)) = R \cap {\sf Im}(f_0)$.
So, $R$ intersects ${\sf Im}(f_0)$.

For the second statement, let $y_0 \in {\sf imC}(f_0)$ and 
$y \in {\sf imC}(f)$ be such that 
$y \le_{\sf pref} y_0 = y w$ (for some $w \in A^*$). We want to show that
$|y_0|$ and $|y|$ are related by some function in ${\cal T}$ that depends 
only on $f$ and $f_0$. 
Since $f_0^{-1}({\sf imC}(f_0)) \subseteq {\sf domC}(f_0)$ and
$f^{-1}({\sf imC}(f)) \subseteq {\sf domC}(f)$ (by Lemma \ref{f_inv_imc}),
there exists $x \in {\sf domC}(f)$ such that $y = f(x)$, and hence
$y_0 = f(x) \, w = f(xw)$; and $xw \in f_0^{-1}({\sf imC}(f_0))$
$\subseteq$ ${\sf domC}(f_0)$. So, $x \in {\sf domC}(f)$ and 
$xw \in {\sf domC}(f_0)$, and $xw \ge_{\sf pref} x$, hence $|x|$ and $|xw|$
are length-related by a function in $\cal T$ (because 
$f_0 \equiv_{\cal T} f$). Moreover, $|f(x)|$ and $|x|$ are polynomially 
related (because of the I/O-balance of $f$), and $|f_0(xw)|$ and $|xw|$ are 
polynomially related (because of the I/O-balance of $f_0$).
Thus, $|y|$ and $|y_0|$ are length-related by a function in $\cal T$.
 \ \ \ $\Box$

\begin{lem} \label{regularity1} 
 \ Let $h, g$ be any right-ideal morphisms such that $hgh \equiv_{\sf bd} h$.
Then \ $hgh \subseteq h$, and \ $hgh \ g \ hgh =  hgh$.
\end{lem} 
{\bf Proof.} For all functions we have 
${\sf Dom}(hgh) \subseteq {\sf Dom}(h)$,
so since $hgh \equiv_{\sf bd} h$, we have $h g h \subseteq h$.
Hence, for all $x \in {\sf Dom}(hgh)$ we have $hgh(x) = h(x)$, and $hg$ is
defined on $h(x)$. Since $hgh(x) = h(x)$ and $hg$ is defined on $h(x)$,
$hg$ is defined  on $hgh(x)$, and we have $hghgh(x) = hgh(x) = h(x)$ for
all $x \in {\sf Dom}(hgh)$. By the same argument, $hg$ is defined on
$hghgh(x)$, on $hgh(x)$, and on $h(x)$, and we have:
$hghghgh(x) = hghgh(x) = hgh(x) = h(x)$. In particular:
$hgh \, g \, hgh(x) =  hgh(x)$ for all $x \in {\sf Dom}(hgh)$.
  \ \ \ $\Box$

\begin{pro} \label{InvAndMutinv} 
 \ Let $F, G$ be any $\, \equiv_{\cal T}$-equivalence classes in 
${\cal RM}^{\sf P}$.  Then we have:
 
\smallskip

\noindent {\bf (Inverses)} 
 \ $FGF = F$ \, iff \, there exist $f \in F$ and $g \in G$ such that
$\, fgf = f$. 

\smallskip

\noindent {\bf (Mutual inverses)} 
 \ $FGF = F$ and $GFG = G$ \, iff \, there exist $f \in F$ and $g \in G$ 
such that $\, fgf = f \, $ and $\, gfg = g \, $.
\end{pro}
{\bf Proof.} 
For the first statement: If $f \in F$ and $g \in G$ satisfy $fgf = f$ then
$FGF = F$ since $\equiv_{\cal T}$ is a congruence, and $F = [f]$, 
$G = [g]$.
Conversely, if $FGF = F$ then for any $h \in F$, $g \in G$, we have 
$hgh \equiv_{\cal T} h$.  Hence, for any $g \in G$, letting 
$f = hgh \in FGF = F$ we have $f g f = f$ (by Lemma \ref{regularity1}).  

For the second statement: The right-to-left implication is obvious since
$\equiv_{\cal T}$ is a congruence. Conversely, $FGF = F$ implies 
$f g_1 f = f$ for some $f \in F$ and $g_1 \in G$ (by the "Inverses" 
statement of the Proposition, that we just proved). 
Let $g = g_1 f g_1 \in G F G = G$.
Then $fgf = f g_1 f g_1 f = f g_1 f = f$, and 
$g f g = g_1 f g_1 f g_1 f g_1 = g_1 f g_1 f g_1 = g_1 f g_1 = g$.
 \ \ \ $\Box$

\begin{lem} \label{regularity1_1} 
 \ Let $P_0, P_1 \subset A^*$ be prefix codes with 
$P_0 A^* \subseteq P_1 A^*$ and $P_0 \equiv_{\cal T} P_1$;  let 
$\, \tau \in {\cal T} \,$ be the function used for 
$P_0 \equiv_{\cal T} P_1$.  Then for every $y_1 \in P_1$ and every 
$t \in A^*$ with $|t| \ge \tau(|y_1|)$: \ $y_1 t \in P_0 A^*$.
\end{lem}
{\bf Proof.} By definition, $\equiv_{\cal T}$ implies $\equiv_{\sf bd}$, so 
$P_0 A^{\omega}$ $=$ $P_1 A^{\omega}$ (by Prop.\ \ref{equiv_boundEndsPROP}). 
Hence for all $y_1 \in P_1$, $t \in A^*$, and $w \in  A^{\omega}$: the end 
$y_1 t w$ intersects $P_0$, i.e., some prefix of $y_1 t w$ is in $P_0$. 
If $|t| \ge \tau(|y_1|)$ then this prefix is a prefix of $y_1 t$ (by the 
definition of $\equiv_{\cal T}$ and the choice of $\tau$). Hence, 
$y_1 t \in P_0 A^*$. 
 \ \ \ $\Box$

\bigskip

\noindent {\bf Terminology:} \ A {\em normal inverse} of a right-ideal 
morphism $f$ is any normal right-ideal morphism $f'$ (i.e.,
$f'({\sf domC}(f')) = {\sf imC}(f')$, by Def.\ \ref{normal}) such that $f'$ 
is an inverse of $f$.

\begin{lem} \label{invSub}
 \ Let $g, f$ be right-ideal morphisms such that $g \subseteq f$, and let
$f'$ be any inverse of $f$ such that
$f'({\sf Im}(g)) \subseteq {\sf Dom}(g)$.
Then $f'$ is also an inverse of $g$.
\end{lem}
{\bf Proof.} For $x \in {\sf Dom}(g)$, $f' g(x)$ is defined, since
$g \subseteq f$ and since $f(x)$ ($ = g(x)$) is defined, and $f'$ is 
defined on $f(x)$. And $g f' g(x)$ is
defined since $f'({\sf Im}(g)) \subseteq {\sf Dom}(g)$.
Hence, $g f' g(x) = f f' f(x) = f(x) = g(x)$, since $g \subseteq f$.
 \ \ \ $\Box$

\bigskip

\noindent Lemma \ref{regularityNorm2} and Theorem \ref{regofRMequ} below 
are only proved for ${\cal M}^{\sf P}_{\sf poly}$; 
it is not clear for what other ${\cal M}^{\sf P}_{\cal T}$ they hold.

For the next Lemma, recall that if $g \in {\cal RM}^{\sf P}$ is 
regular then $g$ has a normal inverse $g' \in {\cal RM}^{\sf P}$; in fact, 
we can choose $g'$ to be injective such that ${\sf Dom}(g') = {\sf Im}(g)$ 
(see Prop.\ \ref{injInverse}, Lemma \ref{injimC}, and Def.\ \ref{normal}).

\begin{mlem} {\bf (inverse of a $\, \equiv_{\sf poly}$-equivalent 
extension).}        \label{regularityNorm2} 
 \ Suppose $f, f_0 \in {\cal RM}^{\sf P}$ are such that $f_0 \subseteq f$
and $f_0 \equiv_{\sf poly} f$.
Suppose also that $f$ is {\em normal}. Then we have:

\smallskip

\noindent (1) \, If $f_0$ is regular then $f$ is regular.

\smallskip

\noindent (2) \, For every injective inverse $f'_0 \in {\cal RM}^{\sf P}$ 
of $f_0$ such that $\, {\sf Dom}(f'_0) = {\sf Im}(f_0)$, there exists an 
{\em injective} inverse $f'_1 \in {\cal RM}^{\sf P}$ of $f$ such that 
$\, {\sf domC}(f_0') \equiv_{\sf poly} {\sf domC}(f'_1)$.

\smallskip

\noindent (3) \, Moreover, $f'_1$ is also an inverse of $f_0$.
But $f'_1$ cannot always be chosen to be an extension of the given $f'_0$.
\end{mlem}
{\bf Proof.} (1) follows from (2), since every regular element $f_0$ of
${\cal RM}^{\sf P}$ has an injective inverse $f_0'$ satisfying 
${\sf Dom}(f_0') = {\sf Im}(f_0)$ (by Prop.\ \ref{injInverse} and Lemma 
\ref{injimC}).

\smallskip

\noindent (2) Let $f'_0$ be an inverse of $f_0$ as assumed, hence
$f'_0({\sf imC}(f_0))$ $\subseteq$ $f_0^{-1}({\sf imC}(f_0))$ $\subseteq$
$ {\sf domC}(f_0) \, $ (the latter ``$\subseteq$'' holds by Lemma
\ref{f_inv_imc}).

\smallskip

\noindent {\sf Claim:} If $f$ is normal then for all $y \in {\sf imC}(f)$
and all $t \in A^*$: \ $\, f^{-1}(yt) = f^{-1}(y) \ t$.

\smallskip

\noindent {\sf Proof of Claim:} $[\subseteq]$:
 \ $x \in f^{-1}(yt)$ iff $f(x) = yt$.
Since $x \in {\sf Dom}(f)$ we have $x = p w$ for $p \in {\sf domC}(f)$,
$xw \in A^*$, so $f(x) = yt = f(p) \, w$; hence $f(p)$ and $y$ are
prefix-comparable. By normality, $f(p) \in {\sf imC}(f)$; hence
$f(p) = y \,$ (since $y \in {\sf imC}(f)$ by assumption, and ${\sf imC}(f)$
is a prefix code). Thus, $f(x) = yt = f(p) \, w = y w$, so $w = t$ (since
$y = f(p)$). Hence $x = pw = pt \in f^{-1}(y) \ t$.

\noindent $[\supseteq]$ (this holds also when $f$ is not normal):
 \ $x \in f^{-1}(y) \ t \, $ implies $f(x) \in$
$f(f^{-1}(y) \, t)$.  Since $f^{-1}(y) \subseteq f^{-1}({\sf imC}(f))$
$\subseteq {\sf domC}(f) \, $ (the latter ``$\subseteq$'' holds by Lemma
\ref{f_inv_imc}), we have $f(f^{-1}(y) \, t) = f(f^{-1}(y)) \, t = \{y t\}$.
Hence, $x \in f^{-1}(yt)$.
 \ \ \ {\sf [End, Proof of Claim]}

\smallskip

Now let $y \in {\sf imC}(f)$, and let $t \in A^*$ be any string such that
$y t \in {\sf imC}(f_0)$. Since $f_0 \subseteq f$ and 
$f_0 \equiv_{\sf poly} f$, we have 
${\sf imC}(f_0) \equiv_{\sf poly} {\sf imC}(f)$ (by Lemma \ref{Im_endequiv}); 
hence, $|t| \le q(|y|)$ for some polynomial $q$. And by Lemma 
\ref{regularity1_1} (with ${\cal T} = {\sf poly}$), we can pick $t$ to be 
$t = 0^{q(|y|)}$; then $t$ can be computed from $y$ in polynomial time.
Since $f_0'(yt) \in f_0^{-1}(yt) \subseteq f^{-1}(yt)$, we have by the
Claim:      \ $f_0'(yt) \in f^{-1}(y) \ t$.
Since $f_1'(y)$ should belong to $f^{-1}(y)$, we define:

\smallskip

 \ \ \ $f_1'(y)$ is is prefix of $f_0'(yt)$ obtained by removing the 
       suffix $t$.

\smallskip

\noindent Then, indeed,  $f_1'(y)  \in f^{-1}(y)$.
In general, for all $y \in {\sf imC}(f)$ and all $z \in A^*$, we define 
$f_1'(yz) = f_1'(y) \ z$. Then  $f_1'(yz) \in f^{-1}(yz)$, hence  $f_1'$ is 
an inverse of $f$.  
By construction, ${\sf domC}(f_1') = {\sf imC}(f)$, hence $f_1'$ is 
injective (by Lemma \ref{injimC}(3)).
And $f_1'(y)$ is polynomial-time computable, since $t = 0^{q(|y|)}$ and    
since $f'_0 \in {\cal RM}^{\sf P}$.
Finally, $f_1'$ is polynomially balanced, since  $f_0'$ is polynomially
balanced and $|t| \le q(|y|)$.

By Lemma \ref{Im_endequiv}, ${\sf imC}(f_0) \equiv_{\sf poly} {\sf imC}(f)$. 
Hence, ${\sf domC}(f_0') = {\sf imC}(f_0) \equiv_{\sf poly} {\sf imC}(f)$ 
$=$ ${\sf domC}(f_1')$, so $\, {\sf domC}(f_0')$ $\equiv_{\sf poly}$ 
${\sf domC}(f_1')$.

\smallskip

\noindent (3) By Lemma \ref{invSub}, $f'_1$ is an inverse of $f_0$.
In order to apply Lemma \ref{invSub} we need to check that
$f'_1({\sf Im}(f_0)) \subseteq {\sf Dom}(f_0)$.  For all
$y z \in {\sf imC}(f_0)$ (with $y  \in {\sf imC}(f)$) we have
$f_1'(y z) \in f^{-1}(y z)$; and
$f^{-1}(y z) = f_0^{-1}(y z)$ when $y z \in {\sf imC}(f_0)$. Moreover,
$f_0^{-1}(y z) \subseteq {\sf Dom}(f_0)$.
Hence, $f'_1({\sf imC}(f_0)) \subseteq {\sf Dom}(f_0)$, and thus
$f'_1({\sf Im}(f_0)) \subseteq {\sf Dom}(f_0)$.

By Lemma \ref{nonExt} below, the inverse $f'_0$ of $f_0$ is not necessarily
a restriction of an inverse of $f$. So $f_0'$ is not always  extendable to
an inverse of $f$.
 \ \ \ $\Box$

\begin{lem} \label{nonExt}
 \ There exist $f, g \in {\cal RM}^{\sf P}$ such that $g \subseteq f$,
$g \equiv_{\sf poly} f$, and $g$ is regular, but such that not every
inverse $g' \in {\cal RM}^{\sf P}$ (not even every injective inverse) of 
$g$ is extendable to an inverse of $f$.
\end{lem}
{\bf Proof.} This is illustrated by the following example:

\smallskip

$f(0) = f(1) = 1$, with $\, {\sf domC}(f) = \{0,1\}$, 
 \ ${\sf imC}(f) = \{1\}$; \ and

\smallskip

$g(00) = g(10) = 10, \ g(01) = g(11) = 11$, with ${\sf domC}(g)$
$=$ $\{00,10,01,11\}$, \ ${\sf imC}(g) = \{10, 11\}$.

\smallskip

\noindent Then every inverse $f'$ of $f$ satisfies either $f'(1) = 0$ or
$f'(1) = 1$. In particular, $f$ has two injective inverses with domain code
$\{1\}$ ($= {\sf imC}(f)$), namely $f'_0$ and $f'_1$, given by $f'_0(1) = 0$
and $f'_1(1) = 1$.

And $g$ has four injective inverses with domain code $\{10, 11\}$
($= {\sf imC}(g)$). Two of them, namely $g'_0$ and $g'_0$,
are restrictions of $f'_0$, respectively $f'_1$, defined by
$g'_0(10) = 00, \ g'_0(11) = 01$, and $g'_1(10) = 10, \ g'_1(11) = 11$.
The two other injective inverses of $g$ with domain code $\{10, 11\}$ are
$g'_2$ and $g'_3$, defined by $g'_2(10) = 00, \ g'_2(11) = 10$, and
$g'_3(10) = 10, \ g'_3(11) = 01$. These are not restrictions of inverses
of $f$, since every inverse $f'$ of $f$ satisfies either $f'(1) = 0$ or
$f'(1) = 1$, hence $f'(10) = 10, f'(11) = 11$,
or $f'(10) = 00, f'(11) = 01$; in either case,
$g'_2$, $g'_3$ are not restrictions of $f'$.
 \ \ \ $\Box$

\begin{thm}  \label{regofRMequ} 
 \ The monoid ${\cal M}^{\sf P}_{\sf poly}$ is regular iff 
$\, {\cal RM}^{\sf P}$ is regular.   
\end{thm}
{\bf Proof.} Obviously, if ${\cal RM}^{\sf P}$ is regular then its 
homomorphic image ${\cal M}^{\sf P}_{\sf poly}$ is regular, since 
$\equiv_{\sf poly}$ is a congruence. 

For the converse we will show that if ${\cal M}^{\sf P}_{\sf poly}$ is 
regular then {\sf fP} is regular; the latter implies that 
${\cal RM}^{\sf P}$ is regular (by  Prop.\ 2.6 of \cite{s1f}). 
For any $f \in {\sf fP}$, let $f^C \in {\cal RM}^{\sf P}$ be the encoding 
of $f$, as defined near the beginning of the Introduction. Then $f^C$ is 
normal (see the Examples after Def.\ \ref{normal}).  
Let $F = [f^C] \in {\cal M}^{\sf P}_{\sf poly}$ be the 
$\equiv_{\sf poly}$-class of $f^C$ in ${\cal RM}^{\sf P}$, and let 
$F' \in {\cal M}^{\sf P}_{\sf poly}$ be an inverse of $F$.
A consequence of $F F' F = F$ in ${\cal M}^{\sf P}_{\sf poly}$ is that for 
all $h \in F$ and all $g \in F'$: $hgh \equiv_{\sf poly} h$. 
Then by Lemma \ref{regularity1}, $hgh \in F$ and $hgh$ is regular with
inverse $g \in F'$. Also, $hgh \subseteq h$. Let $h = f^C$, which is normal. 
Then the Lemma \ref{regularityNorm2}(1) applies since  
$hgh \subseteq h$, $hgh \equiv_{\sf poly} h$ (with $h = f^C$), 
$h$ is normal, and $hgh$ is regular. Hence Lemma \ref{regularityNorm2}(1)
implies that $h = f^C$ is regular in ${\cal RM}^{\sf P}$.
Hence by Prop.\ 3.4(2) in \cite{s1f}, $f$ is regular in {\sf fP}.
 \ \ \ $\Box$

\bigskip

\noindent {\bf Comments:} 
The proof of Theorem \ref{regofRMequ} also shows the following fact:
{\em If all normal elements of ${\cal RM}^{\sf P}$ are regular then
${\cal RM}^{\sf P}$ is regular.} Thus the set of all normal elements of
${\cal RM}^{\sf P}$ plays a crucial role. 
It remains an open question whether we have the following element-wise
properties:
Let $F \in {\cal M}^{\sf P}_{\sf poly}$ (hence 
$F \subset {\cal RM}^{\sf P}$); if $F$ is regular in 
${\cal M}^{\sf P}_{\sf poly}$, does that imply that every $f \in F$ is 
regular in ${\cal RM}^{\sf P}$?
Equivalently, let $f_0, f \in {\cal RM}^{\sf P}$ be such that
$f_0 \subseteq f$, $f_0 \equiv_{\sf poly} f$, and $f_0$ is regular; does 
that imply that $f$ is regular? \ Lemma \ref{regularityNorm2}(1) yields
this statement when $f$ is normal.

\section{A non-regular monoid that maps onto ${\cal M}^{\sf P}_{\sf poly}$}

We show that there is a non-regular submonoid of ${\cal RM}^{\sf P}$ that
maps homomorphically onto ${\cal M}^{\sf P}_{\sf poly}$.  
The fact that some non-regular monoid maps onto 
${\cal M}^{\sf P}_{\sf poly}$ is trivial, by itself, because we could use a 
(finitely generated) free monoid for this.  
However, there is a non-regular submonoid ${\cal RM}^{n + o(n)}$ of 
${\cal RM}^{\sf P}$ such that the following monoid homomorphisms (where 
$\nearrow$ is injective) form a commutative diagram:


\hspace{1.63in} ${\cal RM}^{\sf P}$

\hspace{1.4in} $\nearrow$ \hspace{.15in} $\downarrow$

\vspace{-.15in}

\hspace{1.6in}  \hspace{.15in} $\downarrow$

\smallskip

\hspace{0.7in} ${\cal RM}^{n + o(n)}$
 \ $\twoheadrightarrow$ \ ${\cal M}^{\sf P}_{\sf poly}$ 
 
\bigskip

\noindent The construction of ${\cal RM}^{n + o(n)}$ is intuitive, but
we need some definitions.

We will use the classical Landau symbol $o$.  For two total functions
$t_1, t_2$: ${\mathbb N} \to {\mathbb R}_{\ge 0}$ we say that
``{\em $t_1$ is $o(t_2)$}'' iff there exists a total function $\epsilon$:
${\mathbb N} \to {\mathbb R_{\ge 0}}$ such that
$\lim_{n \to \infty} \epsilon(n) = 0$, and for all $n \in {\mathbb N}$,
$\, t_1(n) \le \epsilon(n) \cdot t_2(n)$.
In particular, a total function $t$: ${\mathbb N} \to {\mathbb N}$ is said
to be $n + o(n)$ iff there exists a total function
$\epsilon$: ${\mathbb N} \to {\mathbb R_{\ge 0}}$ such that
$\lim_{n \to \infty} \epsilon(n) = 0$, and for all
$n \in {\mathbb N}$: \ $t(n) \le n + \epsilon(n) \cdot n$.
Since $n + \epsilon(n) \cdot n = (1 + \epsilon(n)) \cdot n$, we can also
write $(1 + o(1)) \cdot n \,$ for $n + o(n)$. (By the definition of the
Landau symbol, a function $t$: ${\mathbb N} \to {\mathbb R}_{\ge 0}$ is
$o(1)$ iff $\lim_{n \to \infty} t(n) = 0$.)

Clearly, the set of total functions ${\mathbb N} \to {\mathbb N}$ that
are $n + o(n)$ is closed under composition.

An {\em ${\cal RM}^{\sf P}$\!-machine} is a multi-tape Turing machine $M$
with a read-only input-tape that contains the input, and with a
write-only output-tape, such that the input-tape head and the output-tape
head never move left.  The machine has an {\em accept state}; when $M$
halts, the content of the output-tape is a valid output iff $M$ is in the
accept state (when $M$ halts in a non-accept state, the output is
undefined). A convention of this sort is necessary, otherwise there is 
always an output (possibly the empty string).
Let $f_M$ denote the input-output function of $M$. We assume that
for every $x \in {\sf Dom}(f_M)$ and every word $z \in A^*$, the
computation of $M$ on input $xz$ has the following property: the
input-tape head does not start reading $z$ until $f_M(x)$ has been written
on the output tape. (To ``read'' a letter $\ell$ means to make a transition
whose input-tape letter is this letter $\ell$.) This is not the complete
definition of an ${\cal RM}^{\sf P}$\!-machine, but that is all we need
here; the details are given at the beginning of Section 2 in
\cite{BiInfGen}.
We define the following submonoid of ${\cal RM}^{\sf P}$:

\medskip

${\cal RM}^{n + o(n)} \ = \ \{f \in {\cal RM}^{\sf P}: \, f$ can be computed 
    by an ${\cal RM}^{\sf P}$\!-machine whose input-output 

\hspace{1.85in} balance and time-complexity are $n + o(n) \}$.

\medskip

\noindent It should be pointed out that the bound $|x| + o(|x|)$ is only
assumed when $f(x)$ is defined; for $x \not\in {\sf Dom}(f)$, we do not 
assume any time-bound. Of course, there exists also a machine that runs
in polynomial time for all inputs, but then it is not guaranteed that 
the running time is $|x| + o(|x|)$ for accepted inputs.
An ${\cal RM}^{\sf P}$\!-machine whose time and balance on accepted inputs
are $n + o(n)$ is called an {\em ${\cal RM}^{n + o(n)}$-machine}.

Note that ${\cal RM}^{n + o(n)}$ is a strict subset of
${\cal RM}^{\sf lin}$, that consists of the elements of ${\cal RM}^{\sf P}$
that have linear upper-bounds on their balance and their time-complexity
(where by ``linear'' we mean any function of the form $n \mapsto an+b$
for some natural integers $a, b$). Indeed, if $t(.)$ is $n + o(n)$ then
$t(n) \le 2n + c$ for some constant $c$; the strictness of the inclusion
comes from the fact that ${\cal RM}^{\sf lin}$ contains, for example,
functions whose output-length is twice the input-length, and the function 
$n \mapsto 2n$ is not $n + o(n)$.

\begin{lem} \label{RMo_monoid}
 \ \ ${\cal RM}^{n + o(n)}$ is a monoid.
\end{lem}
{\bf Proof.} Let $f_1, f_2 \in {\cal RM}^{n + o(n)}$ and let $M_1, M_2$ be 
${\cal RM}^{n + o(n)}$-machines that compute $f_1$, respectively $f_2$.
Since the set of functions that are $n + o(n)$ is closed under 
composition, the I/O-balance of $f_2 \circ f_1$ is $n + o(n)$.

To compute $f_2 \circ f_1(x)$ in time $n + o(n)$ (where $n = |x|$), we
combine $M_1$ and $M_2$ into an ${\cal RM}^{n + o(n)}$-machine $M$, as 
follows. The output-tape of $M_1$ and the input-tape of $M_2$ are combined 
into one work-tape of $M$; we call this work-tape the intermediate tape. 
On input $x$, the machine $M$ starts simulating $M_1$ and starts writing 
$f_1(x)$ on the intermediate tape; as soon as there is something on this
intermediate tape, $M$ starts the simulation of $M_2$ on $f_1(x)$.
The writing of $f_1(x)$ by $M_1$ takes at most $o(n)$ more steps than it
takes to read $x$; the computation of $f_1(x)$, except for this 
$o(n)$-step delay, is done in parallel (simultaneously) with the reading of 
$x$.  Similarly, when $M_2$ reads $f_1(x)$ as an input, it computes 
$f_2(f_1(x))$ at the same time as it reads $f_1(x)$, except for a
$o(|f_1(x)|)$-step delay; but $o(|f_1(x)|)$ means
$\, \le \epsilon_1(|x|) \cdot |f_1(x)| \, \le \,$
$\epsilon_1(|x|) \cdot (|x| + \epsilon_2(|x|) \cdot |x|)$, and this is
$\, \le \epsilon(|x|) \cdot |x| \, $ (for some functions $\epsilon$
with limit 0); hence, $o(|f_1(x)|)$ is $o(|x|)$.
So when $x \in {\sf Dom}(f_2 \circ f_1)$ the total time taken by $M$
(i.e., $M_1$ and $M_2$ working together, mostly in parallel) is
$|x| + o(|x|)$.
 \ \ \ $\Box$

\begin{pro} \label{RMo_nonreg}
 \ The monoid ${\cal RM}^{n + o(n)}$ is non-regular.
In fact, there exists a real-time function in ${\cal RM}^{n + o(n)}$ that 
has no inverse in ${\cal RM}^{n + o(n)}$. 
\end{pro}
{\bf Proof.} We use the encoding function 
$\, {\sf code}$: $\{0,1, \#\}^* \longmapsto \{00, 01, 11\}^*$,
replacing $0$ by $00$, $1$ by $01$, and $\#$ by $11$ (as discussed at the 
beginning of the Introduction). 
For a string $x$, $x^{\sf rev}$ denotes the string in reverse (i.e., 
backwards) order.  

Consider the right-ideal morphism defined for all $x, w \in \{0,1\}^*$ by  

\smallskip

\hspace{.5in}
$s$: $ \ {\sf code}(x) \ 11 \ 0^{|x|} \, w \ \longmapsto $
$ \ 0^{2 \, |x|} \ 11 \ x^{\sf rev} \, w$ , 

\smallskip

\noindent where $\, {\sf domC}(s)$ $ \ = \ $ 
$\bigcup_{k \ge 0} \, \{00, 01\}^k \, 11 \, 0^k$.
Thus $s$ is injective and length-preserving, and belongs to 
${\cal RM}^{\sf P}$. Moreover, $s$ belongs to 
${\cal RM}^{n + o(n)}$ since an ${\cal RM}^{n + o(n)}$-machine can compute 
$\, s \big({\sf code}(x) \ 11 \ 0^{|x|} \, w \big) \, $ in time 
$\le n + o(n)$, where $n = 2 \, |x| + 2 + |x| + |w|$, as follows:
The machine reads ${\sf code}(x)$ in time $2 \, |x|$, while writing $x$ on a 
work-tape and while writing $0^{2 \, |x|}$ on the output-tape. When $11$ is
encountered in the input, the head of the work-tape is at the right end of 
$x$. The machine copies $11$ to the output-tape, then reads $0^{|x|}$ in the 
input, while copying the work-tape from right to left to the output-tape; 
thus, $x^{\sf rev}$ is written. When the work-tape head reaches the left end 
of the work-tape, the input-tape head reaches $w$ on the input-tape while
copying it to the output-tape. 
Note that the above machine is a real-time Turing machine, with running time 
$\le n + c$ for some constant $c \ge 0$.

We show next that $s$ does not have an inverse in ${\cal RM}^{n + o(n)}$;
hence ${\cal RM}^{n + o(n)}$ is not regular. 
For every inverse $s'$ of $s$ we have 
$s'(0^{2 \, |x|} \ 11 \ x^{\sf rev}) \, = \, {\sf code}(x) \ 11 \ 0^{|x|}$. 
It is easy to see that although $s'$ can be chosen so as to belong to 
${\cal RM}^{\sf P}$, $s'$ cannot be evaluated in time $\le n + o(n)$; here, 
$n = |{\sf code}(x)| + 2 + |x| = 2 \, |x| + 2 + |x|$. 
Indeed, an ${\cal RM}^{n + o(n)}$-machine reads the input
$0^{2 \, |x|} \ 11 \ x^{\sf rev}$ only once, from left to right. While 
$0^{2 \, |x|}$ is being read, $x^{\sf rev}$ has not yet been seen, so no 
letter-pair of ${\sf code}(x)$ can be written on the output-tape; indeed, the 
machine is deterministic, so anything written on the output-tape up to this 
moment would be false for some input $x$. At the moment $x^{\sf rev}$ starts 
being read, the machine has made $2 \, |x| + 2$ steps, and no output has 
been written yet  (except perhaps 
one 0). To write down the output (which has length $n$) will take 
at least $n$ steps from here onward. So the total time will be 
$\ge n + 2 \, |x| + 2 \ge n + n/2$. But $n + n/2$ does not have $n + o(n)$ 
as an upper-bound.
 \ \ \ $\Box$

\begin{pro} \label{RMo_equivPMpoly}
 \ The monoid $\, {\cal RM}^{n + o(n)}\!/\!\equiv_{\sf poly} \, $ is 
isomorphic to $\, {\cal RM}^{\sf P}\!/\!\equiv_{\sf poly}$
($= {\cal M}^{\sf P}_{\sf poly}$).
\end{pro}
{\bf Proof.} We show that the embedding 
$\, [g]_{\sf poly} \in {\cal RM}^{n + o(n)}\!/\!\equiv_{\sf poly}$
$ \ \longmapsto \ $ $[g]_{\sf poly} \in {\cal M}^{\sf P}_{\sf poly} \, $  
is surjective.  More precisely, for every $f \in {\cal RM}^{\sf P}$ with
time-complexity and balance $\le p(.)$ (a polynomial), we construct a
function $F \in {\cal RM}^{n + o(n)}$ such that $F \equiv_{\sf poly} f$.
This is done by a padding argument similar to the one in Lemma
\ref{recTimeConstr}:  $\, F$ is the restriction of
$f$ to $\, \bigcup_{x \in {\sf domC}(f)} \, x \, A^{q(|x|)} \ A^*$,
where $q(.)$ is a fully time-constructible function that satisfies
$(n \cdot p(n))^2 < q(n) \,$ for all $n \in {\mathbb N}$.
The function $q(.)$ exists by Lemma \ref{recTimeConstr}(1); in fact, since
$p(.)$ is a polynomial, $q(.)$ can be chosen to be a fully time-constructible
polynomial of the form $n \mapsto a \, (n + 1)^d$.
Since every $A^{q(|x|)}$ is a maximal prefix code and since $q$ is
polynomial bound, $F \equiv_{\sf poly} f$.

We want to show that $F \in {\cal RM}^{n + o(n)}$. We construct a
Turing machine $M_F$ that on input $x v w$ computes $f(x) \, v w$ in time
$n + o(n)$, where $x \in {\sf domC}(f)$, $v \in A^{q(|x|)}$, $w \in A^*$,
and $n = |xvw|$. On input $z \not\in {\sf Dom}(f)$, $M_F$ should reject,
but we do not care how much time $M_F$ takes in that case.
On input $x v w$, $M_F$ works as follows.
First it finds $x$ as the first prefix of the input that belongs to
${\sf Dom}(f)$, and writes $f(x)$ on the output-tape; this takes time
$\, \le 2 \, |x| \cdot p(|x|) \, $ (see the proof of Lemma
\ref{recTimeConstr}(2)). If no prefix of the input is in ${\sf Dom}(f)$,
$M_F$ rejects.
If $x \in {\sf domC}(f)$ is found, since ${\sf domC}(f)$ is a prefix code,
the end of $x$, and the beginning of $v$, are uniquely determined within the
input $xvw$.
Next, $M_F$ finds $v$, consisting of the next $q(|x|)$ letters, and
concatenates this to the right of $f(x)$ on the output-tape;
since $q$ is fully time-constructible, this can be done in time $q(|x|)$ 
exactly. If the remainder $vw$ of the input has length $< q(|x|)$, $M_F$ 
rejects.
If $v$ is found, the remainder $w$ of the input is copied to the output-tape.

So the total time of the computation (if there is an output) is
$ \, \le \, 2 \, |x| \cdot p(|x|) + q(|x|) + |w|$.  Since $n = |xvw|$,
$|v| = q(|x|)$, and since $\, 2 \, |x| \cdot p(|x|) = \sqrt{q(|x|)}$, the 
total time is $ \, \le \sqrt{n} + n$, which is $\, \le n + o(n)$.  
 \ \ \ $\Box$

\medskip

\noindent As a consequence of Prop.\ \ref{RMo_equivPMpoly} and earlier
results we have:

\begin{cor} \label{WitnessPNP}
 \ \ ${\sf P} \neq {\sf NP}$ iff there exists a function in
${\cal RM}^{n + o(n)}$ that has no inverse in ${\cal RM}^{\sf P}$.
\end{cor}
{\bf Proof.} If some $F \in {\cal RM}^{n + o(n)}$ ($\subset$
${\cal RM}^{\sf P}$) has no inverse in ${\cal RM}^{\sf P}$, then
${\cal RM}^{\sf P}$ is not regular, hence ${\sf P} \neq {\sf NP}$ (by
results from \cite{s1f}, as we saw in the Introduction).

Conversely, suppose every $F \in {\cal RM}^{n + o(n)}$ has some inverse
$F' \in {\cal RM}^{\sf P}$. By Prop.\ \ref{RMo_equivPMpoly}, every element
of ${\cal M}^{\sf P}_{\sf poly}$ is an $\equiv_{\sf poly}$-class of the
form $[F]$ for some $F \in {\cal RM}^{n + o(n)}$.
Since $F F' F = F$, we have $[F] [F'] [F] = [F]$ (since $\equiv_{\sf poly}$
is a congruence by Theorem \ref{polyequivCongr}).
Hence ${\cal M}^{\sf P}_{\sf poly}$ is regular. By Theorem \ref{regofRMequ}
this implies that ${\cal RM}^{\sf P}$ is regular.
 \ \ \ $\Box$

\bigskip

\noindent {\bf Remark.} By Corollary \ref{WitnessPNP}, if
${\sf P} \neq {\sf NP}$ then this is ``witnessed'' by an element of
${\cal RM}^{n + o(n)}$.  Although ${\cal RM}^{n + o(n)}$ is not regular by
itself, its non-regularity in ${\cal RM}^{\sf P}$ is not obvious (and
equivalent to ${\sf P} \neq {\sf NP}$). 

It is not especially surprising that 
${\cal RM}^{n + o(n)}$ is non-regular; ultimately, this is due to the
limitations of tapes as storage devices. By itself, it is not too surprising
either that ${\cal RM}^{n + o(n)}$ is $\, \equiv_{\sf poly}$-equivalent to 
${\cal M}^{\sf P}_{\sf poly}$.
The $\, \equiv_{\sf poly}$-equivalence of ${\cal RM}^{n + o(n)}$ and 
${\cal M}^{\sf P}_{\sf poly}$ is proved by pushing the familiar padding 
argument a little further. The combination of the two facts is interesting, 
however, because $\, \equiv_{\sf poly}$-equivalence means that 
${\cal RM}^{n + o(n)}$ and ${\cal M}^{\sf P}_{\sf poly}$ are very close to 
each other; yet, ${\cal RM}^{n + o(n)}$ is non-regular, while the 
non-regularity of ${\cal M}^{\sf P}_{\sf poly}$ is equivalent to 
${\sf P} \neq {\sf NP}$.

\bigskip

\bigskip

\noindent {\bf Acknowledgement.} The author thanks the referee for a very
thorough and thoughtful reading of the journal version of this paper, which
contributed to many  corrections and improvements.


\bigskip



\end{document}